\documentclass[12pt]{elsart}

\usepackage{graphicx}% Include figure files

\usepackage{epstopdf}% To incorporate .eps illustrations using PDFLaTeX, etc.%¶à¼ÓÈëµÄ°ü
\usepackage{enumerate}

\usepackage{amsmath}
 \usepackage{color}
\usepackage{changepage}
\usepackage{subfigure}
\usepackage{multirow}

\usepackage{amssymb}
\usepackage [latin1]{inputenc}
\usepackage[compress]{cite}

\theoremstyle{plain}% Theorem-like structures provided by amsthm.sty
\newtheorem{theorem}{Theorem}[section]
\newtheorem{lemma}[theorem]{Lemma}

\theoremstyle{definition}
\newtheorem{definition}[theorem]{Definition}

\theoremstyle{remark}
\newtheorem{remark}{Remark}

\begin{document}

\begin{frontmatter}

\title{The normalized Laplacian spectrum of  $n$-polygon graphs and its applications}

%\author{
%\name{Tengjie Chen\textsuperscript{a, b}, Zhenhua Yuan\textsuperscript{a, b} and Junhao Peng\textsuperscript{a, b}\thanks{Junhao Peng Email: pengjh@gzhu.edu.cn}}
%\affil{\textsuperscript{a}School of Mathematics and Information Science, Guangzhou University, Guangzhou 510006, P.R. China;\\
% \textsuperscript{b}Guangdong Provincial Key Laboratory co-sponsored by province and city of Information Security Technology, Guangzhou University, Guangzhou 510006, P.R. China}
%}

\author[lable1,label2]{Tengjie Chen}
\author[lable1,label2]{Zhenhua Yuan}
\author[lable1,label2]{Junhao Peng}
\ead{pengjh@gzhu.edu.cn}

%\ead{}

\address[lable1]{School of Math and Information Science, Guangzhou University, Guangzhou 510006, China.}
\address[label2]{Guangdong Provincial Key Laboratory co-sponsored by province and city of Information Security Technology, Guangzhou University, Guangzhou 510006,  China.}
%\address[label3]{Research Center for Computer Science and Information Technologies, Macedonian Academy of Sciences and Arts, Bul. Krste Misirkov 2, 1000 Skopje, Macedonia.}
%\address[label4]{Institute of Physics \& Astronomy, University of Potsdam, D-14776 Potsdam-Golm, Germany.}
%\address[label5]{Institute of Physics, Faculty of Natural Sciences and Mathematics, Ss.~Cyril and Methodius University, Arhimedova 3, 1000 Skopje, Macedonia.}
%\address[label6]{Faculty of Computer Science and Engineering, Ss. Cyril and Methodius University, \\ P.O. Box 393, 1000 Skopje, Macedonia}
\begin{abstract}
Given an arbitrary connected  $G$, the $n$-polygon graph $\tau_n(G)$ is obtained by adding a path with length $n$ $(n\geq 2)$ to each  edge of graph $G$,  and the iterated $n$-polygon graphs $\tau_n^g(G)$ ($g\geq 0$), is obtained from the iteration $\tau_n^g(G)=\tau_n(\tau_n^{g-1}(G))$, with initial condition $\tau_n^0(G)=G$.
In this paper, a method for calculating the eigenvalues of normalized Laplacian matrix for graph  $\tau_n(G)$ is presented if the eigenvalues of normalized Laplacian matrix for graph  $G$ is given firstly. Then, the  normalized Laplacian spectrums for the graph $\tau_n(G)$ and the graphs $\tau_n^g(G)$ ($g\geq 0$) can also be derived. Finally, as  applications, we calculate the multiplicative degree-Kirchhoff index, Kemeny's constant and the number of spanning trees for the  graph $\tau_n(G)$ and the graphs $\tau_n^g(G)$ by exploring their connections with the  normalized Laplacian spectrum, exact results for these quantities are obtained.

\end{abstract}
\end{frontmatter}
%% main text

\section{Introduction}
%\label{intro}
The spectra of the adjacency matrices, Laplacian and the normalized Laplacian for a graph is of particular important since many important structural and dynamical properties of the graph can be obtained from the eigenvalues and eigenvectors of these matrices~\cite{LO93, cvetkovic2010-book}. For example, the spectra of these matrices  provides information on  degree distribution, the multiplicative degree-Kirchhoff index, community structure, local clustering, total number of links, the number of spanning trees and etc.~\cite{cvetkovic2010-book, chung1997book,  MehatariB-2015-AMC}. Also, first-passage properties of the networks~\cite{Book-Redner-2007, WenHW-2019-IJB}, such as mean first-passage time, the mean hitting time, the  mixing time and the Kemeny's const, can be expressed in terms of the eigenvalues of the Laplacian and normalized Laplacian~\cite{Wang2004-BOOK, cvetkovic2010-book, Brouwer2012-BOOK}. In the past several years, there is a particular interest in the study of spectra for  the normalized Laplacian matrices  of different graphs~\cite{CaversFK-2010LMA, Lin-Zhang-2013-PRE, Zhang-Lin-2015-PRE, huang_li_2015-BAMS, Zhang-Lin-Guo-2015PRE, SHARMA2017-DAM, QiZhang-2019-IEEE, HeLi2018-JCAM, CheonKM2018-LMA}.

Recently, a class of iterated graphs obtained by replacing each edge of a `base graph' with a given structure, have attracted lots  attention~\cite{Dorogovtsev-Goltsev-2002-PRE, Rozenfeld-2006-NJP}, Since these graphs not only exhibit rich structural properties, such as self-similarity, fractal~\cite{Rozenfeld-2006-NJP, Gao_2021-JSTAT, Peng-Agliari-2017-Chaos}, scale-free and small-world properties~\cite{Dorogovtsev-Goltsev-2002-PRE, Peng-Agliari-Zhang-2015-Chaos}, but also show distinctive dynamic properties~\cite{Gao_2021-JSTAT, Zhang-2011-JPA, Peng-Agliari-Zhang-2015-Chaos, Agliari2008-PRE, Peng-Agliari-2017-Chaos}. %constructed recursively

For an arbitrary connected `base graph' $G$, if the normalized Laplacian spectrum of graph $G$ is known, how can we calculate the normalized Laplacian spectrum  and related quantities of graph $\tau(G)$, which is obtained from $G$ by replacing each edge of $G$  with a given structure. It is an interesting topic with wide applications and lots of results were obtained. Examples include the graph $\tau(G)$ obtained by replacing every edges of $G$  with a geometry, such as  triangles~\cite{XIE2016-AMC, Chen2018LMA, WANG2018-AMC}, a quadrilateral~\cite{GuoL2021-LMA, LI2017-AMC, HuangLi2018LMA}, a pentagon~\cite{XuWW-2019JMRA}. However, for more general  graphs $\tau_n(G)$  obtained by replacing every edges of $G$  by a $(n+1)$-polygon ($n\geq 2$), related topic is more difficult and general result for the normalized Laplacian spectrum  $\tau_n(G)$ is still unknown.

To fill this gap, we consider the case while every edges of $G$ are replaced by a $(n+1)$-polygon. We analyze the relation between  the normalized Laplacian spectrum of the graph $\tau_n(G)$  and that of graph $G$ and present detailed  spectrum of the graph $\tau_n(G)$ for any $n\geq 2$. Using our results recursively, we also obtain the normalized Laplacian spectrum of the iterated  graph $\tau_n^g(G)$ ($g\geq 0$), where $\tau_n^0(G)=G$, and   $\tau_n^g(G)=\tau_n(\tau_n^{g-1}(G))$. As  applications, we also analyze the multiplicative degree-Kirchhoff index, the Kemeny's const and the number of spanning trees of the graph $\tau_n^g(G)$. Exactly results for these quantities are also presented.

This paper is organized as follows. First, in Sec.~\ref{Def}, we present related definitions and notations which will appear in the manuscript. In Sec.~\ref{pre}, we present some
classical results we quote and some preliminary results which will be used in the manuscript.  Then, in Sec.~\ref{Main_result}, we present the main results on the normalized Laplacian spectrum of the graph $\tau_n(G)$,  and in Sec.~\ref{Sec:Appli}, we calculate the multiplicative degree-Kirchhoff index, the Kemeny's const and the number of spanning trees of the graph $\tau_n^t(G)$ and $\tau_n^g(G)$. The detailed proof of our results are presented in Sec.\ref{Proof_Lemma121}-\ref{Proof-t34} and Appendix~\ref{Proof_Pro_Rec_an}-\ref{Prov_Pro_an2}.

\section{Definitions and notations}
\label{Def}

 Let $G=(V(G),E(G))$ be a undirected graph with vertex set $V(G)$ and edge set $E(G)$, $N_0=|V(G)|$ represent the total number of  nodes of $G$ and $E_0=|E(G)|$ be the total number of edges of graph $G$.  For any two nodes $i,j\in V(G)$, if there is an edge between node $i$ and $j$ in $E(G)$, we say $i$ is a neighbor of $j$ or  $i$ and $j$ are adjacent (i.e., $i\sim j$). If $e$ is an edge with end-vertices $i$ and $j$,  we say  edge $e$ is incident to nodes $i$  and $j$, or nodes $i$  and $j$ are incident with edge $e$. The degree of vertex $i$, referred to as $d_{i}$, is the number of edges incident to node $i$. %connected
 \begin{definition}\label{defineLA} ~\cite{cvetkovic2010-book}
  For any undirected graph $G$, the Laplacian matrix of $G$ is defined as
 \begin{displaymath}
 L(G) = D(G)- A(G),
\end{displaymath}
where $D(G) = diag(d_1, d_2, \cdots , d_{N_0})$ is the degree diagonal matrix of $G$ and $A(G)=(A_{ij})_{N_{0}\times E_{0} }$ is the adjacency matrix of $G$, with
 \begin{equation} \label{ADJ}
  A_{ij}=\left\{ \begin{array}{ll} 1 & i\sim j\\ 0 & \text{others} \end{array} \right..
\end{equation}
 \end{definition}

  \begin{definition}\label{defineNorL} ~\cite{cvetkovic2010-book, Zhang-2020-book}
  For any undirected graph $G$, the normalized Laplacian of $G$ is defined to be
     \begin{equation} \label{Def_NLaplacian}
     \mathcal{L}_{G}=D(G)^{-\frac{1}{2}}L(G)D(G)^{-\frac{1}{2}}=I-D(G)^{-\frac{1}{2}}A(G)D(G)^{-\frac{1}{2}}.
     \end{equation}
 \end{definition}

     Let $\delta_{ij}$ be the Kronecker delta function. Then, the $(i, j)$-th entry of the matrix $\mathcal{L}_{G}$ can be written as
     \begin{equation}\label{Lij}
      \mathcal{L}_{G}\left ( i,j \right )=\delta_{ij}-\frac{A_{ij}}{\sqrt{d_{i}d_{j}}}.
     \end{equation}
     %where $\mathcal{L}_{G}\left ( i,j \right )$ and $A\left ( i,j \right )$ denote the (i,j)-entry of $\mathcal{L}_{G}$ and $A$ respectively.
     %Since $\mathcal{L}_{G}$ is Hermitian and similar to $I-M=D^{-1}L$,
     \begin{definition}\label{defineNorLS} ~\cite{cvetkovic2010-book}
  For any graph undirected $G$ with $N_0$ vertexes, the normalized Laplacian spectrum of  graph $G$ is defined to be
     \begin{equation} \label{Def_NLaplacianS}
     \sigma(G) =\left \{ \lambda_{1}, \lambda_{2},\cdots ,\lambda_{N_{0}} \right \},
     \end{equation}
     where $\lambda_{i}$ ($i=1,~2,~\cdots,~N_0$) are the eigenvalues of $\mathcal{L}_{G}$.
 \end{definition}

   \begin{definition}\label{defineB} ~\cite{chung1997book}
 Let $G$ be a  undirected  graph, with vertex set $V\left ( G \right )=\left \{ 1,2,\cdots ,N_{0} \right \}$ and edge set $E\left ( G \right )=\left \{ e_{1},e_{2},\cdots ,e_{E_{0}} \right \}$. Then the incidence matrix of graph $G$ is defined by
\begin{displaymath}
B=\left ( b_{ij} \right )_{N_{0}\times E_{0} },
\end{displaymath}
 where
\begin{displaymath}
b_{ij}= \left\{ \begin{array}{ll}
1, & \textrm{if vertex $i$ and edge $e_{j}$ are incident}\\
0, & \textrm{otherwise}
\end{array} \right..
\end{displaymath}
 \end{definition}

\begin{definition}\label{define01}~\cite{Wang2004-BOOK, chung1997book}
Let  $G^{'}$ be a directed graph, $V\left ( G^{'} \right )=\left \{ 1,2,\cdots ,N_{0} \right \}$ and directed edge set $E\left ( G^{'} \right )=\left \{e_{1}^{'},e_{2}^{'},\cdots ,e_{E_{0}}^{'} \right \}$. Then the incidence matrix of directed graph $G^{'}$ is defined by
\begin{displaymath}
B^{'}=\left ( b_{ij}^{'} \right )_{N_{0}\times E_{0} },
\end{displaymath}
 where
\begin{displaymath}
 b_{ij}^{'}= \left\{ \begin{array}{ll}
-1, & \textrm{if vertex $i$ is the source vertex of edge $e_{j}^{'}$}\\
0, & \textrm{if vertex $i$ is not incident on edge $e_{j}^{'}$ }\\
1, & \textrm{if vertex $i$ is the target vertex of edge $e_{j}^{'}$ }
\end{array} \right..
\end{displaymath}
\end{definition}

\begin{definition}\label{define02}~\cite{Wang2004-BOOK, chung1997book}
Let  $G^{'}$ be a directed graph, $G^{'}$ is called as weakly connected if its underlying undirected graph (i.e. a graph obtained by  replacing each directed edge of $G^{'}$ with an undirected edge) is connected. %That is to say, for any two nodes $i$ and $j$ in $V(G^{'})$, there is an undirected path in $G^{'}$ from $i$ to $j$.
\end{definition}

\begin{figure}
\begin{center}
\includegraphics[scale=1.05]{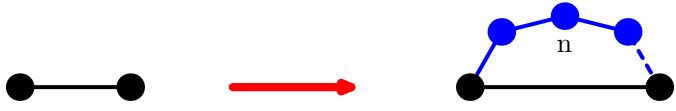}%{huatu-3.eps}
\caption{The construction method of the $n$-polygon graph. The $n$-polygon graph of $G$, denoted by $\tau_n(G)$,  is obtained from $G$ by replacing each edge $ij$ of $G$  with the $(n+1)$-polygon on the right-hand side of the arrow.}
\label{fig:1}       % Give a unique label
\end{center}
\end{figure}

\begin{definition}\label{define1}
The $n$-polygon $(n\geq 2)$ graph of $G$, denoted by $\tau_n(G)$, is the graph obtained by adding a path with length $n$  to each edge  of $G$.
\end{definition}
In other word, $\tau_n$ is a operator which turns each edge of a graph into a $(n+1)$-polygon. Thus, $\tau_n(G)$ is a graph obtained by replacing each edge $ij$ of $G$ with a $(n+1)$-polygon, as shown in Fig.~\ref{fig:1}.
 \begin{definition}\label{define2}
For any $g>0$, the iterated $n$-polygon graph with generation $g$ $(g\geq 1)$, referred to as  $\tau_n^g(G)$, is defined as the graph obtained through the iteration $\tau_n^g(G)=\tau_n(\tau_n^{g-1}(G))$,
with initial condition $\tau^0(G)=G$.
\end{definition}

%\begin{figure}
%\centering
%\includegraphics[width=\linewidth]{huatu-3}
%\caption{The constructing }
%\end{figure}

%%%%%%%%%%%%%%%%%%%%%%%%%%%%%%%%%%%%%%%%%%%%%%%%%%%%%%%%%
% Figure  1
%%%%%%%%%%%%%%%%%%%%%%%%%%%%%%%%%%%%%%%%%%%%%%%%%%%%%%%%%%

%\begin{figure}
%\centering
%\subfloat[An example of an individual figure sub-caption.]{%
%\resizebox*{5cm}{!}{\includegraphics{graph1.eps}}}\hspace{5pt}
%\subfloat[A slightly shorter sub-caption.]{%
%\resizebox*{5cm}{!}{\includegraphics{graph2.eps}}}
%\caption{Example of a two-part figure with individual sub-captions
% showing that captions are flush left and justified if greater
% than one line of text.} \label{sample-figure}
%\end{figure}

 Let $N_g$ and $E_g$  be the total number of vertices and edges of graph $\tau_n^g(G)$ respectively.

For any $g>0$, we have
\begin{equation}\label{ng}
N_{g}=N_{g-1}+(n-1)E_{g-1},~E_{g}=(n+1)E_{g-1}\nonumber.
\end{equation}
%$N_{g}=N_{g-1}+(n-1)E_{g-1}, E_{g}=(n+1)E_{g-1}$\\
Therefore
\begin{equation}\label{ng1}
N_{g}=N_{0}+(n-1)\frac{(n+1)^{g}-1 }{n}E_{0},
\end{equation}
and
\begin{equation}\label{Eg}
 E_{g}=(n+1)^{g}E_{0}.
\end{equation}

\section{Preliminaries}
\label{pre}
 In this section, we present some classical results we quote and some preliminary results which will be used to derive our main results in the manuscript.

 \begin{lemma}\label{L1}~\cite{chung1997book}
	For any connected undirected graph $G$, let $\lambda_{1}\leq \lambda_{2}\leq \cdots \leq \lambda_{N_{0}}$ be the eigenvalues of $\mathcal{L}_{G}$, Then
 %is the normalized Laplacian of $G$ and $0=\lambda_{1}<  \lambda_{2} \leq \cdots \leq \lambda _{N_{0}}$ are the  eigenvalues of $\mathcal{L}_{G}$.
%The normalized Laplacian spectrum of $G$ is $\sigma _{0}=\left \{ 0=\lambda_{1}<  \lambda_{2}\leq \cdots \leq \lambda _{N_{0}} \right \}$. We have\\
\begin{enumerate}[(i)]
  \item $\lambda_{1}=0$, and  for any $2\leq i\leq N_{0}$,  $0< \lambda_{i}\leq 2$;
  \item  $\lambda_{N_{0}}= 2$ if and only if $G$ is bipartite;
  \item $G$ is bipartite if and only if both $\lambda_{i}$ and $2-\lambda_{i}$ are eigenvalues of $\mathcal{L}_{G}$
  %$\mu_{i}$ is an eigenvalue of $\mathcal{L}_{G}$,then the value $2-\mu_{i}$ is also an eigenvalue of $\mathcal{L}_{G}$
   and $m_{\mathcal{L}_{G}}\left ( \lambda_{i} \right )=m_{\mathcal{L}_{G}}\left ( 2-\lambda_{i} \right )$, where $m_{\mathcal{L}_{G}}\left ( \lambda_{i} \right )$ is the multiplicity of $\lambda_{i}$ which is the eigenvalue of $\mathcal{L}_{G}$. % of $\mathcal{L}_{G}$.
\end{enumerate}
	\end{lemma}
% Therefore we can label the eigenvalues of  $\mathcal{L}_{G}$ as $0=\lambda_{1}< \lambda_{2}\leq \cdots \leq \lambda_{N_{0}}$, $\mathcal{L}_{G}$ has a unique $0$ eigenvalue and
 Therefore, all the eigenvalues of $\mathcal{L}_{G}$ are non-negative. %and the normalized Laplacian spectrum  of the graph  $G$ can be labeled as $\sigma _{G} =\left \{ \lambda_{1}, \lambda_{2},\cdots ,\lambda_{N_{0}} \right \}$. %,which $\sigma _{0}$ is also called the normalized Laplacian spectrum of $G$.
	Given the normalized Laplacian spectrum of  graph $G$,  the multiplicative degree-Kirchhoff index, the Kemeny's constant and the number of spanning trees of graph $G$ can be expressed as follows. %$\sigma _{0} =\left \{ 0=\lambda_{1}< \lambda_{2}\leq \cdots \leq \lambda_{N_{0}} \right \}$ is the spectrum on the normalized Laplacian
\begin{lemma}\label{L2}~\cite{Wang2004-BOOK, chung1997book, cvetkovic2010-book}
Let $G$ be a connected  undirected graph with $N_{0}$ vertices and $E_{0}$ edges, $0=\lambda_{1}< \lambda_{2}\leq \cdots \leq \lambda_{N_{0}}$ are the  eigenvalues of the normalized Laplacian $\mathcal{L}_{G}$. Then,
\begin{enumerate}[(i)]
 \item the multiplicative degree-Kirchhoff index of $G$ can be written as
\begin{displaymath}
Kf^{'}\left ( G \right )=2E_{0}\sum_{i=2}^{N_{0}}\frac{1}{\lambda_{i}};
\end{displaymath}
 \item %\cite{Butler2016}
the Kemeny's constant of $G$ is
\begin{displaymath}
K\left ( G \right )=\sum_{i=2}^{N_{0}}\frac{1}{\lambda_{i}};
\end{displaymath}
 \item %\cite{chung1997spectral}
the number $N_{st}\left ( G \right )$ of spanning trees of $G$ is
\begin{displaymath}
N_{st}\left ( G \right )=\frac{1}{2E_{0}}\left(\prod_{i=1}^{N_{0}}d_{i}\right)\prod_{i=2}^{N_{0}}\lambda_{i}.
\end{displaymath}
\end{enumerate}
	\end{lemma}

%\begin{lemma}\label{Vieta}~\cite{Borwein1997-book}
%Let $f(x)=b_{n}x^{n}+b_{n-1}x^{n-1}+\cdots +b_{1}x+b_{0}$ be a polynomial of degree $n$ and  $r_{i}$ $(i=1,~2,~\cdots,~n)$ be the roots of polynomial equation $f(x)=0$.
%Then Vieta's formulas states that\\
%(i)
%\begin{displaymath}
%\prod_{i=1}^{n}r_{i}=\left ( -1 \right )^{n}\frac{b_{0}}{b_{n}},
%\end{displaymath}
%%\begin{displaymath}
%%\sum_{i=1}^{n}{r_{i}}=-\frac{b_{n-1}}{b_{n}},
%%\end{displaymath}
%(ii)
%\begin{displaymath}
%\left(\prod_{j=1}^{n}r_{j}\right)\sum_{i=1}^{n}\frac{1}{r_{i}}=(-1)^{n-1}\frac{b_{1}}{b_{n}}.
%\end{displaymath}
%\end{lemma}

\begin{lemma}\label{Vieta}~\cite{Borwein1997-book}%For a polynomial of the form with roots $r_{i} ,i=1,2,\cdots ,n$,
Let $f(x)=b_{n}x^{n}+b_{n-1}x^{n-1}+\cdots +b_{1}x+b_{0}$ $(b_{n}\neq0)$ be a polynomial of degree $n$ and  $r_{i}$ $(i=1,~2,~\cdots,~n)$ be the roots of $f(x)=0$.
Then Vieta's formulas states that\\
%\begin{equation} \label{D0h}
%  D_0=\left\{ \begin{array}{ll}\frac{n(n^2-1)}{8} & \text{if $n$  is odd}\\ \frac{n^3}{8} & \text{if $n$  is even} \end{array} \right..
%\end{equation}
\begin{displaymath}
\prod_{i=1}^{n}r_{i}=\left ( -1 \right )^{n}\frac{b_{0}}{b_{n}},
\end{displaymath}
and %in the case $r_{i}\neq 0$ $(i=1,~2,~\cdots,~n)$,
\begin{displaymath}
\sum_{i=1}^{n}\frac{1}{r_{i}}=-\frac{b_{1}}{b_{0}}.
\end{displaymath}
\end{lemma}

%\begin{lemma}\label{L3}%\cite{cvetkovic2010introduction}
%Let $B$ be the incidence matrix of a connected graph $G$ with $N_{0}$ vertices.Then
%\begin{displaymath}
%rank\left ( B \right )= \left\{ \begin{array}{ll}
%N_{0}-1& \textrm{if $G$ is bipartite,}\\
%N_{0} & \textrm{if $G$ is non-bipartite. }
%\end{array} \right.
%\end{displaymath}
%\end{lemma}

\begin{lemma}\label{L3}\cite{cvetkovic2010-book}
Let $B$ be the incidence matrix of a connected graph $G$ with $N_{0}$ vertices. Then
\begin{displaymath}
rank\left ( B \right )= \left\{ \begin{array}{ll}
N_{0}-1& \textrm{if $G$ is bipartite,}\\
N_{0} & \textrm{if $G$ is non-bipartite. }
\end{array} \right.
\end{displaymath}
\end{lemma}

\begin{lemma}\label{L03}~\cite{Wang2004-BOOK, chung1997book}
Let ${B}'$ be the incidence matrix of weakly connected directed graph ${G}'$ with $N_{0}$ vertices.Then $rank\left ( {B}' \right )=N_{0}-1$.
\end{lemma}

\begin{lemma}\label{Pro_Rec_an}
Let $u$ be an arbitrary real number, and $\{a_n(\mu)\}_{n\geq -1}$ is a series which defined by the recursive equation $a_{n}(\mu)=2(1-\mu) a_{n-1}(\mu)-a_{n-2}(\mu)$, with initial conditions $a_{-1}(\mu)=0$ and $a_{0}(\mu)=1$. Then %$a_{1}(\mu)=2(1-\mu)$, $a_{2}=\beta ^{2}-1$, specially\\
%\iffalse
%       \begin{equation}\label{a1}
%       a_{n}=(\frac{1}{2}-\frac{i\beta }{2\sqrt{4-\beta ^{2}}})\left (\frac{\beta +i\sqrt{4-\beta ^{2}}}{2}\right )^{n}+(\frac{1}{2}+\frac{i\beta }{2\sqrt{4-\beta ^{2}}})\left (\frac{\beta -i\sqrt{4-\beta ^{2}}}{2}\right )^{n}.
%       \end{equation}
%Let $\beta =2\left ( 1-\mu \right ),\mu\in \left ( 0,2 \right )$,we have\\
%       \begin{equation}\label{a2}
%       a_{n}=\left [ \frac{1}{2}-\frac{i\left ( 1-\mu \right )}{\sqrt{2\mu-\mu^{2}}} \right ] \left (1-\mu+i\sqrt{2\mu-\mu^{2}} \right )^{n}+\left [ \frac{1}{2}+\frac{i\left ( 1-\mu \right )}{\sqrt{2\mu-\mu^{2}}} \right ] \left (1-\mu-i\sqrt{2\mu-\mu^{2}}\right )^{n}.
%        \end{equation}
%\fi
%By the recursive of $a_{n}$ and $1+a_{n}$,we get\\
\begin{enumerate}[(i)]
  %\item  for any  $n\geq 0$, $a_n(\mu)$  is an $n$th order polynomial in $\mu$;   %a  polynomial in $\mu$ of degree $n$
  \item for any $n\geq 1$,  $a_n(2-\mu)=(-1)^na_n(\mu)$;
  \item if $\mu=0$
  \begin{equation}\label{mu=0}
  a_{n}(0)=n+1,
  \end{equation}
  and if $\mu=2$,
  \begin{equation}\label{mu=2}
  a_{n}(2)=\left ( -1 \right )^{n}\left ( n+1 \right );
  \end{equation}
  \item  for any  $n\geq 0$, $a_n(\mu)$  is an $n$th order polynomial in $\mu$, which can be written as
  $$a_{n}(\mu)=\sum_{i=0}^na_{n}^{(i)}\mu^i,$$
  %define $a_{n}^{\left ( 0 \right )}$,$a_{n}^{\left ( 1 \right )}$ and $a_{n}^{\left ( n \right )}$ as the constant term, coefficients of first order and coefficients of $n$th order, respectively of $a_{n}$. Then
 where $a_{n}^{(i)}$ is the coefficient of $\mu^i$; further more, %we find
  \begin{equation}\label{n=0}
  a_{n}^{\left ( 0 \right )}=n+1,
  \end{equation}
  \begin{equation}\label{n=1}
   a_{n}^{\left ( 1 \right )}=-\frac{n^{3}+3n^{2}+2n}{3},
  \end{equation}
  \begin{equation}\label{n=n}
  a_{n}^{\left ( n \right )}=\left ( -1 \right )^{n}2^{n};
  \end{equation}
  \item for any $n\geq 2$, $a_n(\mu)$ and $1+a_{n}(\mu)$ can be expanded as
  \begin{equation}\label{a3}
  a_{n}(\mu)= \left\{ \begin{array}{ll}
  \left ( a_{\frac{n-1}{2}+1}(\mu)-a_{\frac{n-1}{2}-1}(\mu) \right )a_{\frac{n-1}{2}}(\mu) & \textrm{$n$ is  odd }\\
  \left ( a_{\frac{n}{2}}(\mu)-a_{\frac{n}{2}-1}(\mu) \right )\left ( a_{\frac{n}{2}}(\mu)+a_{\frac{n}{2}-1}(\mu) \right ) & \textrm{$n$ is  even }
  \end{array} \right.,
  \end{equation}
  and
  \begin{equation}\label{a4}
  1+a_{n}(\mu)= \left\{ \begin{array}{ll}
  \left ( a_{\frac{n+1}{2}-1}(\mu)-a_{\frac{n+1}{2}-2}(\mu) \right )\left ( a_{\frac{n+1}{2}}(\mu)+a_{\frac{n+1}{2}-1}(\mu) \right )& \textrm{$n$ is odd}\\
  \left ( a_{\frac{n}{2}}(\mu)-a_{\frac{n}{2}-2}(\mu) \right )a_{\frac{n}{2}}(\mu) & \textrm{$n$ is even }
  \end{array} \right..
  \end{equation}
  \end{enumerate}
\end{lemma}
%Note: in this manuscript, we often use `$a_{k}$' to represent `$a_{k}(\mu)$' to lighten the notations.
The proof of  Lemma~\ref{Pro_Rec_an} is presented in Appendix \ref{Proof_Pro_Rec_an}.
%drop the dependence on $\mu$ to lighten the notations of  $a_{n}(\mu)$.
%  \begin{equation}\label{a3}
%a_{n}(\mu)= \left\{ \begin{array}{ll}
%\left ( a_{\frac{n-1}{2}+1}-a_{\frac{n-1}{2}-1} \right )a_{\frac{n-1}{2}} & \textrm{$n\geq 3$ and $n$ is an odd number }\\
%\left ( a_{\frac{n}{2}}-a_{\frac{n}{2}-1} \right )\left ( a_{\frac{n}{2}}+a_{\frac{n}{2}-1} \right ) & \textrm{$n\geq 2$ and $n$ is an even number}
%\end{array} \right.,
%\end{equation}
%and
%\begin{equation}\label{a4}
%1+a_{n}(\mu)= \left\{ \begin{array}{ll}
%\left ( a_{\frac{n+1}{2}-1}-a_{\frac{n+1}{2}-2} \right )\left ( a_{\frac{n+1}{2}}+a_{\frac{n+1}{2}-1} \right )& \textrm{$n\geq 3$ and $n$ is an odd number}\\
%\left ( a_{\frac{n}{2}}-a_{\frac{n}{2}-2} \right )a_{\frac{n}{2}} & \textrm{$n\geq 2$ and $n$ is an even number }
%\end{array} \right..
%\end{equation}
%{\color{red}
%
%
%}
\begin{lemma} \label{Pro_an1}  %\label{lambda_0_2_1}
\begin{enumerate}[(i)]
If $n$ $(n\geq 3)$ is an odd number, and $\{a_k(\mu)\}_{k\geq 0}$ is a series defined in Lemma \ref{Pro_Rec_an}, we have the following results.
  \item Let $\mu$ be an arbitrary  root of equation $a_{\frac{n-1}{2}}(\mu)-a_{\frac{n-1}{2}-1}(\mu)=0$. Then
  \begin{equation}
1-\frac{ a_{\frac{n-1}{2}+1}\left (\mu \right )-a_{\frac{n-1}{2}-1}\left ( \mu\right ) }{ a_{\frac{n-1}{2}}\left (\mu \right )-a_{\frac{n-1}{2}-2}\left ( \mu \right ) }=2,
\end{equation}
\begin{equation}
a_{n-2}\left ( \mu\right )=1.
\end{equation}
  \item Let $\mu$ be an arbitrary  root of equation $a_{\frac{n-1}{2}}(\mu)+a_{\frac{n-1}{2}-1}(\mu)=0$. Then
  \begin{equation}
1-\frac{ a_{\frac{n-1}{2}+1}\left (\mu \right )-a_{\frac{n-1}{2}-1}\left ( \mu\right ) }{ a_{\frac{n-1}{2}}\left (\mu \right )-a_{\frac{n-1}{2}-2}\left ( \mu \right ) }=0,
\end{equation}
\begin{equation}\label{AN_2C1}
a_{n-2}\left ( \mu\right )=-1.
\end{equation}
  \item  Let $\mu$ be an arbitrary  root of equation $a_{\frac{n-1}{2}}(\mu)=0$. Then
  \begin{equation}\label{ANC1}
1-\frac{ a_{\frac{n-1}{2}+1}\left (\mu \right )-a_{\frac{n-1}{2}-1}\left ( \mu\right ) }{ a_{\frac{n-1}{2}}\left (\mu \right )-a_{\frac{n-1}{2}-2}\left ( \mu \right ) }\in \left ( -\infty ,0 \right )\cup \left ( 2,+\infty  \right ),
\end{equation}
\begin{equation}\label{ANC2}
a_{n-1}\left ( \mu\right )=-1 ~~\textrm{and~~} a_{n-2}\left ( \mu  \right )=-2\left ( 1- \mu  \right ).
\end{equation}
  \end{enumerate}
\end{lemma}
The proof of the Lemma is presented in Appendix \ref{Prov_Pro_an1}.

\begin{lemma} \label{lambda_0_2_2}
\begin{enumerate}[(i)]
If $n$ $(n\geq 2)$  is an even number, and  $\{a_k(\mu)\}_{k\geq 0}$ is a series defined in Lemma \ref{Pro_Rec_an}, we have the following results.
  \item Let $\mu$ be an arbitrary  root of equation $a_{\frac{n}{2}}(\mu)-a_{\frac{n-2}{2}}(\mu)=0$. Then
\begin{equation}\label{even1}
1-\frac{ a_{\frac{n}{2}}\left (\mu\right )-a_{\frac{n}{2}-1}\left ( \mu \right ) }{ a_{\frac{n}{2}-1}\left ( \mu \right )-a_{\frac{n}{2}-2}\left ( \mu \right ) }=2,
\end{equation}
\begin{equation}\label{even2}
a_{n-2}\left (\mu \right )=1.
\end{equation}
 \item Let $\mu$ be an arbitrary  root of equation $a_{\frac{n}{2}-1}(\mu)=0$. Then
\begin{equation}\label{even3}
1-\frac{ a_{\frac{n}{2}}\left (\mu\right )-a_{\frac{n}{2}-1}\left ( \mu \right ) }{ a_{\frac{n}{2}-1}\left ( \mu \right )-a_{\frac{n}{2}-2}\left ( \mu \right ) }=0,
\end{equation}
\begin{equation}\label{even4}
a_{n-2}\left (\mu \right )=-1.
\end{equation}
  \item  Let $\mu$ be an arbitrary  root of equation of equation $a_{\frac{n}{2}}(\mu)+a_{\frac{n}{2}-1}(\mu)=0$. Then
\begin{equation}\label{even5}
1-\frac{ a_{\frac{n}{2}}\left (\mu\right )-a_{\frac{n}{2}-1}\left ( \mu \right ) }{ a_{\frac{n}{2}-1}\left ( \mu \right )-a_{\frac{n}{2}-2}\left ( \mu \right ) }\in \left ( -\infty ,0 \right )\cup \left ( 2,+\infty  \right ),
\end{equation}
\begin{equation}\label{even6}
a_{n-1}\left ( \mu\right )=-1 ~~\textrm{and~~} a_{n-2}\left ( \mu  \right )=-2\left ( 1- \mu  \right ).
\end{equation}
  \end{enumerate}
\end{lemma}
The proof of the Lemma is presented in Appendix \ref{Prov_Pro_an2}.

\section{Main Results}
  \label{Main_result}
 % In this section, we firstly present the relation between  the eigenvalues of $\mathcal{L}_{G}$
%  For any connected graph $G$ with $N_{0}$ vertices,
%  if the  normalized Laplacian spectrum $\sigma _{G} =\left \{ \lambda_{1}, \lambda_{2},\cdots ,\lambda_{N_{0}} \right \}$ is known,  we will present the normalized Laplacian spectrum of  $n$-polygon graph $\tau_n(G)$, defined  by \texttt{Definition}~\ref{define1} in this section.

  For any connected undirected graph $G$ with $N_{0}$ vertices,   if the  normalized Laplacian spectrum $\sigma _{G} =\left \{ \lambda_{1}, \lambda_{2},\cdots ,\lambda_{N_{0}} \right \}$ is known,  we will present the normalized Laplacian spectrum of  $n$-polygon graph $\tau_n(G)$, defined  by \texttt{Definition}~\ref{define1} in this section. Firstly, in \textrm{Lemma} \ref{l21}, we present the general relation between the eigenvalues of $\mathcal{L}_{G}$ and  the eigenvalues   of $\mathcal{L}_{\tau_n(G)}$, where $\mathcal{L}_{G}$ and $\mathcal{L}_{\tau_n(G)}$ are the normalized Laplacian  matrices for graph $G$ and graph $\tau_n(G)$ respectively. Two \texttt{Remarks} follow \textrm{Lemma} \ref{l21} to clarify the detailed formulas of the relations in the case $n$ is odd and in the case $n$ is  even. Then, in \texttt{Theorem} \ref{t21} and \texttt{Theorem} \ref{t22}, we present the way to derive  the eigenvalues of $\mathcal{L}_{\tau_n(G)}$. By using \texttt{Theorem} \ref{t21} and \texttt{Theorem} \ref{t22} recursively, we can also derive  the eigenvalues of graph $\tau_n^g(G)$ ($n \geq 2, g\geq 1$), defined in \textrm{Definition}~\ref{define2}, which is the iterated $n$-polygon graph of $G$ with generation $g$.%Finally, in \texttt{Theorem} \ref{t23}, we present the normalized Laplacian spectrum of  graph $\tau_n(G)$ if the  normalized Laplacian spectrum of $G$ is given.

% 	For the graph $\tau_n(G)$ of $G$, the normalized Laplacian of it is written as $\mathcal{L}_{\tau_n(G)}$.Let the degree of the vertex $i\in N\left ( \tau_n(G) \right )$ be $d_{i}^{'}$. $A_{\tau_n(G)}$ is the adjacency matrix of and $P_{\tau_n(G)}=D_{\tau_n(G)}^{-\frac{1}{2}}A_{\tau_n(G)}D_{\tau_n(G)}^{-\frac{1}{2}}$,where $D_{\tau_n(G)}$ is the degree matrix of $\tau_n(G)$.In order to keep accordance,the normalized Lapacian of $G$ is denoted by $\mathcal{L}_{G}$ and let $P_{G}=D^{-\frac{1}{2}}AD^{-\frac{1}{2}}$.

\begin{lemma}\label{l21}
For any integer $n$ $(n\geq 2)$, let $\mu$ be a real number such that $a_{n-1}(\mu)\neq 0$ and $a_{n-1}(\mu)+1\neq 0$, where $\{a_{k}(\mu)\}$ is a series defined in Lemma \ref{Pro_Rec_an}. Then, $\mu$ is an eigenvalue of $\mathcal{L}_{\tau_n(G)}$ with  multiplicity $k$ ($k>0$) if and only if $\lambda\equiv1-\frac{a_{n}(\mu)}{1+a_{n-1}(\mu)}$ is an eigenvalue of $\mathcal{L}_{G}$ with  multiplicity $k$ and $\lambda\neq 0$, $\lambda\neq 2$.
\end{lemma}

The proof of this Lemma is presented in Sec.~\ref{Proof_Lemma121}.

\begin{remark}\label{Remark-Oddn}
   If $n$ is an odd number,  replacing  $a_{n}(\mu)$ and ${1+a_{n-1}(\mu)}$ from  Eqs.\eqref{a3} and \eqref{a4} respectively, we  get
%\begin{equation}\label{a5}
%\frac{a_{n}}{1+a_{n-1}}= \frac{\left ( a_{\frac{n-1}{2}+1}-a_{\frac{n-1}{2}-1} \right )a_{\frac{n-1}{2}}}{\left ( a_{\frac{n-1}{2}}-a_{\frac{n-1}{2}-2} \right )a_{\frac{n-1}{2}}}
%=\frac{ a_{\frac{n-1}{2}+1}-a_{\frac{n-1}{2}-1} }{ a_{\frac{n-1}{2}}-a_{\frac{n-1}{2}-2} }
%\end{equation}
\begin{eqnarray}\label{a5}
\lambda\equiv1-\frac{a_{n}(\mu)}{1+a_{n-1}(\mu)}&=& 1-\frac{\left ( a_{\frac{n-1}{2}+1}(\mu)-a_{\frac{n-1}{2}-1}(\mu) \right )a_{\frac{n-1}{2}}(\mu)}{\left ( a_{\frac{n-1}{2}}(\mu)-a_{\frac{n-1}{2}-2}(\mu) \right )a_{\frac{n-1}{2}}(\mu)}\nonumber\\
&=&1-\frac{ a_{\frac{n-1}{2}+1}(\mu)-a_{\frac{n-1}{2}-1}(\mu) }{ a_{\frac{n-1}{2}}(\mu)-a_{\frac{n-1}{2}-2}(\mu) },
\end{eqnarray}
and ${1+a_{n-1}(\mu)}\neq 0$  is equivalent to $a_{\frac{n-1}{2}}-a_{\frac{n-1}{2}-2}\neq 0$ and $a_{\frac{n-1}{2}}\neq 0$.
\end{remark}
\begin{remark}\label{Remark-Evenn}
   If $n$ is an even number,  replacing  $a_{n}(\mu)$ and ${1+a_{n-1}(\mu)}$ from  Eqs.\eqref{a3} and \eqref{a4} respectively, we  get
\begin{eqnarray}\label{a6}
\lambda\equiv1-\frac{a_{n}(\mu)}{1+a_{n-1}(\mu)}&=&1-\frac{\left ( a_{\frac{n}{2}}(\mu)-a_{\frac{n}{2}-1}(\mu) \right )\left ( a_{\frac{n}{2}}(\mu)+a_{\frac{n}{2}-1} (\mu)\right )}{\left ( a_{\frac{n}{2}-1}(\mu)-a_{\frac{n}{2}-2}(\mu) \right )\left ( a_{\frac{n}{2}}(\mu)+a_{\frac{n}{2}-1}(\mu) \right )}\nonumber\\
&=&1-\frac{ a_{\frac{n}{2}}(\mu)-a_{\frac{n}{2}-1}(\mu) }{ a_{\frac{n}{2}-1}(\mu)-a_{\frac{n}{2}-2}(\mu) },
\end{eqnarray}
and ${1+a_{n-1}(\mu)}\neq 0$  is equivalent to $a_{\frac{n}{2}-1}(\mu)-a_{\frac{n}{2}-2}(\mu)\neq 0$ and $ a_{\frac{n}{2}}(\mu)+a_{\frac{n}{2}-1}(\mu) \neq 0$.
\end{remark}

 \begin{theorem}\label{t21}
		Let $G$ be a  connected graph with $N_{0}$ vertices and $E_{0}$ edges, and $\tau_n(G)$ be the n-polygon graph of $G$, where $n \geq 3$ and $n$ is  odd. The eigenvalues for the normalized Laplacian  $\mathcal{L}_{\tau_n(G)}$ can be obtained in the following way.\\
        (i) $0$ is the eigenvalue of $\mathcal{L}_{\tau_n(G)}$ with the multiplicity $1$. If $G$ is bipartite, $2$ is the eigenvalue of $\mathcal{L}_{G\left ( n \right )}$ with the multiplicity $1$.\\
        (ii)Let $\mu$  be an arbitrary root of equation $a_{\frac{n-1}{2}}(\mu)=0$, then $\mu$ is a eigenvalue of $\mathcal{L}_{\tau_n(G)}$ with  multiplicity $N_{0}$.\\
        (iii)Let $\mu$  be an arbitrary root of equation $a_{\frac{n-1}{2}}(\mu)+a_{\frac{n-1}{2}-1}(\mu) =0$,  then $\mu$ is a eigenvalue of $\mathcal{L}_{\tau_n(G)}$ with  multiplicity $E_{0}-N_{0}+1$.\\
        (iv)In the case $G$ is non-bipartite, let $\mu$  be an arbitrary root of equation $a_{\frac{n-1}{2}}(\mu)-a_{\frac{n-1}{2}-1}(\mu) =0$,  then $\mu$ is a eigenvalue of  $\mathcal{L}_{\tau_n(G)}$ with multiplicity  $E_{0}-N_{0}$.\\
        (v)In the case  $G$ is bipartite, let $\mu$  be an arbitrary root of equation $a_{\frac{n-1}{2}}(\mu)-a_{\frac{n-1}{2}-1}(\mu) =0$,  then $\mu$ is a eigenvalue of  $\mathcal{L}_{\tau_n(G)}$ with multiplicity  $E_{0}-N_{0}+1$.\\
        (vi)Let $\lambda $ be an arbitrary eigenvalue of $\mathcal{L}_{G}$ such that $\lambda \neq 0$ and $\lambda \neq 2$,  $\mu_{i}\left ( \lambda \right )$, $( i=1,2,\cdots,\frac{n+1}{2})$ be the roots of equation
        \begin{equation}\label{b9}
         1-\frac{ a_{\frac{n-1}{2}+1}\left ( x \right )-a_{\frac{n-1}{2}-1}\left ( x \right ) }{ a_{\frac{n-1}{2}}\left ( x \right )-a_{\frac{n-1}{2}-2}\left ( x \right ) }=\lambda.
         \end{equation}
        Then,  $\mu_{i}\left ( \lambda \right )$, $( i=1,2,\cdots,\frac{n+1}{2})$  are eigenvalues of $\mathcal{L}_{\tau_n(G)}$ with $m_{\mathcal{L}_{G}}\left ( \lambda \right )=m_{\mathcal{L}_{\tau_n(G)}}\left ( \mu_{i}\left ( \lambda \right )  \right ),i=1,2,\cdots,\frac{n+1}{2}$.\\
	\end{theorem}

The proof of Theorem \ref{t21} is presented in Sec.~\ref{Proof-Theo-t21}.

\begin{theorem}\label{t22}
		Let $G$ be a  connected graph with $N_{0}$ vertices and $E_{0}$ edges, and $\tau_n(G)$ be the n-polygon graph of $G$, where $n \geq 2$ and $n$ is  even. The eigenvalues for the normalized Laplacian  $\mathcal{L}_{\tau_n(G)}$ can be obtained in the following way.\\
        (i) $0$ is the eigenvalue of $\mathcal{L}_{\tau_n(G)}$ with the multiplicity $1$.\\
       (ii)Let $\mu$  be an arbitrary root of equation $a_{\frac{n}{2}}(\mu)+a_{\frac{n}{2}-1}(\mu)=0$, then $\mu$ is eigenvalue of $\mathcal{L}_{\tau_n(G)}$ with the multiplicity $N_{0}$.\\
       (iii)Let $\mu$  be an arbitrary root of equation $a_{\frac{n}{2}-1}(\mu) =0$,  then $\mu$ is a eigenvalue of $\mathcal{L}_{\tau_n(G)}$ with  multiplicity $E_{0}-N_{0}+1$.\\
       (iv)In the case $G$ is non-bipartite, let $\mu$  be an arbitrary root of equation $a_{\frac{n}{2}}(\mu)-a_{\frac{n}{2}-1}(\mu) =0$,  then $\mu$ is a eigenvalue of  $\mathcal{L}_{\tau_n(G)}$ with multiplicity  $E_{0}-N_{0}$.\\
       (v)In the case  $G$ is bipartite, let $\mu$  be an arbitrary root of equation $a_{\frac{n}{2}}(\mu)-a_{\frac{n}{2}-1}(\mu) =0$,  then $\mu$ is a eigenvalue of  $\mathcal{L}_{\tau_n(G)}$ with multiplicity  $E_{0}-N_{0}+1$.\\
        (vi)Let $\lambda $ be an arbitrary eigenvalue of $\mathcal{L}_{G}$ such that $\lambda \neq 0$ and $\lambda \neq 2$,  $\mu_{i}\left ( \lambda \right )$, $( i=1,2,\cdots,\frac{n}{2})$ be the roots of equation
        \begin{equation}\label{bb9}
         1-\frac{ a_{\frac{n}{2}}\left ( x \right )-a_{\frac{n}{2}-1}\left ( x \right ) }{ a_{\frac{n}{2}-1}\left ( x \right )-a_{\frac{n}{2}-2}\left ( x \right ) }=\lambda.
         \end{equation}
        Then,  $\mu_{i}\left ( \lambda \right )$, $( i=1,2,\cdots,\frac{n}{2})$  are eigenvalues of $\mathcal{L}_{\tau_n(G)}$ with $m_{\mathcal{L}_{G}}\left ( \lambda \right )=m_{\mathcal{L}_{\tau_n(G)}}\left ( \mu_{i}\left ( \lambda \right )  \right ),i=1,2,\cdots,\frac{n}{2}$.\\
	\end{theorem}

The proof of the Theorem is presented in Sec.~\ref{Proof-Theo-t22}.
\begin{remark}\label{Remark-Oddn}
By using \texttt{Theorem} \ref{t21} and \texttt{Theorem} \ref{t22},  we can obtain the complete normalized Laplacian spectrum $\sigma(\tau_n(G))$ of   $n$-polygon graph $\tau_n(G)$. We can also derive the complete normalized Laplacian spectrum of graph $\tau_n^g(G)$ ($n \geq 2, g\geq 1$), defined in \textrm{Definition}~\ref{define2}, which is the iterated $n$-polygon graph of $G$ with generation $g$, by using \texttt{Theorem} \ref{t21} and \texttt{Theorem} \ref{t22} recursively. In fact, let $n=2$ and $n=4$ in \texttt{Theorem} \ref{t22} respectively, we can recover the results obtained in Ref.~\cite{XIE2016-AMC} and Ref.~\cite{XuWW-2019JMRA} respectively; let $n=3$ in \texttt{Theorem} \ref{t21}, we can recover the results obtained in Ref.~\cite{HuangLi2018LMA}.
\end{remark}

\section{Applications}
\label{Sec:Appli}
Recalling that the multiplicative degree-Kirchhoff index, the Kemeny's const and the number of spanning trees can be expressed as functions of the normalized Laplacian spectrum, as shown in Lemma~\ref{L2}, we obtain the following Theorems.
\begin{theorem}\label{t32}
Let $G$ be a  connected  undirected graph with $N_{0}$ vertices and $E_{0}$ edges, $\tau_n(G)$ ($n \geq 2$) be the $n$-polygon graph of $G$ and $\tau_n^g(G)$ ($n \geq 2, g\geq 1$), defined in \textrm{Definition}~\ref{define2}, be the iterated $n$-polygon graph of $G$ with generation $g$, and $Kf^{'}\left (G \right )$, $Kf^{'}\left (\tau_n(G) \right )$, $Kf^{'}\left (\tau_n^g(G)\right )$  represent  the multiplicative degree-Kirchhoff indexes of  graph  $G$, $\tau_n(G)$, $\tau_n^g(G)$ respectively. Then, for any $n \geq 2$,
%Let $G$ be a  connected graph with $N_{0}$ vertices and $E_{0}$ edges and $\tau_n(G)$ be the n-polygon graph of $G$, where $n \geq 3$ and $n$ is  odd.
%	Let $G$ be a connected graph with. The multiplicative degree-Kirchhoff index $Kf^{'}\left ( G\left ( n \right ) \right )$ of the n-gon graph $G\left ( n \right )$ is determined by the multiplicative degree-Kirchhoff index of $Kf^{'}\left ( G \right )$ the initial graph $G$ as follows \\
\begin{align}\label{d0}
Kf^{'}\left ( \tau_n(G) \right )=&\left ( n^{2}+n \right )Kf^{'}\left ( G \right )+\frac{2}{3}\left ( n+1 \right )\left ( n^{2}-1 \right )E_{0}^{2}-\frac{2}{3}\left ( n^{2}-1 \right )E_{0}N_{0}\notag \\
&-\frac{1}{6} \left ( n^{2}-1 \right )\left ( n-2 \right )E_{0},
\end{align}
and for $n \geq 2$ and $g\geq 1$,
\begin{align}\label{KFNG}
Kf^{'}\left ( \tau_n^g(G) \right )=&\left ( n^{2}+n \right )^{g}Kf^{'}\left ( G \right )-\frac{1}{3}\left ( n-2\right )\left ( n+1 \right )^{g}\left ( n^{g}-1 \right )E_{0}\notag\\
&+\frac{2\left ( n-1 \right )}{3n}\left ( n+1 \right )^{g}\left [ \left ( n+1 \right )^{g}\left ( n^{2}+1 \right )-n^{g+2}-1 \right ]E_{0}^{2}\notag\\
&-\frac{2}{3}\left ( n+1 \right )^{g}\left ( n^{g}-1 \right )E_{0}N_{0}.
\end{align}
	\end{theorem}
The proof of \textrm{Theorem}~\ref{t32} is presented in Sec.~\ref{Proof-t32}.  %Combing Theorem \ref{t32} and Lemma \ref{L2}.

\begin{theorem}\label{t33}
Let $G$ be a  connected  undirected graph with $N_{0}$ vertices and $E_{0}$ edges, and $K\left (G \right )$, $K\left (\tau_n(G) \right )$, $K\left (\tau_n^g(G)\right )$  represent  the Kemeny's constant of  graph  $G$, $\tau_n(G)$, $\tau_n^g(G)$ respectively. Then, for any $n \geq 2$, %	The Kemeny's constant for random walks on $G\left ( n \right )$ can be expressed by\\
\begin{displaymath}
K\left ( \tau_n(G) \right )=nK\left ( G \right )+\frac{1}{3}\left ( n^{2}-1 \right )E_{0}-\frac{1}{3}\left ( n-1 \right )N_{0}-\frac{1}{12} \left ( n-1 \right )\left ( n-2 \right ),
\end{displaymath}
and for $n \geq 2$ and $g\geq 1$,
\begin{align*}
 K\left (\tau_n^g(G) \right )
=& n^{g}K\left ( G \right )   -\frac{1}{3}\left ( n^{g}-1 \right )N_{0}-\frac{1}{12} \left ( n-2 \right )\left ( n^{g}-1 \right )+\\
&\left [ \frac{\left ( n-1 \right )\left ( n^{2}+1 \right )}{3n}\left ( n+1 \right )^{g}+\frac{1}{3}\left ( -n^{3}+n^{2}+1 \right )n^{g-1}-\frac{1}{3}\right ]E_{0}.
\end{align*}
	\end{theorem}
 %The proof of \textrm{Theorem}~\ref{t33} is presented in Sec.~\ref{Proof-t33}.
This Theorem can be obtained directly from  Theorem \ref{t32} and Lemma \ref{L2}.
\begin{theorem}\label{t34}
%Let $G$ be a  connected graph with $N_{0}$ vertices and $E_{0}$ edges, $\tau_n(G)$ ($n \geq 2$) be the $n$-polygon graph of $G$ and $\tau_n^g(G)$ ($n \geq 2, g\geq 1$), defined in \textrm{Definition}~\ref{define2}, be the iterated $n$-polygon graph of $G$ with generation $g$, and
Let $G$ be a  connected  undirected graph with $N_{0}$ vertices and $E_{0}$ edges and $N_{st}\left (G \right )$, $N_{st}\left (\tau_n(G) \right )$, $N_{st}\left (\tau_n^g(G)\right )$  represent  the number  of spanning trees of  graph  $G$, $\tau_n(G)$, $\tau_n^g(G)$ respectively. Then,  for $n \geq 2$,
\begin{equation}\label{d4}
N_{st}\left ( \tau_n(G) \right )=\left (n+1  \right )^{N_{0}-1} n^{E_{0}-N_{0}+1}\cdot  N_{st}\left ( G \right ),
\end{equation}
and for $n \geq 2$ and $g\geq 1$,
\begin{equation}\label{FNSTNG}
N_{st}\left (\tau_n^g(G)\right )=\left (n+1  \right )^{\frac{\left ( n-1 \right )\left [ \left ( n+1 \right )^{g}-ng-1 \right ]}{n^{2}}E_{0}+g N_{0}-g}\cdot n^{\frac{\left ( n+1 \right )^{g}+\left ( n-1 \right )ng-1}{n^{2}}E_{0}-g N_{0}+g} \cdot N_{st}\left ( G \right ).
\end{equation}
\end{theorem}
 The proof of \textrm{Theorem}~\ref{t34} is presented in Sec.~\ref{Proof-t34}.
 
\section{The proof of  Lemma~\ref{l21} }
\label{Proof_Lemma121}
First,  we  prove the  ``only if" part of this Lemma. \\
For an arbitrary  real number $\mu$, if  $a_{n-1}(\mu)\neq 0$ and $a_{n-1}(\mu)+1\neq 0$ and
$\mu$ is an eigenvalue of $\mathcal{L}_{\tau_n(G)}$ with  multiplicity $k$ ($k>0$), we will show that $\lambda\equiv1-\frac{a_{n}(\mu)}{1+a_{n-1}(\mu)}$ is an eigenvalue of $\mathcal{L}_{G}$ with  multiplicity $k$.

Let $\vec{v}= \left ( v_{1},v_{2},\cdots ,v_{N_{1}} \right )^{T}$ be an eigenvector with respect to the eigenvalue $\mu$ of $\tau_n(G)$. Thus %, with the entry $v_{i}$ corresponds to the vertex $i$ of graph $\tau_n(G)$ $\mu $ be an eigenvalue of $\mathcal{L}_{\tau_n(G)}$ and
		\begin{equation}\label{eq1}
			\mathcal{L}_{\tau_n(G)}\vec{v}=\mu\vec{v}.%=\left ( I-P_{\tau_n(G)} \right )\vec{v}
		\end{equation}

Replacing every entries of $\mathcal{L}_{\tau_n(G)}$ from Eq.~(\ref{Lij}), we find, for any vertex $q$ of graph $\tau_n(G)$, the corresponding entry $v_{q}$ in $\vec{v}$ satisfies %$q\in V\left ( \tau_n(G) \right )$
        \begin{equation}\label{eq2}
        \left ( 1-\mu \right )v_{q}=\sum_{t=1}^{N_{1}}\frac{A_{\tau_n(G)}(q, t)}{\sqrt{d_{q}^{'}d_{t}^{'}}}v_{t}, %\sum_{t=1}^{N_{1}}P_{\tau_n(G)}\left ( q,t \right )v_{t}=
        \end{equation}
where $A_{\tau_n(G)}(q, t)$  is the $(q, t)$ entry for the adjacency matrix of $\tau_n(G)$, $d_{q}^{'}$ and $d_{t}^{'}$ are the degree for vertexes $q$ and $t$ in graph $\tau_n(G)$ respectively.

Let $V(G)$ and $V(\tau_n(G))$ be the vertexes set of graph $G$ and $\tau_n(G)$ respectively and $V_{N}=\{v: v\in V(\tau_n(G)), v \notin V(G)\}$. We have $V(G)\subset V(\tau_n(G))$ and $V(\tau_n(G))=V_{N}\cup V(G)$.	We  will show the relation between the  corresponding entries of $\vec{v}$ for nodes in $V(G)$ and those for nodes in $V_{N}$.
	%Let $V_{N}$ be the set of all the newly added vertices in $\tau_n(G)$ and $V(G)$ be the set of the vertices inherited from $G$.That is,$V\left ( \tau_n(G) \right )=V_{N}\cup V(G)$.

For an arbitrary edge of graph $G$, let $i$ and $j$ be the two ends of the edge. Recalling the construction of the graph $\tau_n(G)$, as shown in Fig.~\ref{fig:1}, the edge is replacing with a $(n+1)$-polygon, whose nodes are labeled as $i$, $i_{1}^{'}$, $i_{2}^{'}$, $\cdots$, $i_{n-2}^{'}$, $i_{n-1}^{'}$, $j$.

   %For any vertex $i\in V(G)$, which is also a vertex of graph $\tau_n(G)$,
   Let $N_{O}^i\subseteq  V(G)$ denote the set for neighbors of vertex $i$ in graph $G$ and $N_{N}^i\subseteq  V_{N}$ denote the set of the new neighbors of vertex $i$ in graph $\tau_n(G)$. Note that $d_{i}^{'}=2d_i$ and $d_{j}^{'}=2d_j$. We can rewritte Eq.~\eqref{eq2} as%,can ben we have By the construction $\tau_n(G)$ of and
        \begin{eqnarray}\label{eq3}
        \left ( 1-\mu \right )v_{i}&=&\sum_{i_{1}^{'}\in N_N^i}\frac{1}{\sqrt{d_{i}^{'}d_{i_{1}^{'}}^{'}}}v_{i_{1}^{'}}+\sum_{j\in N_{O}^i}\frac{1}{\sqrt{d_{i}^{'}d_{j}^{'}}}v_{j}\nonumber\\
        &=&\sum_{i_{1}^{'}\in N_N^i}\frac{1}{2\sqrt{d_{i}}}v_{i_{1}^{'}}+\sum_{j\in N_{O}^i}\frac{1}{2\sqrt{d_{i}d_{j}}}v_{j}.
        \end{eqnarray}
       % where $i_{1}^{'}\in N_N^i$ and $j\in V(G)$ are neighbor vertices of $i$ in $\tau_n(G)$.\\
%        \begin{equation}\label{eq3}
%        \left ( 1-\mu \right )v_{i}=\sum_{i_{1}^{'}\in N_N^i}\frac{1}{\sqrt{d_{i}^{'}d_{i_{1}^{'}}^{'}}}v_{i_{1}^{'}}+\sum_{j\in V(G)}\frac{1}{\sqrt{d_{i}^{'}d_{j}^{'}}}v_{j}=\sum_{i_{1}^{'}\in N_N^i}\frac{1}{2\sqrt{d_{i}}}v_{i_{1}^{'}}+\sum_{j\in V(G)}\frac{1}{2\sqrt{d_{i}d_{j}}}v_{j} .
%        \end{equation}
%        where $i_{1}^{'}\in N_N^i$ and $j\in V(G)$ are neighbor vertices of $i$ in $\tau_n(G)$.\\
        Similarly, for any vertex $i_{1}^{'}\in N_N^i$, %it follows
        \begin{eqnarray}\label{eq4}
        \left ( 1-\mu \right )v_{i_{1}^{'}}&=&\frac{1}{\sqrt{d_{i_{1}^{'}}^{'}d_{i_{2}^{'}}^{'}}}v_{i_{2}^{'}}+
        \frac{1}{\sqrt{d_{i_{1}^{'}}^{'}d_{i}^{'}}}v_{i}\nonumber\\
        &=&\frac{1}{2}v_{i_{2}^{'}}+\frac{1}{2\sqrt{d_{i}}}v_{i},
        \end{eqnarray}
        %where $i_{2}^{'}\in V_{N}$ and $i\in V(G)$ are neighbors  of $i_{1}^{'}$ in $\tau_n(G)$.\\
       for the vertex $i_{n-1}^{'}\in N_{N}^j$,% which is adjacent to $i_{n-2}^{'}\in V_{N}$ and $j\in V(G)$,we obtain
        \begin{equation}\label{eq8}
        \left ( 1-\mu \right )v_{i_{n-1}^{'}}=\frac{1}{2}v_{i_{n-2}^{'}}+\frac{1}{2\sqrt{d_{j}}}v_{j},
        \end{equation}
       and  for the vertex $i_{k}^{'}$ $(2\leq k\leq n-2)$, %\in V_{N} which is adjacent to $i_{1}^{'}\in N_N^i$ and $i_{3}^{'}\in V_{N}$,we obtain
        \begin{equation}\label{eq5}
        \left ( 1-\mu \right )v_{i_{k}^{'}}=\frac{1}{2}\left ( v_{i_{k-1}^{'}}+v_{i_{k+1}^{'}} \right ).
        \end{equation}

Eliminating the  variables  $v_{i_{2}^{'}},v_{i_{3}^{'}},\cdots ,v_{i_{n-2}^{'}}$ in Eqs.~\eqref{eq4} and \eqref{eq5}, we have
      % So these are $ n-1 $ equations,and our goal is to get rid of $v_{i_{2}^{'}},v_{i_{3}^{'}},\cdots ,v_{i_{n-2}^{'}}$.We have

       \begin{equation}\label{eq12}
        a_{n-2}(\mu)v_{i_{1}^{'}}=v_{i_{n-1}^{'}}+a_{n-3}(\mu)\frac{v_{i}}{\sqrt{d_{i}}},
        \end{equation}
  
        where  $\{{a_{n}}(\mu)\}$ is just the series defined in Lemma~\ref{Pro_Rec_an}.

Similarly,  eliminating the  variables  $v_{i_{2}^{'}},v_{i_{3}^{'}},\cdots ,v_{i_{n-1}^{'}}$ in Eqs.~\eqref{eq4}, \eqref{eq8} and \eqref{eq5}, we get

%combining Eq.\eqref{eq9}, Eq.\eqref{eq11} and Eqs.\eqref{eq10}.  {eq3} {eq4} {eq8} {eq5} {eq2} And our goal is to get rid of $v_{i_{2}^{'}},v_{i_{3}^{'}},\cdots ,v_{i_{n-1}^{'}}$.We have
        \begin{equation}\label{eq15}
        a_{n-1}(\mu)v_{i_{1}^{'}}=a_{n-2}(\mu)\frac{1}{\sqrt{d_{i}}}v_{i}+\frac{1}{\sqrt{d_{j}}}v_{j}.
        \end{equation}
%where  $\{{a_{n}}(\mu)\}$ is just the series defined in Lemma~\ref{Pro_Rec_an}.

 Therefore, in the case $a_{n-1}(\mu)\neq 0$, multiplying both sides of the Eq.~(\ref{eq3}) by a factor $a_{n-1}(\mu)$, and replacing  $a_{n-1}(\mu)v_{i_{1}^{'}}$ from Eq.\eqref{eq15},  we have    %Combining Eq.\eqref{eq3} and Eq.\eqref{eq15}, for$a_{n-1}(\mu)\neq 0$, it follows
    \begin{equation}\label{eq16}
        2\left ( 1-\mu \right )a_{n-1}(\mu)v_{i}=a_{n-2}(\mu)v_{i}+\sum_{j\in N_{O}^i}\frac{1+a_{n-1}(\mu)}{\sqrt{d_{i}d_{j}}}v_{j},
     \end{equation}
   which is equivalent to %     Therefore,the equation
    \begin{equation}\label{eq1622}
       a_{n}(\mu)v_{i}=\sum_{j\in N_{O}^i}\frac{1+a_{n-1}(\mu)}{\sqrt{d_{i}d_{j}}}v_{j}.
    \end{equation}
    In the case $a_{n-1}(\mu)\neq 0$ and $1+a_{n-1}(\mu)\neq 0$, Eq.~(\ref{eq1622}) is also equal to
        \begin{equation}\label{eq17}
        \frac{a_{n}(\mu)}{1+a_{n-1}(\mu)}v_{i}=\sum_{j\in N_{O}^i}\frac{1}{\sqrt{d_{i}d_{j}}}v_{j}=\sum_{j\in V{(G)}}\frac{A_{ij}}{\sqrt{d_{i}d_{j}}}v_{j},
        \end{equation}
        which implies that $\frac{a_{n}(\mu)}{1+a_{n-1}(\mu)}$ is just an eigenvalue of the matrix $D(G)^{-\frac{1}{2}}A(G)D(G)^{-\frac{1}{2}}$.
        Therefore, in the case $a_{n-1}(\mu)\neq 0$ and $1+a_{n-1}(\mu)\neq 0$, $\lambda\equiv1-\frac{a_{n}(\mu)}{1+a_{n-1}(\mu)}$ is an eigenvalue of the matrix $\mathcal{L}_{G}=I-D(G)^{-\frac{1}{2}}A(G)D(G)^{-\frac{1}{2}}$.
    %    hold for $a_{n-1}\neq 0,1+a_{n-1}\neq 0$.\\

    Further more, Eqs.~\eqref{eq4}, \eqref{eq8} and \eqref{eq5} also show that the entries $v_{i_{k}^{'}}$ ($k=1$, $2$, $\cdots$, $n-1$) of $\vec{v}$ are completely decided by entries $v_i$ and $v_j$. This is to say, the  dimension of solution space for  linear equation $\mathcal{L}_{\tau_n(G)}\vec{v}=\mu\vec{v}$ is the same as that for linear equation $\mathcal{L}_{G}\vec{v'}=\lambda \vec{v'}$. Therefore, $\lambda\equiv1-\frac{a_{n}(\mu)}{1+a_{n-1}(\mu)}$ is an eigenvalue of the matrix $\mathcal{L}_{G}$  with the same multiplicity as the eigenvalue $\mu $ of $\mathcal{L}_{\tau_n(G)}$.

As shown in Lemma \ref{Pro_an1} and Lemma  \ref{lambda_0_2_2},    if $a_{n-1}(\mu)\neq 0$ and $1+a_{n-1}(\mu)\neq 0$, $\lambda\equiv1-\frac{a_{n}(\mu)}{1+a_{n-1}(\mu)}\neq 0$ and $\lambda\equiv1-\frac{a_{n}(\mu)}{1+a_{n-1}(\mu)}\neq 2$.

 Therefore, we obtain the  ``only if" part of this Lemma. Then we prove the ``if" part of this Lemma.

 Let $\lambda$ be an arbitrary eigenvalue of $\mathcal{L}_{G}$ with  multiplicity $k$ $(k>0)$, and
 $\vec{v'}= \left ( v_{1},v_{2},\cdots ,v_{N_{0}} \right )^{T}$ be an eigenvector with respect to the eigenvalue $\lambda$ of $\mathcal{L}_{G}$, i.e.,
		\begin{equation}
			\mathcal{L}_{G}\vec{v'}=\lambda \vec{v'},
		\end{equation}
   which can be rewritten as     %for any vertex $q\in V\left ( G\left ( n \right ) \right )$,Eq.\eqref{eq1} indicates that
        \begin{equation}\label{eq-LG}
        \left ( 1-\lambda \right )v_{i}=\sum_{j\in V{(G)}}\frac{A_{ij}}{\sqrt{d_{i}d_{j}}}v_{j}=\sum_{j\in N_{O}^i}\frac{1}{\sqrt{d_{i}d_{j}}}v_{j},
        \end{equation}
for any node $i\in V{(G)}$.  Therefore for any $\mu$ which satisfies $\lambda=1-\frac{a_{n}(\mu)}{1+a_{n-1}(\mu)}$, Eq.~(\ref{eq17}) holds and $1+a_{n-1}(\mu)\neq 0$.

Note that $\lambda\neq 0$ and $\lambda \neq 2$. As proved in Lemma \ref{Pro_an1} and Lemma  \ref{lambda_0_2_2}, we get $a_{n-1}(\mu)\neq 0$,, and for any node $i\in V{(G)}$,  Eq.~(\ref{eq16}) holds.

 For any edge with ends $i$ and $j$,  let $v_i$ and $v_j$ be the corresponding two entries in $\vec{v'}$ and  let
        \begin{equation} \label{Reik1}
        v_{i_{1}^{'}}=\frac{a_{n-2}(\mu)}{a_{n-1}(\mu)}\frac{1}{\sqrt{d_{i}}}v_{i}+\frac{1}{a_{n-1}(\mu)\sqrt{d_{j}}}v_{j},
        \end{equation}
        \begin{equation}\label{Reik2}
        v_{i_{2}^{'}}=2\left ( 1-\mu \right )v_{i_{1}^{'}}- \frac{1}{\sqrt{d_i}}v_i,
        \end{equation}
   and  for any $k$ $(2\leq k\leq n-2)$, %the vertex $i_{k}^{'}$ %\in V_{N} which is adjacent to $i_{1}^{'}\in N_N^i$ and $i_{3}^{'}\in V_{N}$,we obtain
        \begin{equation}\label{Reik}
        v_{i_{k+1}^{'}}=2\left ( 1-\mu \right )v_{i_{k}^{'}}- v_{i_{k-1}^{'}}.
        \end{equation}
In this way, we can obtain a $N_1$-dimensional vector $\vec{v}$, whose entries satisfy
 Eqs.~\eqref{eq3}-\eqref{eq5}. Then, for any node $q\in V(\tau_n(G))$, the entry $v_q$ in $\vec{v}$ satisfies Eq.~\eqref{eq2}.

Therefore,  the vector $\vec{v}$ satisfies Eq.~\eqref{eq1} and $\mu$ is  an  eigenvalue of  $\mathcal{L}_{\tau_n(G)}$. Note that the corresponding eigenvectors  $\vec{v}$  of $\mu$ is in one-to-one correspondence with  eigenvectors  $\vec{v^{'}}$ of $\lambda$. Thus the eigenvalue $\mu $ of $\mathcal{L}_{\tau_n(G)}$ has the same multiplicity as  eigenvalue $\lambda$ of $\mathcal{L}_{G}$.

This ends the proof.

\section{Proof of Theorem \ref{t21}}
\label{Proof-Theo-t21}

	%\begin{proof}
(i)It is obvious from Lemma \ref{L1}.\\
(ii) For any root $\mu$ of equation $a_{\frac{n-1}{2}}(\mu)=0$, we can obtain from Lemma \ref{Pro_an1} (see Eqs.\eqref{ANC2}) that %$\hat\mu_{i}$ and \eqref{b18}$(i=1,2,\cdots,\frac{n-1}{2})$
%$$a_{n-1}\left ( \hat\mu_{i}  \right )=-1, ~\textrm{and}~a_{n-2}\left ( \hat\mu_{i}  \right )=-2\left ( 1- \hat\mu_{i}  \right ).$$
$$a_{n-1}\left ( \mu  \right )=-1, ~\textrm{and}~a_{n-2}\left ( \mu  \right )=-2\left ( 1- \mu  \right ).$$
%\begin{align}\label{b17}
%a_{n-1}\left ( \hat\mu_{i}  \right )
%&=\left ( a_{\frac{n-1}{2}}\left ( \hat\mu_{i}  \right )-a_{\frac{n-1}{2}-1}\left ( \hat\mu_{i}  \right ) \right )\left ( a_{\frac{n-1}{2}}\left ( \hat\mu_{i}  \right )+a_{\frac{n-1}{2}-1}\left ( \hat\mu_{i}  \right ) \right )\notag\\
%&=-\left ( a_{\frac{n-1}{2}-1}\left ( \hat\mu_{i}  \right ) \right )^{2}\notag\\
%&=-1
%\end{align}
%Combining Eqs.\eqref{b15} and \eqref{b17},we have\\
%\begin{equation}\label{b18}
%a_{n-2}\left ( \hat\mu_{i}  \right )=-2\left ( 1- \hat\mu_{i}  \right )
%\end{equation}$\mu=\hat\mu_{i},i=1,2,\cdots,\frac{n-1}{2}$ and combining Eqs.\eqref{b17} and \eqref{b18}

Inserting the two equations  into Eq.\eqref{eq15}, we have\\
 \begin{equation}\label{b19}
       v_{i_{1}^{'}}=\frac{2(1-\mu)}{\sqrt{d_{i}}}v_{i}-\frac{1}{\sqrt{d_{j}}}v_{j}.
 \end{equation}
%  \begin{equation}\label{b19}
%        -v_{i_{1}^{'}}=\frac{-2(1-\hat\mu_{i})}{\sqrt{d_{i}}}v_{i}+\frac{1}{\sqrt{d_{j}}}v_{j}.
% \end{equation}
%Replacing $\mu$ with $\hat\mu_{i}$ in Eq.\eqref{eq3} and substituting $v_{i_{1}^{'}}$ from  Eq.\eqref{b19},  we get%and $\mu =\hat\mu_{i},i=1,2,\cdots,\frac{n-1}{2}$ into Eq.\eqref{eq3},we obtain\\
Therefore, for any node $i\in V(G)$,
\begin{eqnarray}\label{b20}
&&\sum_{i_{1}^{'}\in N_N^i}\frac{1}{\sqrt{d_{i}^{'}d_{i_{1}^{'}}^{'}}}v_{i_{1}^{'}}+\sum_{j\in V(G)}\frac{1}{\sqrt{d_{i}^{'}d_{j}^{'}}}v_{j}\nonumber\\
&=& \sum_{i_{1}^{'}\in N_N^i}\frac{1}{2\sqrt{d_{i}}}\left [ \frac{2(1-\mu)}{\sqrt{d_{i}}}v_{i}-\frac{1}{\sqrt{d_{j}}}v_{j} \right ]+\sum_{j\in N_{O}^i}\frac{1}{2\sqrt{d_{i}d_{j}}}v_{j}\nonumber\\
&=&\sum_{j\in N_{O}^i}\left [ -\frac{1}{2\sqrt{d_{i}d_{j}}}v_{j}+\frac{2(1-\mu)}{2{d_{i}}}v_{i} \right ]+\sum_{j\in N_{O}^i}\frac{1}{2\sqrt{d_{i}d_{j}}}v_{j}\nonumber\\
&=&\sum_{j\in N_{O}^i}\frac{(1-\mu)}{{d_{i}}}v_{i}\nonumber\\
&=&(1-\mu)v_{i},
\end{eqnarray}
where $N_N^i$ and $N_{O}^i$, which represent the  neighbors of node $i$ in $V_N$ and $V(G)$ respectively,  are defined in Lemma \ref{l21}.  %if $\mu=\hat\mu_{i}$,
Thus, Eq.\eqref{eq3} holds no matter what $v_{i}$ is. %This is to say, Eq. \eqref{b20}

Let $\vec{v'}= \left ( v_{1}, v_{2},\cdots , v_{N_{0}} \right )^{T}$ be an arbitrary vector. Eq. \eqref{b20} inform us that  Eq.\eqref{eq3} holds if $\mu$ is a root of equation $a_{\frac{n-1}{2}}(\mu)=0$. For  an arbitrary edge of graph $G$, let nodes $i$ and $j$ be the two ends of the edge. Considering the $(n+1)$-polygon as shown in the right-hand side of Fig.~\ref{fig:1}, and calculating $v_{i_{1}^{'}}$, $v_{i_{2}^{'}}$, $\cdots$, $v_{i_{n-1}^{'}}$ by using Eqs.\eqref{Reik1}, \eqref{Reik2}, \eqref{Reik}. In this way, we obtain
 a $N_1$-dimensional vector $\vec{v}$, whose entries satisfy  Eqs.~\eqref{eq3}-\eqref{eq5}. Then, for any node $q\in V(\tau_n(G))$, the entry $v_q$ in $\vec{v}$ satisfies Eq.~\eqref{eq2}.

Therefore,  the vector $\vec{v}$ satisfies Eq.~\eqref{eq1} and $\mu$ is  an  eigenvalue of  $\mathcal{L}_{\tau_n(G)}$. Further more, there are strict one-to-one correspondences be
the  eigenvectors  $\vec{v}$  of $\mu$ and  $\vec{v'}$, which is an arbitrary vector in $N_{0}$-dimensional space. Thus the  multiplicity of the eigenvalue $\mu $ is $N_{0}$.

(iii)  Let $\mu$  be an arbitrary root of equation $a_{\frac{n-1}{2}}(\mu)+a_{\frac{n-1}{2}-1}(\mu)=0$. Then
$$a_{n-1}\left (\mu \right )=[a_{\frac{n-1}{2}}(\mu)-a_{\frac{n-1}{2}-1}(\mu)][a_{\frac{n-1}{2}}(\mu)+a_{\frac{n-1}{2}-1}(\mu)]=0,$$
 and we can obtain from Lemma \ref{Pro_an1} (see Eq. \eqref{AN_2C1}) that
$$a_{n-2}\left (\mu \right )=-1.$$
Therefore, $$a_{n-3}\left (\mu \right )=2(1-\mu )a_{n-2}\left (\mu \right )-a_{n-1}\left (\mu \right )=-2\left ( 1-\mu \right ).$$
Replacing $a_{n-1}$ and $a_{n-2}$ with $0$ and $-1$ in Eq.\eqref{eq15} respectively, for any two nodes $i$ and $j$ of graph $G$, if $i\sim j$,
        \begin{equation}\label{b25}
        \frac{v_{i}}{\sqrt{d_{i}}}=\frac{v_{j}}{\sqrt{d_{j}}}.
        \end{equation}
Since $G$ is a connected graph, Eq.~(\ref{b25}) shows that $\frac{v_{i}}{\sqrt{d_{i}}}$ is a const for any $i\in V(G)$. Let $\theta\equiv\frac{v_{i}}{\sqrt{d_{i}}},~\textrm{for any} ~i\in  V(G)$. Eqs. \eqref{eq12} and \eqref{eq3} can be rewritten as\\
%because of the connectivity of $G$.For any edge $ij\in E(G)$,

%After substituting $\mu=\check\mu_{i},i=\frac{n-1}{2}+1,\frac{n-1}{2}+2,\cdots,n-1$ into Eq.\eqref{eq15},and combining Eqs.\eqref{b07} and \eqref{b09},we get
%       \begin{equation}\label{b25}
%       \frac{v_{i}}{\sqrt{d_{i}}}=\frac{v_{j}}{\sqrt{d_{j}}},\quad i\sim j,\quad i\in V(G),\quad j\in V(G).
%       \end{equation}
%Substituting $\mu$ into Eqs. \eqref{eq12} and \eqref{eq3}, respectively, we have\\i_{1}^{'}\not\sim i_{n-1}^{'},
       \begin{equation}\label{b26}
       v_{i_{1}^{'}}+v_{i_{n-1}^{'}}=2\left ( 1-\mu \right )\theta,\quad i_{1}^{'}\in N_N^i,i_{n-1}^{'}\in N_{N}^j.
       \end{equation}
       \begin{equation}\label{b27}
       \sum_{i_{1}^{'}\in N_N^i}v_{i_{1}^{'}}=\left [ 2\left ( 1-\mu \right )-1 \right ]\theta d_{i},\quad i\in V(G).
       \end{equation}
Calculating $\sum_{i\in V(G)}\sum_{i_{1}^{'}\in N_N^i}v_{i_{1}^{'}}$ by using Eq.\eqref{b26}, we have%i_{1}^{'}\not\sim i_{n-1}^{'}i\in V(G),  j\in V(G)
       \begin{equation}\label{b28}
       \sum_{i\in V(G)}\sum_{i_{1}^{'}\in N_N^i}v_{i_{1}^{'}}=\frac{1}{2}\sum_{i\sim j,  i_{1}^{'}\in N_{N}^i, i_{n-1}^{'}\in N_{N}^j}\left ( v_{i_{1}^{'}}+v_{i_{n-1}^{'}} \right )=2\left ( 1-\mu \right )\theta E_{0}.
       \end{equation}
On the other hand, Calculating $\sum_{i\in V(G)}\sum_{i_{1}^{'}\in N_N^i}v_{i_{1}^{'}}$ by using Eq.\eqref{b27}, we get \\
       \begin{equation}\label{b29}
       \sum_{i\in V(G)}\sum_{i_{1}^{'}\in N_N^i}v_{i_{1}^{'}}=\left [ 2\left ( 1-\mu \right )-1 \right ]\theta \sum_{i\in V(G)}d_{i}=\left [ 4\left ( 1-\mu \right )-2 \right ]\theta E_{0}.
       \end{equation}
Therefore, $\left [ 4\left ( 1-\mu \right )-2 \right ]\theta=2\left ( 1-\mu \right )\theta$, which leads to $\mu=0$ or $\theta=0$.

If $\mu=0$, replacing $a_{\frac{n-1}{2}}(\mu)$ and $a_{\frac{n-1}{2}-1}(\mu)$ from Eq.\eqref{mu=0}, we have  $$a_{\frac{n-1}{2}}(\mu)+a_{\frac{n-1}{2}-1}(\mu)=\frac{n-1}{2}+1+\frac{n-1}{2}-1+1=n \neq 0.$$
Thus, $\mu \neq 0$, and  $\theta=0$, which leads to% from Eqs.\eqref{b28} and \eqref{b29} $v_{i}=0$, %Therefore, the eigenvectors $\vec{v}=\left ( v_{1},v_{2},\cdots ,v_{N} \right )^{T}$ with respect to $\mu$ can be completely obtained by equations below\\
       \begin{equation}\label{b30}
       v_{i}=0,%\phantom{1} i\in V_{O};
       \end{equation}
for any node $i$ of graph $G$.

Substituting $v_{i}$ with $0$,  Eq.\eqref{eq3} can be rewritten as
       \begin{equation}\label{b31}
       \sum_{i_{1}^{'}\in N_N^i}v_{i_{1}^{'}}=0,
       \end{equation}
       for any node $i\in V(G)$.

Similarly,  replacing $v_{i}$ and $a_{n-2}\left (\mu \right )$ with $0$ and $-1$ respectively in  Eq. \eqref{eq12}, we get
       \begin{equation}\label{b32}
       v_{i_{1}^{'}}+v_{i_{n-1}^{'}}=0.
       \end{equation}

Further more, if $n$ is odd and $v_{i}=0$, in any $(n+1)$-polygon of $\tau_n(G)$ as shown in Fig.~\ref{fig:1}, we find, Eqs. \eqref{eq4}, \eqref{eq8} and \eqref{eq5} are equivalent to
\begin{align}\label{b33}
&v_{i_{t+1}^{'}}+v_{i_{n-t-1}^{'}}=0 , \quad t=0, ~1, ~2,~\cdots,~\frac{n-1}{2}-1,\\
&a_{t} ( \mu_{i} )v_{i_{1}^{'}}=v_{i_{t+1}^{'}}, \quad t=1, ~2,~\cdots,~\frac{n-1}{2}-1.\label{RVii}
\end{align}

Therefore,  Eq.~(\ref{eq2}) holds if and only if Eqs. \eqref{b30}, \eqref{b31}, \eqref{b33} and \eqref{RVii} hold for any $(n+1)$-polygon of graph $\tau_n(G)$. %This is to say,
%an vector $\vec{v}=\left ( v_{1},v_{2},\cdots ,v_{N_1} \right )^{T}$ is the eigenvector of $\mathcal{L}_{\tau_n(G)}$  corresponding to $\mu$ if and only if every entry $v_i$ for nodes in

Note that the total number of $(n+1)$-polygons in graph $\tau_n(G)$ is  $E_0$. For an arbitrary $(n+1)$-polygon, $P_r$, of graph $\tau_n(G)$, let $v_{i_{1}^{'}}=-v_{i_{n-1}^{'}}=x_{r}$,  and let $X=(x_{1}$, $x_{2}$, $\cdots$, $x_{E_0})^{T}$. We find Eqs.\eqref{b31} and \eqref{b32} are equivalent to  ${B}'X=0$, where ${B}'$  is the incident matrix of weakly connected directed graph ${G}'$ (see Definition \ref{define01} and \ref{define02}).

Similarly, for an arbitrary $(n+1)$-polygon, $P_r$, of graph $\tau_n(G)$, let
$$v_{i_{t+1}^{'}}=-v_{i_{n-t-1}^{'}}=y_{r}^{\left ( t \right )},~t=1,~2,~\ldots,~\frac{n-1}{2}-1,$$
and let $Y^{(t)}=(y_{1}^{(t)}$, $y_{2}^{(t)}$, $\cdots$, $y_{{E_0}}^{(t)})^{T}$,  we find,  Eqs. \eqref{b33} and \eqref{RVii} can be rewritten as
$$Y^{(t)}-a_{t} (\mu )X=0, \quad t=1, ~2,~\cdots,~\frac{n-1}{2}-1.$$
Therefore, there are strict one-to-one correspondences between the roots of Eq.~(\ref{eq1})  and roots of the equation

\begin{align}\label{matrix2}
\begin{pmatrix}
I & 0 & \cdots & 0 & -a_{1}\left (\mu \right )I\\
0 & I & \cdots & 0 & -a_{2}\left (\mu \right )I\\
\vdots  & \vdots & \ddots & \vdots  & \vdots\\
0 & 0 & \cdots & I & -a_{\frac{n-1}{2}-1}\left ( \mu \right )I\\
0 & 0 & \cdots & 0 & {B}'
\end{pmatrix}
\begin{pmatrix}
Y^{\left ( 1 \right )}\\
Y^{\left ( 2 \right )}\\
\vdots \\
Y^{\left ( \frac{n-1}{2}-1 \right )}\\
X
\end{pmatrix}
=
\begin{pmatrix}
0\\
0\\
\vdots \\
0\\
0
\end{pmatrix},
\end{align}
 where $I$ is an identity matrix of order $E_{0}$ and ${B}'$ is the incidence matrix of weakly connected directed graph ${G}'$ whose underlying undirected graph (all edges replaced by undirected edges) is $G$.

 Because $rank({B}')=N_{0}-1$ (see Lemma \ref{L03}). Therefore the rank for the coefficient matrix of linear equation (\ref{matrix2}) is $\frac{n-3}{2}E_0+N_0$ and the dimension for the vector space spanned by the roots of  linear equation (\ref{matrix2}) is $E_0-N_0+1$. Thus the dimension for the eigenspace of $\mu$ is also $E_0-N_0+1$ and the multiplicity  for eigenvalue $\mu$ is $E_0-N_0+1$.

(iv)
Let $\mu$  be an arbitrary root of equation $a_{\frac{n-1}{2}}(\mu)-a_{\frac{n-1}{2}-1}(\mu)=0$.  Then
$$a_{n-1}\left (\mu \right )=[a_{\frac{n-1}{2}}(\mu)-a_{\frac{n-1}{2}-1}(\mu)][a_{\frac{n-1}{2}}(\mu)+a_{\frac{n-1}{2}-1}(\mu)]=0,$$
 and we can obtain from Lemma \ref{Pro_an1} (see Eq. \eqref{b08}) that
$$a_{n-2}\left (\mu \right )=1.$$
Replacing $a_{n-1}$ and $a_{n-1}$ with $0$ and $1$ in Eq.\eqref{eq15}, for any two nodes $i$ and $j$ of graph $G$, if $i\sim j$,
%We can obtain from Lemma \ref{Pro_an1} (see Eqs.\eqref{b07} and \eqref{b08})
%that
%$$a_{n-1}\left (\mu \right )=0~\textrm{and}~a_{n-2}\left (\mu \right )=1.$$
%Plugging the two equations into Eq.\eqref{eq15},
        \begin{equation}\label{b21}
        \frac{v_{i}}{\sqrt{d_{i}}}=-\frac{v_{j}}{\sqrt{d_{j}}}.
        \end{equation}

Note that $G$ is non-bipartite, there is at least an odd cycle $C$ in $G$~\cite{Brouwer2012-BOOK, Jeribi2015-book}. Let $i_{1}$, $i_{2}$, $\cdots$, $i_{k}$ ($k$ is odd) be the nodes series of cycle $C$.
Eq.~(\ref{b21}) informs us that %We obtain Suppose is of length $s$ which its vertices are $i_{1},i_{2},\cdots ,i_{s}$ in turn.
$$\frac{v_{i_{1}}}{\sqrt{d_{i_{1}}}}=-\frac{v_{i_{2}}}{\sqrt{d_{i_{2}}}}=\frac{v_{i_{3}}}{\sqrt{d_{i_{3}}}}=\cdots =\frac{v_{i_{k}}}{\sqrt{d_{i_{k}}}}=-\frac{v_{i_{1}}}{\sqrt{d_{i_{1}}}}.$$
Therefore $v_{i_{1}}=v_{i_{2}}=\cdots=v_{i_{k}}=0$.  Since $G$ is connected,  Eq.~(\ref{b21}) informs us,
  \begin{equation}\label{b22}
       v_{i}=0.
 \end{equation}
 for any node $i$ of graph $G$.
%  \begin{equation}\label{b22}
%       v_{i}=0,\phantom{1} i\in V(G).
% \end{equation}

Substituting $v_{i}$ with $0$,  Eq.\eqref{eq3} can be rewritten as%  equivalent to
  \begin{equation}\label{b24}
       \sum_{i_{1}^{'}\in N_N^i}v_{i_{1}^{'}}=0,
   \end{equation}
  for any node $i$ of graph $G$.% and $N_N^i$ are set for neighbors of $i$ which .for any edge $ij\in E(G)$.

Similarly,  replacing $v_{i}$ and $a_{n-2}\left (\mu \right )$ with $0$ and $1$ respectively in  Eq. \eqref{eq12}, we get
  \begin{equation}\label{b23}
       v_{i_{1}^{'}}-v_{i_{n-1}^{'}}=0. %\quad i_{1}^{'}\not\sim i_{n-1}^{'},\quad i_{1}^{'},i_{n-1}^{'}\in V_{N},
  \end{equation}
Further more, if $n$ is odd and $v_{i}=0$ for any node $i\in V(G)$, we find, Eqs. \eqref{eq4}, \eqref{eq8} and \eqref{eq5} are equivalent to
\begin{align}\label{b34}
&v_{i_{t+1}^{'}}-v_{i_{n-t-1}^{'}}=0 , \quad t=0, ~1, ~2,~\cdots,~\frac{n-1}{2}-1,\\
&a_{t} ( \mu_{i} )v_{i_{1}^{'}}=v_{i_{t+1}^{'}}, \quad t=1, ~2,~\cdots,~\frac{n-1}{2}-1.\label{RVi}
\end{align}

Therefore,  Eq.~(\ref{eq2}) holds if and only if Eqs. \eqref{b22}, \eqref{b24}, \eqref{b34} and \eqref{RVi} hold for any $(n+1)$-polygon of graph $\tau_n(G)$. %This is to say,
%an vector $\vec{v}=\left ( v_{1},v_{2},\cdots ,v_{N_1} \right )^{T}$ is the eigenvector of $\mathcal{L}_{\tau_n(G)}$  corresponding to $\mu$ if and only if every entry $v_i$ for nodes in

Note that the total number of $(n+1)$-polygons in graph $\tau_n(G)$ is  $E_0$. For an arbitrary $(n+1)$-polygon, $P_r$, of graph $\tau_n(G)$, let $v_{i_{1}^{'}}=v_{i_{n-1}^{'}}=x_{r}$,  and let $X=(x_{1}$, $x_{2}$, $\cdots$, $x_{E_0})^{T}$. We find Eqs.\eqref{b24} and \eqref{b23} are equivalent to $BX=0$, where $B$  is the incident matrix of graph $G$.

Similarly, for an arbitrary $(n+1)$-polygon, $P_r$, of graph $\tau_n(G)$, let
$$v_{i_{t+1}^{'}}=v_{i_{n-t-1}^{'}}=y_{r}^{\left ( t \right )},~t=1,~2,~\ldots,~\frac{n-1}{2}-1,$$
and let $Y^{(t)}=(y_{1}^{(t)}$, $y_{2}^{(t)}$, $\cdots$, $y_{{E_0}}^{(t)})^{T}$,  we find,  Eqs. \eqref{b34} and \eqref{RVi} can be rewritten as
$$Y^{(t)}-a_{t} (\mu )X=0, \quad t=1, ~2,~\cdots,~\frac{n-1}{2}-1.$$
Therefore, there are strict one-to-one correspondences between the roots of Eq.~(\ref{eq1})  and roots of the equation
\begin{align}\label{matrix1}
\begin{pmatrix}
I & 0 & \cdots  & 0 & -a_{1}\left ( \mu \right )I\\
0 & I & \cdots  & 0 & -a_{2}\left ( \mu \right )I\\
\vdots & \vdots & \ddots & \vdots  & \vdots\\
0 & 0 & \cdots & I & -a_{\frac{n-1}{2}-1}\left ( \mu \right )I\\
0 & 0 & \cdots & 0 & B
\end{pmatrix}
\begin{pmatrix}
Y^{\left ( 1 \right )}\\
Y^{\left ( 2 \right )}\\
\vdots \\
Y^{\left ( \frac{n-1}{2}-1 \right )}\\
X
\end{pmatrix}
=
\begin{pmatrix}
0\\
0\\
\vdots \\
0\\
0
\end{pmatrix},
\end{align}
where $B$ is the incident matrix of $G$ and $I$ is identity matrix of order $E_{0}$.

Because $G$ is non-bipartite, $rank(B)=N_{0}$ (see Lemma \ref{L3}). Therefore the rank for the coefficient matrix of linear equation (\ref{matrix1}) is $\frac{n-3}{2}E_0+N_0$ and the dimension for the vector space spanned by the roots of  linear equation (\ref{matrix1}) is $E_0-N_0$. Thus the dimension for the eigenspace of $\mu$ is also $E_0-N_0$ and the multiplicity  for eigenvalue $\mu$ is $E_0-N_0$.

(v) Let $\mu$ be an arbitrary  root of equation $a_{\frac{n-1}{2}}(\mu)-a_{\frac{n-1}{2}-1}(\mu)=0$. That is to say
 $$a_{\frac{n-1}{2}}(\mu)=a_{\frac{n-1}{2}-1}(\mu).$$
As shown in  Lemma \ref{Pro_Rec_an}, for any $n\geq 1$,  $a_n(2-\mu)=(-1)^na_n(\mu)$.
Therefore,
  \begin{eqnarray}\label{Second_m_GFPT}
      a_{\frac{n-1}{2}}(2-\mu)
      &=&(-1)^{\frac{n-1}{2}}a_{\frac{n-1}{2}}(\mu)\nonumber\\
     &=&(-1)^{\frac{n-1}{2}}a_{\frac{n-1}{2}-1}(\mu)\nonumber\\
     &=&(-1)^{\frac{n-1}{2}}\times (-1)^{\frac{n-1}{2}-1}a_{\frac{n-1}{2}-1}(2-\mu)\nonumber\\
     &=&(-1)^{n-2}a_{\frac{n-1}{2}-1}(2-\mu)\nonumber\\
     &=&(-1)a_{\frac{n-1}{2}-1}(2-\mu).
  \end{eqnarray}
Therefore, $2-\mu$  is a root of equation $a_{\frac{n-1}{2}}(\mu)+a_{\frac{n-1}{2}-1}(\mu)=0$.

Therefore  $2-\mu$ is a  eigenvalue $\tau_n(G)$ with multiplicity $E_{0}-N_{0}+1$ (see (iii) of this Theorem).

 Note that $G$ is bipartite, then $\tau_n(G)$ is also bipartite.   We can obtain from Lemma \ref{Pro_Rec_an} that $\mu$ is also a  eigenvalue $\tau_n(G)$ with multiplicity $E_{0}-N_{0}+1$.

 (vi) It is just the result of Lemma \ref{l21}.

The proof is completed.

\section{Proof of Theorem \ref{t22}}
\label{Proof-Theo-t22}
This theorem can be proved in the same way as Theorem \ref{t21}.\\
(i) It is obvious from Lemma \ref{L1}.\\
(ii) Similar to the proof of Theorem\ref{t21}(ii). For any root $\mu$ of equation $a_{\frac{n}{2}}(\mu)+a_{\frac{n}{2}-1}(\mu)=0$, we can obtain from Lemma \ref{lambda_0_2_2} (see Eq.\eqref{even6}) that
$$ a_{n-1}\left ( \mu \right )-1, ~\textrm{and}~ a_{n-2}\left ( \mu  \right )=-2\left ( 1- \mu  \right ).$$
Inserting the two equations  into Eq.\eqref{eq15}, we have\\
 \begin{equation}\label{c19}
        v_{i_{1}^{'}}=\frac{2(1-\mu)}{\sqrt{d_{i}}}v_{i}-\frac{1}{\sqrt{d_{j}}}v_{j}.
        \end{equation}
Therefore, for any node $i\in V(G)$,
%After substituting Eq.\eqref{c19} and $\mu$ into Eq.\eqref{eq3},we obtain\\
\begin{eqnarray}\label{c20}
&&\sum_{i_{1}^{'}\in N_N^i}\frac{1}{\sqrt{d_{i}^{'}d_{i_{1}^{'}}^{'}}}v_{i_{1}^{'}}+\sum_{j\in V(G)}\frac{1}{\sqrt{d_{i}^{'}d_{j}^{'}}}v_{j}\nonumber\\
&=& \sum_{i_{1}^{'}\in N_N^i}\frac{1}{2\sqrt{d_{i}}}\left [ \frac{2(1-\mu)}{\sqrt{d_{i}}}v_{i}-\frac{1}{\sqrt{d_{j}}}v_{j} \right ]+\sum_{j\in N_{O}^i}\frac{1}{2\sqrt{d_{i}d_{j}}}v_{j}\nonumber\\
&=&\sum_{j\in N_{O}^i}\left [ -\frac{1}{2\sqrt{d_{i}d_{j}}}v_{j}+\frac{2(1-\mu)}{2{d_{i}}}v_{i} \right ]+\sum_{j\in N_{O}^i}\frac{1}{2\sqrt{d_{i}d_{j}}}v_{j}\nonumber\\
&=&\sum_{j\in N_{O}^i}\frac{(1-\mu)}{{d_{i}}}v_{i}\nonumber\\
&=&(1-\mu)v_{i},
\end{eqnarray}
%\begin{align}\label{c20}
%&\sum_{i_{1}^{'}\in N_N^i}\frac{1}{\sqrt{d_{i}^{'}d_{i_{1}^{'}}^{'}}}v_{i_{1}^{'}}+\sum_{j\in V(G)}\frac{1}{\sqrt{d_{i}^{'}d_{j}^{'}}}v_{j}\notag\\
%=&\sum_{i_{1}^{'}\in N_N^i}\frac{1}{2\sqrt{d_{i}}}\left [ \frac{2(1-\mu )}{\sqrt{d_{i}}}v_{i}-\frac{1}{\sqrt{d_{j}}}v_{j} \right ]+\sum_{j\in V_{O}}\frac{1}{2\sqrt{d_{i}d_{j}}}v_{j}\notag\\
%=&\sum_{j\in V_{O}}\left [ -\frac{1}{2\sqrt{d_{i}d_{j}}}v_{j}+\frac{2(1-\mu)}{2{d_{i}}}v_{i} \right ]+\sum_{j\in V_{O}}\frac{1}{2\sqrt{d_{i}d_{j}}}v_{j}\notag\\
%=&\sum_{j\in V_{O}}\frac{(1-\mu)}{{d_{i}}}v_{i}\notag\\
%=&(1-\mu)v_{i}
%\end{align}
%Eq.\eqref{c20} indicates Eq.\eqref{eq3} is an identical equation when $\mu=\tilde\mu _{i},i=1,2,\cdots,\frac{n}{2}$ and Eq.\eqref{c19} hold.Therefore,the eigenvectors associated with $\mu=\tilde\mu_{i},i=1,2,\cdots,\frac{n}{2}$ are completely determined by any $v_{i}$ and $v_{j}$ .We get $m_{\mathcal{L}_{G\left ( n \right )}}\left ( \tilde\mu_{i} \right )=N_{0},i=1,2,\cdots,\frac{n}{2}$.\\

where $N_N^i$ and $N_{O}^i$, which represent the  neighbors of node $i$ in $V_N$ and $V(G)$ respectively,  are defined in Lemma \ref{l21}.  %if $\mu=\hat\mu_{i}$,
Thus, Eq.\eqref{eq3} holds no matter what $v_{i}$ is. %This is to say, Eq. \eqref{b20}

Let $\vec{v'}= \left ( v_{1}, v_{2},\cdots , v_{N_{0}} \right )^{T}$ be an arbitrary vector. Eq. \eqref{c20} inform us that  Eq.\eqref{eq3} holds if $\mu$ is a root of equation $a_{\frac{n-1}{2}}(\mu)=0$. For  an arbitrary edge of graph $G$, let nodes $i$ and $j$ be the two ends of the edge. Considering the $(n+1)$-polygon as shown in the right-hand side of Fig.~\ref{fig:1}, and calculating $v_{i_{1}^{'}}$, $v_{i_{2}^{'}}$, $\cdots$, $v_{i_{n-1}^{'}}$ by using Eqs.\eqref{Reik1}, \eqref{Reik2}, \eqref{Reik}. In this way, we obtain a $N_1$-dimensional vector $\vec{v}$, whose entries satisfy  Eqs.~\eqref{eq3}-\eqref{eq5}. Then, for any node $q\in V(\tau_n(G))$, the entry $v_q$ in $\vec{v}$ satisfies Eq.~\eqref{eq2}.

Therefore,  the vector $\vec{v}$ satisfies Eq.~\eqref{eq1} and $\mu$ is  an  eigenvalue of  $\mathcal{L}_{\tau_n(G)}$. Further more, there are strict one-to-one correspondences be
the  eigenvectors  $\vec{v}$  of $\mu$ and  $\vec{v'}$, which is an arbitrary vector in $N_{0}$-dimensional space. Thus the  multiplicity of the eigenvalue $\mu $ is $N_{0}$.

(iii) %Similar to the proof of Theorem\ref{t21}(iii). \\
Let $\mu$ be an arbitrary root of equation $a_{\frac{n}{2}-1}(\mu)=0$.Then
$$a_{n-1}\left (\mu \right )=[a_{\frac{n}{2}}(\mu)-a_{\frac{n}{2}-2}(\mu)]a_{\frac{n}{2}-1}(\mu)=0,$$
 and we can obtain from Lemma \ref{lambda_0_2_2} (see Eq. \eqref{even4}) that
$$a_{n-2}\left (\mu \right )=-1.$$
Therefore,
$$a_{n-3}\left (\mu \right )=2(1-\mu )a_{n-2}\left (\mu \right )-a_{n-1}\left (\mu \right )=-2\left ( 1-\mu \right ).$$
Replacing $a_{n-1}$ and $a_{n-2}$ with $0$ and $-1$ in Eq.\eqref{eq15} respectively, for any two nodes $i$ and $j$ of graph $G$, if $i\sim j$,
        \begin{equation}\label{d25}
        \frac{v_{i}}{\sqrt{d_{i}}}=\frac{v_{j}}{\sqrt{d_{j}}}.
        \end{equation}
Since $G$ is a connected graph, Eq.~(\ref{d25}) shows that $\frac{v_{i}}{\sqrt{d_{i}}}$ is a const for any $i\in V(G)$. Let $\theta\equiv\frac{v_{i}}{\sqrt{d_{i}}},~\textrm{for any} ~i\in  V(G)$. Eqs. \eqref{eq12} and \eqref{eq3} can be rewritten as\\
%because of the connectivity of $G$.For any edge $ij\in E(G)$,

%After substituting $\mu=\check\mu_{i},i=\frac{n-1}{2}+1,\frac{n-1}{2}+2,\cdots,n-1$ into Eq.\eqref{eq15},and combining Eqs.\eqref{b07} and \eqref{b09},we get
%       \begin{equation}\label{b25}
%       \frac{v_{i}}{\sqrt{d_{i}}}=\frac{v_{j}}{\sqrt{d_{j}}},\quad i\sim j,\quad i\in V(G),\quad j\in V(G).
%       \end{equation}
%Substituting $\mu$ into Eqs. \eqref{eq12} and \eqref{eq3}, respectively, we have\\i_{1}^{'}\not\sim i_{n-1}^{'},
       \begin{equation}\label{d26}
       v_{i_{1}^{'}}+v_{i_{n-1}^{'}}=2\left ( 1-\mu \right )\theta,\quad i_{1}^{'}\in N_N^i,i_{n-1}^{'}\in N_{N}^j.
       \end{equation}
       \begin{equation}\label{d27}
       \sum_{i_{1}^{'}\in N_N^i}v_{i_{1}^{'}}=\left [ 2\left ( 1-\mu \right )-1 \right ]\theta d_{i},\quad i\in V(G).
       \end{equation}
Calculating $\sum_{i\in V(G)}\sum_{i_{1}^{'}\in N_N^i}v_{i_{1}^{'}}$ by using Eq.\eqref{d26}, we have%i_{1}^{'}\not\sim i_{n-1}^{'}i\in V(G),  j\in V(G)
       \begin{equation}\label{d28}
       \sum_{i\in V(G)}\sum_{i_{1}^{'}\in N_N^i}v_{i_{1}^{'}}=\frac{1}{2}\sum_{i\sim j,  i_{1}^{'}\in N_{N}^i, i_{n-1}^{'}\in N_{N}^j}\left ( v_{i_{1}^{'}}+v_{i_{n-1}^{'}} \right )=2\left ( 1-\mu \right )\theta E_{0}.
       \end{equation}
On the other hand, Calculating $\sum_{i\in V(G)}\sum_{i_{1}^{'}\in N_N^i}v_{i_{1}^{'}}$ by using Eq.\eqref{d27}, we get \\
       \begin{equation}\label{d29}
       \sum_{i\in V(G)}\sum_{i_{1}^{'}\in N_N^i}v_{i_{1}^{'}}=\left [ 2\left ( 1-\mu \right )-1 \right ]\theta \sum_{i\in V(G)}d_{i}=\left [ 4\left ( 1-\mu \right )-2 \right ]\theta E_{0}.
       \end{equation}
Therefore, $\left [ 4\left ( 1-\mu \right )-2 \right ]\theta=2\left ( 1-\mu \right )\theta$, which leads to $\mu=0$ or $\theta=0$.

If $\mu=0$, replacing $a_{\frac{n}{2}-1}(\mu)$ from Eq.\eqref{mu=0}, we have  $$a_{\frac{n}{2}-1}(\mu)=\frac{n}{2}-1+1=\frac{n}{2} \neq 0.$$
Thus, $\mu \neq 0$, and  $\theta=0$, which leads to% from Eqs.\eqref{b28} and \eqref{b29} $v_{i}=0$, %Therefore, the eigenvectors $\vec{v}=\left ( v_{1},v_{2},\cdots ,v_{N} \right )^{T}$ with respect to $\mu$ can be completely obtained by equations below\\
       \begin{equation}\label{d30}
       v_{i}=0,%\phantom{1} i\in V_{O};
       \end{equation}
for any node $i$ of graph $G$.

Substituting $v_{i}$ with $0$,  Eq.\eqref{eq3} can be rewritten as
       \begin{equation}\label{d31}
       \sum_{i_{1}^{'}\in N_N^i}v_{i_{1}^{'}}=0,
       \end{equation}
       for any node $i\in V(G)$.

Similarly,  replacing $v_{i}$ and $a_{n-2}\left (\mu \right )$ with $0$ and $-1$ respectively in  Eq. \eqref{eq12}, we get
       \begin{equation}\label{d32}
       v_{i_{1}^{'}}+v_{i_{n-1}^{'}}=0.
       \end{equation}

Further more, if $n$ is odd and $v_{i}=0$, in any $(n+1)$-polygon of $\tau_n(G)$ as shown in Fig.~\ref{fig:1}, we find, Eqs. \eqref{eq4}, \eqref{eq8} and \eqref{eq5} are equivalent to
\begin{align}\label{d33}
&v_{i_{t+1}^{'}}+v_{i_{n-t-1}^{'}}=0 , \quad t=0, ~1, ~2,~\cdots,~\frac{n}{2}-1,\\
&a_{t} ( \mu_{i} )v_{i_{1}^{'}}=v_{i_{t+1}^{'}}, \quad t=1, ~2,~\cdots,~\frac{n}{2}-1.\label{RVii}
\end{align}

Therefore,  Eq.~(\ref{eq2}) holds if and only if Eqs. \eqref{d30}, \eqref{d31}, \eqref{d33} and \eqref{RVii} hold for any $(n+1)$-polygon of graph $\tau_n(G)$. %This is to say,
%an vector $\vec{v}=\left ( v_{1},v_{2},\cdots ,v_{N_1} \right )^{T}$ is the eigenvector of $\mathcal{L}_{\tau_n(G)}$  corresponding to $\mu$ if and only if every entry $v_i$ for nodes in

Note that the total number of $(n+1)$-polygons in graph $\tau_n(G)$ is  $E_0$. For an arbitrary $(n+1)$-polygon, $P_r$, of graph $\tau_n(G)$, let $v_{i_{1}^{'}}=-v_{i_{n-1}^{'}}=x_{r}$,  and let $X=(x_{1}$, $x_{2}$, $\cdots$, $x_{E_0})^{T}$. We find Eqs.\eqref{d31} and \eqref{d32} are equivalent to  ${B}'X=0$, where ${B}'$  is the incident matrix of weakly connected directed graph ${G}'$ (see Definition \ref{define01} and \ref{define02}).

Similarly, for an arbitrary $(n+1)$-polygon, $P_r$, of graph $\tau_n(G)$, let
$$v_{i_{t+1}^{'}}=-v_{i_{n-t-1}^{'}}=y_{r}^{\left ( t \right )},~t=1,~2,~\ldots,~\frac{n}{2}-1,$$
and let $Y^{(t)}=(y_{1}^{(t)}$, $y_{2}^{(t)}$, $\cdots$, $y_{{E_0}}^{(t)})^{T}$,  we find,  Eqs. \eqref{d33} and \eqref{RVii} can be rewritten as
$$Y^{(t)}-a_{t} (\mu )X=0, \quad t=1, ~2,~\cdots,~\frac{n}{2}-1.$$
Therefore, there are strict one-to-one correspondences between the roots of Eq.~(\ref{eq1})  and roots of the equation

\begin{align}\label{matrix33}
\begin{pmatrix}
I & 0 & \cdots & 0 & -a_{1}\left (\mu \right )I\\
0 & I & \cdots & 0 & -a_{2}\left (\mu \right )I\\
\vdots  & \vdots & \ddots & \vdots  & \vdots\\
0 & 0 & \cdots & I & -a_{\frac{n}{2}-2}\left ( \mu \right )I\\
0 & 0 & \cdots & 0 & {B}'
\end{pmatrix}
\begin{pmatrix}
Y^{\left ( 1 \right )}\\
Y^{\left ( 2 \right )}\\
\vdots \\
Y^{\left ( \frac{n}{2}-2 \right )}\\
X
\end{pmatrix}
=
\begin{pmatrix}
0\\
0\\
\vdots \\
0\\
0
\end{pmatrix},
\end{align}
 where $I$ is an identity matrix of order $E_{0}$ and ${B}'$ is the incidence matrix of weakly connected directed graph ${G}'$ whose underlying undirected graph (all edges replaced by undirected edges) is $G$.

 Because $rank({B}')=N_{0}-1$ (see Lemma \ref{L03}). Therefore the rank for the coefficient matrix of linear equation (\ref{matrix33}) is $\frac{n-4}{2}E_0+N_0$ and the dimension for the vector space spanned by the roots of  linear equation (\ref{matrix33}) is $E_0-N_0+1$. Thus the dimension for the eigenspace of $\mu$ is also $E_0-N_0+1$ and the multiplicity  for eigenvalue $\mu$ is $E_0-N_0+1$.

(iv)
% Similar to the proof of Theorem\ref{t21}(iv). \\
Let $\mu$  be an arbitrary root of equation $a_{\frac{n}{2}}(\mu)-a_{\frac{n}{2}-2}(\mu)=0$.  Then
$$a_{n-1}\left (\mu \right )=[a_{\frac{n}{2}}(\mu)-a_{\frac{n}{2}-2}(\mu)]a_{\frac{n}{2}-1}(\mu)=0,$$
 and we can obtain from Lemma \ref{lambda_0_2_2} (see Eq. \eqref{even2}) that
$$a_{n-2}\left (\mu \right )=1.$$
Replacing $a_{n-1}$ and $a_{n-1}$ with $0$ and $1$ in Eq.\eqref{eq15}, for any two nodes $i$ and $j$ of graph $G$, if $i\sim j$,
        \begin{equation}\label{d21}
        \frac{v_{i}}{\sqrt{d_{i}}}=-\frac{v_{j}}{\sqrt{d_{j}}}.
        \end{equation}

Note that $G$ is non-bipartite, there is at least an odd cycle $C$ in $G$~\cite{Brouwer2012-BOOK}. Let $i_{1}$, $i_{2}$, $\cdots$, $i_{k}$ ($k$ is odd) be the nodes series of cycle $C$.
Eq.~(\ref{d21}) informs us that %We obtain Suppose is of length $s$ which its vertices are $i_{1},i_{2},\cdots ,i_{s}$ in turn.
$$\frac{v_{i_{1}}}{\sqrt{d_{i_{1}}}}=-\frac{v_{i_{2}}}{\sqrt{d_{i_{2}}}}=\frac{v_{i_{3}}}{\sqrt{d_{i_{3}}}}=\cdots =\frac{v_{i_{k}}}{\sqrt{d_{i_{k}}}}=-\frac{v_{i_{1}}}{\sqrt{d_{i_{1}}}}.$$
Therefore $v_{i_{1}}=v_{i_{2}}=\cdots=v_{i_{k}}=0$.  Since $G$ is connected,  Eq.~(\ref{d21}) informs us,
  \begin{equation}\label{d22}
       v_{i}=0.
 \end{equation}
 for any node $i$ of graph $G$.
%  \begin{equation}\label{b22}
%       v_{i}=0,\phantom{1} i\in V(G).
% \end{equation}

Substituting $v_{i}$ with $0$,  Eq.\eqref{eq3} can be rewritten as%  equivalent to
  \begin{equation}\label{d24}
       \sum_{i_{1}^{'}\in N_N^i}v_{i_{1}^{'}}=0,
   \end{equation}
  for any node $i$ of graph $G$.% and $N_N^i$ are set for neighbors of $i$ which .for any edge $ij\in E(G)$.

Similarly,  replacing $v_{i}$ and $a_{n-2}\left (\mu \right )$ with $0$ and $1$ respectively in  Eq. \eqref{eq12}, we get
  \begin{equation}\label{d23}
       v_{i_{1}^{'}}-v_{i_{n-1}^{'}}=0. %\quad i_{1}^{'}\not\sim i_{n-1}^{'},\quad i_{1}^{'},i_{n-1}^{'}\in V_{N},
  \end{equation}
Further more, if $n$ is odd and $v_{i}=0$ for any node $i\in V(G)$, we find, Eqs. \eqref{eq4}, \eqref{eq8} and \eqref{eq5} are equivalent to
\begin{align}\label{d34}
&v_{i_{t+1}^{'}}-v_{i_{n-t-1}^{'}}=0 , \quad t=0, ~1, ~2,~\cdots,~\frac{n}{2}-1,\\
&a_{t} ( \mu_{i} )v_{i_{1}^{'}}=v_{i_{t+1}^{'}}, \quad t=1, ~2,~\cdots,~\frac{n}{2}-1.\label{RVi}
\end{align}

Therefore,  Eq.~(\ref{eq2}) holds if and only if Eqs. \eqref{d22}, \eqref{d24}, \eqref{d34} and \eqref{RVi} hold for any $(n+1)$-polygon of graph $\tau_n(G)$. %This is to say,
%an vector $\vec{v}=\left ( v_{1},v_{2},\cdots ,v_{N_1} \right )^{T}$ is the eigenvector of $\mathcal{L}_{\tau_n(G)}$  corresponding to $\mu$ if and only if every entry $v_i$ for nodes in

Note that the total number of $(n+1)$-polygons in graph $\tau_n(G)$ is  $E_0$. For an arbitrary $(n+1)$-polygon, $P_r$, of graph $\tau_n(G)$, let $v_{i_{1}^{'}}=v_{i_{n-1}^{'}}=x_{r}$,  and let $X=(x_{1}$, $x_{2}$, $\cdots$, $x_{E_0})^{T}$. We find Eqs.\eqref{d24} and \eqref{d23} are equivalent to $BX=0$, where $B$  is the incident matrix of graph $G$.

Similarly, for an arbitrary $(n+1)$-polygon, $P_r$, of graph $\tau_n(G)$, let
$$v_{i_{t+1}^{'}}=v_{i_{n-t-1}^{'}}=y_{r}^{\left ( t \right )},~t=1,~2,~\ldots,~\frac{n}{2}-1,$$
and let $Y^{(t)}=(y_{1}^{(t)}$, $y_{2}^{(t)}$, $\cdots$, $y_{{E_0}}^{(t)})^{T}$,  we find,  Eqs. \eqref{d34} and \eqref{RVi} can be rewritten as
$$Y^{(t)}-a_{t} (\mu )X=0, \quad t=1, ~2,~\cdots,~\frac{n-1}{2}-1.$$
Therefore, there are strict one-to-one correspondences between the roots of Eq.~(\ref{eq1})  and roots of the equation

\begin{align}\label{matrix3}
\begin{pmatrix}
I & 0 & \cdots & 0 & -a_{1}\left ( \breve\mu_{i} \right )I\\
0 & I & \cdots & 0 & -a_{2}\left ( \breve\mu_{i} \right )I\\
\vdots  & \vdots & \ddots    & \vdots & \vdots\\
0 & 0 & \cdots & I & -a_{\frac{n}{2}-1}\left ( \breve\mu_{i} \right )I\\
0 & 0 & \cdots & 0 & B
\end{pmatrix}
\begin{pmatrix}
Y^{\left ( 1 \right )}\\
Y^{\left (2 \right )}\\
\vdots \\
Y^{\left ( \frac{n}{2}-1 \right )}\\
X
\end{pmatrix}
=
\begin{pmatrix}
0\\
0\\
\vdots \\
0\\
0
\end{pmatrix},
\end{align}
where $B$ is the incident matrix of $G$ and $I$ is identity matrix of order $E_{0}$.

Because $G$ is non-bipartite, $rank(B)=N_{0}$ (see Lemma \ref{L3}). Therefore the rank for the coefficient matrix of linear equation (\ref{matrix3}) is $\frac{n-2}{2}E_0+N_0$ and the dimension for the vector space spanned by the roots of  linear equation (\ref{matrix3}) is $E_0-N_0$. Thus the dimension for the eigenspace of $\mu$ is also $E_0-N_0$ and the multiplicity  for eigenvalue $\mu$ is $E_0-N_0$.

(v)
Let $\mu$  be an arbitrary root of equation $a_{\frac{n}{2}}(\mu)-a_{\frac{n}{2}-2}(\mu)=0$.  Then
$$a_{n-1}\left (\mu \right )=[a_{\frac{n}{2}}(\mu)-a_{\frac{n}{2}-2}(\mu)]a_{\frac{n}{2}-1}(\mu)=0,$$
 and we can obtain from Lemma \ref{lambda_0_2_2} (see Eq. \eqref{even2}) that
$$a_{n-2}\left (\mu \right )=1.$$
Replacing $a_{n-1}$ and $a_{n-2}$ with $0$ and $1$ in Eq.\eqref{eq15}, for any two nodes $i$ and $j$ of graph $G$, if $i\sim j$,
        \begin{equation}\label{f21}
        \frac{v_{i}}{\sqrt{d_{i}}}=-\frac{v_{j}}{\sqrt{d_{j}}}.
        \end{equation}
Note that $G$ is bipartite whose  vertex set  can be partitioned into two disjoint sets $X$ and $Y$, i.e., $V(G)=X\cup Y$. By  Eq.\eqref{f21},we have
\begin{equation}\label{f22}
\frac{v_{i}}{\sqrt{d_{i}}}=
\left\{\begin{matrix}
\theta ,i\in X\\
-\theta ,i\in Y
\end{matrix}\right.
\end{equation}
If $i\in X$ and $i\in Y$, Eq.\eqref{eq3} can be rewritten as, respectively,
       \begin{equation}\label{f23}
       \sum_{i_{1}^{'}\in N_N^i}v_{i_{1}^{'}}=\left [ 2\left ( 1-\mu \right )+1 \right ]\theta d_{i},\quad i\in X,
       \end{equation}
       \begin{equation}\label{f23}
       \sum_{i_{1}^{'}\in N_N^i}v_{i_{1}^{'}}=-\left [ 2\left ( 1-\mu \right )+1 \right ]\theta d_{i},\quad i\in Y.
       \end{equation}
If $i\in X$, Eq.\eqref{eq12} can be rewritten as
       \begin{equation}\label{f24}
       v_{i_{1}^{'}}-v_{i_{n-1}^{'}}=2\left ( 1-\mu \right )\theta,\quad i_{1}^{'}\in N_N^i,i_{n-1}^{'}\in N_{N}^j.
       \end{equation}
By calculating the difference of new neighbor nodes of nodes in bipartite graph $G$, we have
       \begin{equation}\label{f25}
       \sum_{i\in X}\sum_{i_{1}^{'}\in N_N^i}v_{i_{1}^{'}}-\sum_{i\in Y}\sum_{i_{1}^{'}\in N_N^i}v_{i_{1}^{'}}=\sum_{i\in X}\left ( v_{i_{1}^{'}}-v_{i_{n-1}^{'}} \right ).
       \end{equation}
Therefore,
       \begin{equation}\label{f26}
       \left [ 4\left ( 1-\mu \right )+2 \right ]\theta  E_{0}=2\left ( 1-\mu \right )\theta E_{0}.
       \end{equation}
Thus, $\theta=0$, which leads to
\begin{equation}\label{f27}
v_{i}=0,
\end{equation}
for any node $i$ of graph $G$.

Similarity to the proof of (iv), we find, in the case  $G$ is bipartite, there are strict one-to-one correspondences between the roots of Eq.~(\ref{eq1})  and roots of the equation Eq.~(\ref{matrix3}).
%in the case $G$ is bipartite, The proof of the remainder is similar to (iv).Note that $G$ is bipartite
However, different with non-bipartite graph,  here $rank(B)=N_{0}-1$ (see Lemma \ref{L3}).  Then, the dimension for the eigenspace of $\mu$ is $E_0-N_0+1$ and the multiplicity  for eigenvalue $\mu$ is $E_0-N_0+1$.

(vi)It is just the result of Lemma \ref{l21}.

The proof is completed.
\section{Proof of \textrm{Theorem}~\ref{t32}}
\label{Proof-t32}

 %\begin{proof}
Firstly, we prove  Eq. \eqref{d0} holds in the case $n$ is an odd number.
%\textrm{Theorem}~\ref{t32} is true if $n$ is an odd number.

Let $\lambda$ $(\lambda\neq 0~\textrm{and}~\lambda\neq 2)$ be an arbitrary eigenvalue of $\mathcal{L}_{G}$ and
$\mu_{i}\left ( \lambda \right )(i=1,~2,~\cdots,~\frac{n+1}{2})$ be the roots of equation
        \begin{equation}\label{fr-odd}
         1-\frac{ a_{\frac{n-1}{2}+1}\left ( x \right )-a_{\frac{n-1}{2}-1}\left ( x \right ) }{ a_{\frac{n-1}{2}}\left ( x \right )-a_{\frac{n-1}{2}-2}\left ( x \right ) }=\lambda.
         \end{equation}
       %  where $\lambda\in (0,2)$. \\
Since $a_{\frac{n-1}{2}}\left ( x \right )-a_{\frac{n-1}{2}-2}\left ( x \right ) \neq 0$, Eq.\eqref{fr-odd} can be rewritten as
\begin{equation}\label{f2}
f(x)=a_{\frac{n+1}{2}}\left ( x \right )-\left ( 1 - \lambda \right )a_{\frac{n-1}{2}}\left ( x \right )-a_{\frac{n-3}{2}}\left ( x \right )+\left ( 1 - \lambda \right )a_{\frac{n-5}{2}}\left ( x \right )=0,
\end{equation}
which is a  the polynomial equation of degree $\frac{n+1}{2}$. \\
Let $b_i$ be the  coefficient of $x^i$ $(i=1,~2,~\cdots,~\frac{n+1}{2})$ for polynomial $f(x)$.
According to Lemma \ref{Pro_Rec_an}(iv), we  get
\begin{align}\label{f3}
b_{0}=&a_{\frac{n+1}{2}}^{\left ( 0 \right )}\left ( x \right )-\left ( 1 - \lambda \right )a_{\frac{n-1}{2}}^{\left ( 0 \right )}\left ( x \right )-a_{\frac{n-3}{2}}^{\left ( 0 \right )}\left ( x \right )+\left ( 1 - \lambda \right )a_{\frac{n-5}{2}}^{\left ( 0 \right )}\left ( x \right )\notag \\
=&\frac{n+1}{2}+1-\left ( 1 - \lambda \right )\left ( \frac{n-1}{2}+1 \right )-\left ( \frac{n-3}{2}+1 \right )+\left ( 1 - \lambda \right )\left ( \frac{n-5}{2}+1 \right ) \notag \\
=&2\lambda,
\end{align}
\begin{align}\label{f4}
b_{1}=&a_{\frac{n+1}{2}}^{\left ( 1 \right )}\left ( x \right )-\left ( 1 - \lambda \right )a_{\frac{n-1}{2}}^{\left ( 1 \right )}\left ( x \right )-a_{\frac{n-3}{2}}^{\left ( 1 \right )}\left ( x \right )+\left ( 1 - \lambda \right )a_{\frac{n-5}{2}}^{\left ( 1 \right )}\left ( x \right )\notag\\
=&-\frac{\left ( \frac{n+1}{2}\right )^{3}+3\left ( \frac{n+1}{2}\right )^{2}+2\left ( \frac{n+1}{2}\right )}{3}+\left ( 1 - \lambda \right )\frac{\left ( \frac{n-1}{2}\right )^{3}+3\left ( \frac{n-1}{2}\right )^{2}+2\left ( \frac{n-1}{2}\right )}{3}\notag\\
&+\frac{\left ( \frac{n-3}{2} \right )^{3}+3\left ( \frac{n-3}{2} \right )^{2}+2\left ( \frac{n-3}{2} \right )}{3}-\left ( 1 - \lambda \right )\frac{\left ( \frac{n-5}{2} \right )^{3}+3\left ( \frac{n-5}{2} \right )^{2}+2\left ( \frac{n-5}{2} \right )}{3}\notag\\
=&-\frac{\left ( n-1 \right )^{2}}{2}\lambda-2n,
\end{align}
and
\begin{equation}\label{f5}
b_{\frac{n+1}{2}}=a_{\frac{n+1}{2}}^{\left (\frac{n+1}{2} \right )}\left ( x \right )=\left ( -1 \right )^{\frac{n+1}{2}}2^{\frac{n+1}{2}}.
\end{equation}

%By Lemma \ref{Pro_Rec_an}(iii), we can know $\mu_{i}\left ( \lambda \right )\neq 0,i=1,2,\cdots ,\frac{n+1}{2}$.
By using Lemma \ref{Vieta}, we have %and Eqs.\eqref{f3}-\eqref{f4}
% \begin{equation}\label{b10}
%\prod_{i=1}^{\frac{n+1}{2}}{\mu_{i}\left ( \lambda \right ) }=\frac{n}{\lambda }+\left ( \frac{n-1}{2} \right )^{2}.
%\end{equation}
\begin{equation}\label{b10}
\sum_{i=1}^{\frac{n+1}{2}}\frac{1}{\mu_{i}\left ( \lambda \right ) }=-\frac{b_{1}}{b_{0}}=\frac{n}{\lambda }+\left ( \frac{n-1}{2} \right )^{2}.
\end{equation}

%Using the same reason of Eq. \eqref{b10} above, we can easily get the following result:
Similarly, let $\check\mu_{i}( i=1,2,\cdots,\frac{n-1}{2})$ be the roots of equation $a_{\frac{n-1}{2}}-a_{\frac{n-1}{2}-1}=0$, $\check\mu_{i}( i=\frac{n-1}{2}+1,\frac{n-1}{2}+2,\cdots,n-1)$ be the roots of equation $a_{\frac{n-1}{2}}+a_{\frac{n-1}{2}-1}=0$, and $\hat\mu_{i}( i=1,2,\cdots,\frac{n-1}{2})$ be the roots of equation $a_{\frac{n-1}{2}}=0$. Rewriting these equations as polynomial equation and using Lemma \ref{Vieta}, we have
\begin{equation}\label{b03}
        \sum_{i=1}^{\frac{n-1}{2}}\frac{1}{\check\mu_{i}}=\frac{n^{2}-1}{4},
        \end{equation}
\begin{equation}\label{b05}
        \sum_{i=\frac{n-1}{2}+1}^{n-1}\frac{1}{\check\mu_{i}}=\frac{n^{2}-1}{12},
        \end{equation}
and
\begin{equation}\label{b13}
\sum_{i=1}^{\frac{n-1}{2}}\frac{1}{\hat\mu_{i}}=\frac{n^{2}+2n-3}{12}.
\end{equation}

{%\color{red}
Let $0=\lambda_{1}<\lambda_{2}\leq \lambda_{3}\leq\cdots\leq\lambda_{N_{0}}$ be the eigenvalues of the normalized Laplacian $\mathcal{L}_{G}$. If $G$ is non-bipartite,
Theorem \ref{t21} informs us that the normalized Laplacian spectrum of  graph $\tau_n(G)$ can be written as
%the eigenvalues of the normalized Laplacian $\mathcal{L}_{\tau_n(G)}$ can be listed as:
\begin{align}\label{Sp-Tn-Odd}
\sigma\left ( \tau_n(G) \right )=&\bigcup_{j=2}^{N_0}\left[\bigcup_{i=1}^{\frac{n+1}{2}}\{\mu_{i}\left ( \lambda_j \right )\}\right] \bigcup_{i=1}^{\frac{n-1}{2}}\underbrace{\left \{ \hat\mu _{i},\hat\mu_{i},\cdots ,\hat\mu_{i}\right \}}_{N_{0}}\bigcup_{i=1}^{\frac{n-1}{2}}\underbrace{\left \{ \check\mu_{i},\check\mu_{i},\cdots ,\check\mu_{i} \right \}}_{E_{0}-N_{0}} \notag\\
& \bigcup_{i=\frac{n-1}{2}+1}^{n-1}\underbrace{\left \{ \check\mu_{i},\check\mu_{i},\cdots ,\check\mu_{i} \right \}}_{E_{0}-N_{0}+1}\cup\left \{ 0 \right \}.
\end{align}

%¸ø³ö¶þ²¿Í¼»ò·Ç¶þ²¿Í¼µÄÖ¤Ã÷£¬ÁíÒ»ÖÖÀàËÆ
Therefore, in the case $n$ is an odd number and $G$ is non-bipartite, % and note that $\mu_{i}\left ( 2 \right )=\check\mu_{i},i=1,2,\cdots,\frac{n-1}{2}$ and $\mu_{\frac{n+1}{2}}\left ( 2 \right )=2$. Whether or not $G$ is bipartite, according to Theorem \ref{t21},
by using Lemma \ref{L2}(i) and Eqs \eqref{b10}-\eqref{b13},  we have
\begin{align}\label{d1}
&Kf^{'}\left ( \tau_n(G) \right )= 2E_{1}\sum_{u\in \sigma(\tau_n(G)),~\mu\neq0}\frac{1}{\mu}\notag \\
=&2E_{1}\left [ \sum_{i=2}^{N_{0}}\sum_{j=1}^{\frac{n+1}{2}}\frac{1}{\mu_{j}\left ( \lambda_{i} \right )}+\sum_{i=1}^{\frac{n-1}{2}}\frac{1}{\hat\mu_{i}}N_{0}+\sum_{i=1}^{\frac{n-1}{2}}\frac{1}{\check\mu _{i}}\left ( E_{0}-N_{0} \right )+\sum_{i=\frac{n-1}{2}+1}^{n-1}\frac{1}{\check\mu_{i}}\left ( E_{0}-N_{0}+1 \right ) \right ] \notag \\
=&2E_{1}\left [ \sum_{i=2}^{N_{0}}\left ( \frac{n}{\lambda_{i}}+\left ( \frac{n-1}{2} \right )^{2} \right )+\frac{n^{2}+2n-3}{12}N_{0}+\frac{n^{2}-1}{4}\left ( E_{0}-N_{0} \right )
\right. \notag \\
&\left.
+\frac{n^{2}-1}{12}\left ( E_{0}-N_{0}+1 \right ) \right ] \notag \\
=&2\left ( n+1 \right )E_{0}\left [ n\sum_{i=2}^{N_{0}} \frac{1}{\lambda _{i}}+\frac{n^{2}-1}{3}E_{0}-\frac{n-1}{3}N_{0}-\frac{\left ( n-1 \right )\left ( n-2 \right )}{6} \right ] \notag \\
=&\left ( n^{2}+n \right )Kf^{'}\left ( G \right )+\frac{2}{3}\left ( n+1 \right )\left ( n^{2}-1 \right )E_{0}^{2}-\frac{2}{3}\left ( n^{2}-1 \right )E_{0}N_{0} \notag \\
&-\frac{1}{6} \left ( n^{2}-1 \right )\left ( n-2 \right )E_{0}.
\end{align}
%Similarly, we can also verify Eq. \eqref{d1} holds in the case $n$ is an odd number and $G$ is bipartite.
%Therefore, \textrm{Theorem}~\ref{t32} is true in the case $n$ is an odd number.
Thus  Eq. \eqref{d0} holds in the case $n$ is an odd number and $G$ is non-bipartite.

If $G$ is bipartite, Theorem \ref{t21} informs us that the normalized Laplacian spectrum of  graph $\tau_n(G)$ can be written as
%the eigenvalues of the normalized Laplacian $\mathcal{L}_{\tau_n(G)}$ can be listed as:
\begin{align}\label{Sp-Tn-Odd-b}
\sigma\left ( \tau_n(G) \right )=&\bigcup_{j=2}^{N_0-1}\left[\bigcup_{i=1}^{\frac{n+1}{2}}\{\mu_{i}\left ( \lambda_j \right )\}\right] \bigcup_{i=1}^{\frac{n-1}{2}}\underbrace{\left \{ \hat\mu _{i},\hat\mu_{i},\cdots ,\hat\mu_{i}\right \}}_{N_{0}}\bigcup_{i=1}^{\frac{n-1}{2}}\underbrace{\left \{ \check\mu_{i},\check\mu_{i},\cdots ,\check\mu_{i} \right \}}_{E_{0}-N_{0}+1} \notag\\
& \bigcup_{i=\frac{n-1}{2}+1}^{n-1}\underbrace{\left \{ \check\mu_{i},\check\mu_{i},\cdots ,\check\mu_{i} \right \}}_{E_{0}-N_{0}+1}\cup\left \{ 0 \right \}\cup\left \{ 2 \right \}.
\end{align}
Similarly,  we can also verify Eq. \eqref{d0} holds in the case $n$ is an odd number and $G$ is bipartite.
Therefore, Eq. \eqref{d0} holds in the case $n$ is an odd number.
}

%Similarly, let $\check\mu_{i}( i=1,2,\cdots,\frac{n-1}{2})$ be the roots of equation $a_{\frac{n-1}{2}}-a_{\frac{n-1}{2}-1}=0$. We have
%\begin{equation}\label{b03}
%        \sum_{i=1}^{\frac{n-1}{2}}\frac{1}{\check\mu_{i}}=\frac{n^{2}-1}{4}.
%        \end{equation}
%
%Let $\check\mu_{i}( i=\frac{n-1}{2}+1,\frac{n-1}{2}+2,\cdots,n-1)$ be the roots of equation $a_{\frac{n-1}{2}}+a_{\frac{n-1}{2}-1}=0$.
%\begin{equation}\label{b05}
%        \sum_{i=\frac{n-1}{2}+1}^{n-1}\frac{1}{\check\mu_{i}}=\frac{n^{2}-1}{12}.
%        \end{equation}
%
%Let $\hat\mu_{i}( i=1,2,\cdots,\frac{n-1}{2})$ be the roots of equation $a_{\frac{n-1}{2}}=0$.
%\begin{equation}\label{b13}
%\sum_{i=1}^{\frac{n-1}{2}}\frac{1}{\hat\mu_{i}}=\frac{n^{2}+2n-3}{12}.
%\end{equation}

Then, we prove Eq. \eqref{d0} holds in the case $n$ is an even number.
%\textrm{Theorem}~\ref{t32} is true if $n$ is an even number.

Let $\lambda$ $(\lambda\neq 0~\textrm{and}~\lambda\neq 2)$ be an arbitrary eigenvalue of $\mathcal{L}_{G}$ and
$\mu_{i}\left ( \lambda \right )(i=1,~2,~\cdots,~\frac{n}{2})$ be the roots of equation
 \begin{equation}\label{bb9}
         1-\frac{ a_{\frac{n}{2}}\left ( x \right )-a_{\frac{n}{2}-1}\left ( x \right ) }{ a_{\frac{n}{2}-1}\left ( x \right )-a_{\frac{n}{2}-2}\left ( x \right ) }=\lambda.
         \end{equation}
Similar to the proof of Eq.~(\ref{b10}), we have
\begin{equation}\label{c10}
\sum_{i=1}^{\frac{n}{2}}\frac{1}{\mu_{i}\left ( \lambda \right ) }=\frac{n}{\lambda}+ \frac{n^{2}-2n}{4}.
\end{equation}
Similarly, let $\check\mu_{i}( i=1,2,\cdots,\frac{n}{2})$ be the roots of equation $a_{\frac{n}{2}}-a_{\frac{n}{2}-2}=0$, $\check\mu_{i}( i=\frac{n}{2}+1,\frac{n}{2}+2,\cdots,n-1)$ be the roots of equation  $a_{\frac{n}{2}-1}=0$, and $\hat\mu_{i}( i=1,2,\cdots,\frac{n}{2})$ be the roots of equation $ a_{\frac{n}{2}}+a_{\frac{n}{2}-1}=0$. We have
%Let $\mu_{i}\left ( \lambda \right )( i=1,2,\cdots,\frac{n+1}{2})$ be the roots of equation $a_{\frac{n}{2}}-a_{\frac{n}{2}-2}=0$.
\begin{equation}\label{c03}
        \sum_{i=1}^{\frac{n}{2}}\frac{1}{\breve\mu_{i}}=\frac{n^{2}}{4},
        \end{equation}
%Let $\breve\mu_{i} ( i=\frac{n}{2}+1,\frac{n}{2}+2,\cdots,n-1)$ be the roots of equation $a_{\frac{n}{2}-1}=0$.
\begin{equation}\label{c05}
        \sum_{i=\frac{n}{2}+1}^{n-1}\frac{1}{\breve\mu_{i}}=\frac{n^{2}-4}{12},
        \end{equation}
and
%Let $\hat\mu_{i}( i=1,2,\cdots,\frac{n}{2})$ be the roots of equation $ a_{\frac{n}{2}}+a_{\frac{n}{2}-1}=0$.
\begin{equation}\label{c13}
\sum_{i=1}^{\frac{n}{2}}\frac{1}{\hat\mu_{i}}=\frac{n^{2}+2n}{12}.
\end{equation}

Let $0=\lambda_{1}<\lambda_{2}\leq \lambda_{3}\leq\cdots\leq\lambda_{N_{0}}$ be the eigenvalues of the normalized Laplacian $\mathcal{L}_{G}$. Similar to the proof of Eq.~\eqref{d1}, in the case $n$ is an even number and $G$ is non-bipartite, we have
%Therefore,When $n$ is even number and note that $\mu_{i}\left ( 2 \right )=\breve\mu _{i},i=1,2,\cdots,\frac{n}{2}$. Whether or not $G$ is bipartite,
%By using Theorem \ref{t22}, Lemma \ref{L2}(i) and Eqs \eqref{c10}-\eqref{c13},
\begin{align}\label{d2}
&Kf^{'}\left ( \tau_n(G) \right ) \notag \\
=& 2E_{1}\left [ \sum_{i=2}^{N_{0}}\sum_{j=1}^{\frac{n}{2}}\frac{1}{\mu_{j}\left ( \lambda_{i} \right )}+\sum_{i=1}^{\frac{n}{2}}\frac{1}{\tilde\mu_{i}}N_{0}+\sum_{i=1}^{\frac{n}{2}}\frac{1}{\breve\mu _{i}}\left ( E_{0}-N_{0} \right )+\sum_{i=\frac{n}{2}+1}^{n-1}\frac{1}{\breve\mu_{i}}\left ( E_{0}-N_{0}+1 \right ) \right ]\notag \\
=& 2E_{1}\left [ \sum_{i=2}^{N_{0}}\left ( \frac{n}{\lambda_{i}}+ \frac{n^{2}-2n}{4}  \right )+\frac{n^{2}+2n}{12}N_{0}+\frac{n^{2}}{4}\left ( E_{0}-N_{0} \right )+\frac{n^{2}-4}{12}\left ( E_{0}-N_{0}+1 \right ) \right ] \notag \\
=& 2\left ( n+1 \right )E_{0}\left [ n\sum_{i=2}^{N_{0}} \frac{1}{\lambda _{i}}+\frac{n^{2}-1}{3}E_{0}-\frac{n-1}{3}N_{0}-\frac{\left ( n-1 \right )\left ( n-2 \right )}{6} \right ] \notag \\
=& \left ( n^{2}+n \right )Kf^{'}\left ( G \right )+\frac{2}{3}\left ( n+1 \right )\left ( n^{2}-1 \right )E_{0}^{2}-\frac{2}{3}\left ( n^{2}-1 \right )E_{0}N_{0}\notag \\
&-\frac{1}{6} \left ( n^{2}-1 \right )\left ( n-2 \right )E_{0}.
\end{align}
Thus  Eq. \eqref{d0} holds in the case $n$ is an even number and $G$ is non-bipartite.
Similarly, we can also verify Eq. \eqref{d0} holds in the case $n$ is an even number and $G$ is bipartite.
Therefore, Eq. \eqref{d0} holds for any  $n$ $(n\geq 2)$. %is an even number.

Next, we will prove Eq. \eqref{KFNG} holds for any $n\geq 2$ and $g\geq 1$. % $n$ $(n\geq 2)$ and $g$ $(g\geq 1)$ .

%Thus, Eq.\eqref{d1} is the same as Eq.\eqref{d2}.\\
Recalling the definition  of the iterated $n$-polygon graph  with generation $g$ (see Definition \ref{define2}), for any $g\geq 1$, % $\tau_n^g(G)$ is defined as the graph obtained through the iteration
$\tau_n^g(G)=\tau_n(\tau_n^{g-1}(G))$.
Therefore, we can obtain from from Eq.\eqref{d0} that, for any $g\geq 1$,

%\indent From Eq.\eqref{d0} and the definition of the $(n+1)$-polygon iterative graph we can get the recursive relation
\begin{align}\label{d3}
Kf^{'}\left ( \tau_n^g(G) \right )=&\left ( n^{2}+n \right )Kf^{'}\left ( \tau_{n}^{g-1}(G) \right )+\frac{2}{3}\left ( n+1 \right )\left ( n^{2}-1 \right )E_{g-1}^{2}\notag \\
&-\frac{2}{3}\left ( n^{2}-1 \right )E_{g-1}N_{g-1}-\frac{1}{6} \left ( n^{2}-1 \right )\left ( n-2 \right )E_{g-1}.
\end{align}
Using Eq.\eqref{d3} recursively and replacing $N_g$ from Eq.\eqref{ng1}, we have
%Therefore, by Eqs.\eqref{ng1} and \eqref{d3}, we have
\begin{align*}
&Kf^{'}\left ( \tau_n^g(G) \right )\\
=&\left ( n^{2}+n \right )^{g}Kf^{'}\left ( G \right )+\frac{2}{3}\left ( n+1 \right )\left ( n^{2}-1\right )\sum_{i=0}^{g-1}\left ( n^{2}+n \right )^{i}E_{g-1-i}^{2}-\frac{2}{3}\left ( n^{2}-1\right )\cdot \\
&\sum_{i=0}^{g-1}\left ( n^{2}+n \right )^{i}E_{g-1-i}N_{g-1-i}-\frac{1}{3}\left ( n^{2}-1 \right )\left ( n-2\right )\sum_{i=0}^{g-1}\left ( n^{2}+n \right )^{i}E_{g-1-i}\\
=&\left ( n^{2}+n \right )^{g}Kf^{'}\left ( G \right )+\frac{2}{3}\left ( n+1 \right )\left ( n^{2}-1\right )\left [ \left ( n+1 \right )^{2g-1}-n^{g}\left ( n+1 \right )^{g-1} \right ]E_{0}^{2}-\frac{2}{3}\left ( n^{2}
\right.\\
&\left.
-1\right ) \left \{ \frac{\left ( n+1 \right )^{g-1}\left ( n^{g}-1 \right )}{n-1}E_{0}N_{0}+ \frac{\left ( n+1 \right )^{g-1}\left [ \left ( n-1 \right )\left ( n+1 \right )^{g}-n^{g+1}+1 \right ]}{n}E_{0}^{2}\right \}\\
&-\frac{1}{3}\left ( n^{2}-1 \right )\left ( n-2\right )\frac{\left ( n+1 \right )^{g-1}\left ( n^{g}-1 \right )}{n-1}E_{0}\\
=&\left ( n^{2}+n \right )^{g}Kf^{'}\left ( G \right )+\frac{2\left ( n-1 \right )}{3n}\left ( n+1 \right )^{g}\left [ \left ( n+1 \right )^{g}\left ( n^{2}+1 \right )-n^{g+2}-1 \right ]E_{0}^{2} \\
&-\frac{1}{3}\left ( n-2\right )\left ( n+1 \right )^{g}\left ( n^{g}-1 \right )E_{0}-\frac{2}{3}\left ( n+1 \right )^{g}\left ( n^{g}-1 \right )E_{0}N_{0}.
\end{align*}
Then we obtain  Eq. \eqref{KFNG}.
%\end{proof}

\section{Proof of \textrm{Theorem}~\ref{t34}}
\label{Proof-t34}
%\begin{proof}
Firstly, we prove Eq. \eqref{d4} holds in the case $n$ is an odd number.
%\textrm{Theorem}~\ref{t32} is true if $n$ is an even number.

Let $\lambda$ $(\lambda\neq 0~\textrm{and}~\lambda\neq 2)$ be an arbitrary eigenvalue of $\mathcal{L}_{G}$ and
$\mu_{i}\left ( \lambda \right )(i=1,~2,~\cdots,~\frac{n}{2})$ be the roots of equation
%Let $\mu_{i}\left ( \lambda \right )( i=1,2,\cdots,\frac{n+1}{2})$ be the roots of equation
        \begin{equation}\label{Rf1}
         1-\frac{ a_{\frac{n-1}{2}+1}\left ( x \right )-a_{\frac{n-1}{2}-1}\left ( x \right ) }{ a_{\frac{n-1}{2}}\left ( x \right )-a_{\frac{n-1}{2}-2}\left ( x \right ) }=\lambda.
         \end{equation}
%By Lemma \ref{Pro_Rec_an}(iii), we can know $\mu_{i}\left ( \lambda \right )\neq 0,i=1,2,\cdots ,\frac{n+1}{2}$.
Note that Eq.\eqref{Rf1} can be rewritten as Eq. \eqref{f2}, which is a polynomial equation, with some coefficients shown as Eqs.\eqref{f3} and \eqref{f5}.
By using Lemma \ref{Vieta}, we have
\begin{equation}\label{b11}
\prod_{i=1}^{\frac{n+1}{2}}\mu_{i}\left ( \lambda \right )=\left ( -1 \right )^{n}\frac{b_{0}}{b_{\frac{n+1}{2}}}=\frac{\lambda }{2^{\frac{n-1}{2}}}.
\end{equation}

%Using the same reason of Eq. \eqref{b11} above, we can easily get the following result:
Similarly, let  $\check\mu_{i}( i=1,2,\cdots,\frac{n-1}{2})$ be the roots of equation $a_{\frac{n-1}{2}}-a_{\frac{n-1}{2}-1}=0$, $\check\mu_{i}( i=\frac{n-1}{2}+1,\frac{n-1}{2}+2,\cdots,n-1)$ be the roots of equation $a_{\frac{n-1}{2}}+a_{\frac{n-1}{2}-1}=0$, and $\hat\mu_{i}( i=1,2,\cdots,\frac{n-1}{2})$ be the roots of equation $a_{\frac{n-1}{2}}=0$. Rewriting these equations as polynomial equation and using Lemma \ref{Vieta},
we have
\begin{equation}\label{b04}
        \prod_{i=1}^{\frac{n-1}{2}}\check\mu_{i}=\frac{1}{2^{\frac{n-1}{2}}},
        \end{equation}
 %Let $\check\mu_{i}( i=\frac{n-1}{2}+1,\frac{n-1}{2}+2,\cdots,n-1)$ be the roots of equation $a_{\frac{n-1}{2}}+a_{\frac{n-1}{2}-1}=0$.
\begin{equation}\label{b06}
        \prod_{i=\frac{n-1}{2}+1}^{n-1}\check\mu_{i}=\frac{n}{2^{\frac{n-1}{2}}},
        \end{equation}
and%Let $\hat\mu_{i}( i=1,2,\cdots,\frac{n-1}{2})$ be the roots of equation $a_{\frac{n-1}{2}}=0$.
\begin{equation}\label{b14}
\prod_{i=1}^{\frac{n-1}{2}}\hat\mu_{i}=\frac{{n+1}}{2^{\frac{n+1}{2}}}.
\end{equation}

Let $0=\lambda_{1}<\lambda_{2}\leq \lambda_{3}\leq\cdots\leq\lambda_{N_{0}}$ be the eigenvalues of the normalized Laplacian $\mathcal{L}_{G}$. If $G$ is non-bipartite,
Note that the normalized Laplacian spectrum $\sigma( \tau_n(G))$ of  graph $\tau_n(G)$ are shown in  Eq. \eqref{Sp-Tn-Odd}, we have
%the eigenvalues of the normalized Laplacian $\mathcal{L}_{\tau_n(G)}$ can be listed as:
%\begin{align}\label{Sp-Tn-Odd}
%\sigma\left ( \tau_n(G) \right )=&\bigcup_{j=1}^{N_0}\left[\bigcup_{i=1}^{\frac{n+1}{2}}\{\mu_{i}\left ( \lambda_j \right )\}\right] \bigcup_{i=1}^{\frac{n-1}{2}}\underbrace{\left \{ \hat\mu _{i},\hat\mu_{i},\cdots ,\hat\mu_{i}\right \}}_{N_{0}}\bigcup_{i=1}^{\frac{n-1}{2}}\underbrace{\left \{ \check\mu_{i},\check\mu_{i},\cdots ,\check\mu_{i} \right \}}_{E_{0}-N_{0}} \notag\\
%& \bigcup_{i=\frac{n-1}{2}+1}^{n-1}\underbrace{\left \{ \check\mu_{i},\check\mu_{i},\cdots ,\check\mu_{i} \right \}}_{E_{0}-N_{0}+1}\cup\left \{ 0 \right \}.
%\end{align}
%When $n$ is odd number and note that $\mu_{i}\left ( 2 \right )=\check\mu_{i},i=1,2,\cdots,\frac{n-1}{2}$ and $\mu_{\frac{n+1}{2}}\left ( 2 \right )=2$. Whether or not $G$ is bipartite, by Theorem \ref{t21}, Lemma \ref{L2}(iii) and Eqs \eqref{b11}-\eqref{b14}, we get
\begin{align}\label{d6}
&\prod_{\mu\in\sigma( \tau_n(G)), \mu\neq 0 }\mu\notag\\
&=\left ( \prod_{i=1}^{\frac{n-1}{2}}\hat\mu_{i} \right )^{N_{0}}\times \left ( \prod_{i=1}^{\frac{n-1}{2}}\check\mu_{i} \right )^{E_{0}-N_{0}}\times \left ( \prod_{i=\frac{n-1}{2}+1}^{n-1}\check\mu_{i} \right )^{E_{0}-N_{0}+1}\times \prod_{i=2}^{N_{0}}\prod_{j=1}^{\frac{n-1}{2}}\mu_{j}\left ( \lambda_{i} \right )\notag \\
&=\left ( \frac{n+1}{2^{\frac{n+1}{2}}} \right )^{N_{0}}\times \left ( \frac{1}{2^{\frac{n-1}{2}}} \right )^{E_{0}-N_{0}}\times \left ( \frac{n}{2^{\frac{n-1}{2}}} \right )^{E_{0}-N_{0}+1}\times \prod_{i=2}^{N_{0}}\frac{\lambda_{i}}{2^{\frac{n-1}{2}}}\notag \\
&=\frac{n^{E_{0}-N_{0}+1}}{\left ( 2^{n-1} \right )^{E_{0}}}\left ( \frac{n+1}{2} \right )^{N_{0}}\times \prod_{i=2}^{N_{0}}\lambda_{i}.
\end{align}
Therefore,
\begin{equation}\label{d5}
\frac{N_{st}\left ( G\left ( n \right ) \right )}{N_{st}\left ( G \right )}=\frac{2^{N_{1}} \prod_{\mu\in\sigma( \tau_n(G)), \mu\neq 0 }\mu}{\left (n+1  \right )\prod_{i=2}^{N_{0}}\lambda_{i}}=n^{E_{0}-N_{0}+1}(n+1)^{N_0-1}.
\end{equation}
Thus, Eq. \eqref{d4} holds in the case $n$ is an odd number and  $G$ is non-bipartite.

In the case $n$ is an odd number and  $G$ is bipartite, recalling the normalized Laplacian spectrum $\sigma( \tau_n(G))$ of  graph $\tau_n(G)$ are shown in  Eq. \eqref{Sp-Tn-Odd-b},    we can also verify that Eq. \eqref{d4} holds.  Therefore Eq. \eqref{d4} holds in the case $n$ is an odd number.

Then, we prove Eq. \eqref{d4} holds in the case $n$ is an even number.
%\textrm{Theorem}~\ref{t32} is true if $n$ is an even number.

Let $\lambda$ $(\lambda\neq 0~\textrm{and}~\lambda\neq 2)$ be an arbitrary eigenvalue of $\mathcal{L}_{G}$ and
$\mu_{i}\left ( \lambda \right )(i=1,~2,~\cdots,~\frac{n}{2})$ be the roots of equation
 \begin{equation}\label{bb9}
         1-\frac{ a_{\frac{n}{2}}\left ( x \right )-a_{\frac{n}{2}-1}\left ( x \right ) }{ a_{\frac{n}{2}-1}\left ( x \right )-a_{\frac{n}{2}-2}\left ( x \right ) }=\lambda.
         \end{equation}
Similar to the proof of Eq.~(\ref{b11}), we have
\begin{equation}\label{c11}
\prod_{i=1}^{\frac{n}{2}}\mu_{i}\left ( \lambda \right )=\frac{\lambda}{2^{\frac{n}{2}}}.
\end{equation}
Similarly, let $\check\mu_{i}( i=1,2,\cdots,\frac{n}{2})$ be the roots of equation $a_{\frac{n}{2}}-a_{\frac{n}{2}-2}=0$, $\check\mu_{i}( i=\frac{n}{2}+1,\frac{n}{2}+2,\cdots,n-1)$ be the roots of equation  $a_{\frac{n}{2}-1}=0$, and $\hat\mu_{i}( i=1,2,\cdots,\frac{n}{2})$ be the roots of equation $ a_{\frac{n}{2}}+a_{\frac{n}{2}-1}=0$. We have

%For the case $n$ is an even number,  %Using the same reason of Eq. \eqref{b11} above, we can easily get the following result:

%Let $\mu_{i}\left ( \lambda \right )( i=1,2,\cdots,\frac{n}{2})$ be the roots of equation
% \begin{equation}\label{bb9}
%         1-\frac{ a_{\frac{n}{2}}\left ( x \right )-a_{\frac{n}{2}-1}\left ( x \right ) }{ a_{\frac{n}{2}-1}\left ( x \right )-a_{\frac{n}{2}-2}\left ( x \right ) }=\lambda.
%         \end{equation}
%\begin{equation}\label{c11}
%\prod_{i=1}^{\frac{n}{2}}\mu_{i}\left ( \lambda \right )=\frac{\lambda}{2^{\frac{n}{2}}}.
%\end{equation}

%Let $\mu_{i}\left ( \lambda \right )( i=1,2,\cdots,\frac{n+1}{2})$ be the roots of equation $ a_{\frac{n}{2}}-a_{\frac{n}{2}-2}=0$.
\begin{equation}\label{c04}
        \prod_{i=1}^{\frac{n}{2}}\check\mu_{i}=\frac{1}{2^{\frac{n}{2}-1}},
        \end{equation}
%Let $\breve\mu_{i} ( i=\frac{n}{2}+1,\frac{n}{2}+2,\cdots,n-1)$ be the roots of equation $a_{\frac{n}{2}-1}=0$.
\begin{equation}\label{c06}
        \prod_{i=\frac{n}{2}+1}^{n-1}\breve\mu_{i}=\frac{n}{2^{\frac{n}{2}}},
        \end{equation}
and%Let $\hat\mu_{i}( i=1,2,\cdots,\frac{n}{2})$ be the roots of equation $ a_{\frac{n}{2}}+a_{\frac{n}{2}-1}=0$.
\begin{equation}\label{c14}
\prod_{i=1}^{\frac{n}{2}}\hat\mu_{i}=\frac{{n+1}}{2^{\frac{n}{2}}}.
\end{equation}

Calculating the  normalized Laplacian spectrum $\sigma( \tau_n(G))$ of  graph $\tau_n(G)$ by using Theorem \ref{t22} in the case  $n$ is even number,
we can also get Eq.~\eqref{d5}.
%\begin{equation}\label{d5}
%\frac{N_{st}\left ( G\left ( n \right ) \right )}{N_{st}\left ( G \right )}=\frac{2^{N_{1}} \prod_{\mu\in\sigma( \tau_n(G)), \mu\neq 0 }\mu}{\left (n+1  \right )\prod_{i=2}^{N_{0}}\lambda_{i}}=n^{E_{0}-N_{0}+1}(n+1)^{N_0-1}.
%\end{equation}
Thus, Eq. \eqref{d4} holds in the case $n$ is an even number.

%\begin{align}\label{d7}
%&\prod_{\mu\in\sigma( \tau_n(G)), \mu\neq 0 }\mu\notag\\
%&=\left ( \prod_{i=1}^{\frac{n}{2}}\tilde\mu_{i} \right )^{N_{0}}\times \left ( \prod_{i=1}^{\frac{n}{2}}\breve\mu_{i} \right )^{E_{0}-N_{0}}\times \left ( \prod_{i=\frac{n}{2}+1}^{n-1}\breve\mu_{i} \right )^{E_{0}-N_{0}+1}\times \prod_{i=2}^{N_{0}}\prod_{j=1}^{\frac{n}{2}}\mu_{j}\left (\lambda_{i} \right )\notag \\
%&=\left ( \frac{n+1}{2^{\frac{n}{2}}} \right )^{N_{0}}\times \left ( \frac{1}{2^{\frac{n}{2}-1}} \right )^{E_{0}-N_{0}}\times \left ( \frac{n}{2^{\frac{n}{2}}} \right )^{E_{0}-N_{0}+1}\times \prod_{i=2}^{N_{0}}\frac{\lambda_{i}}{2^{\frac{n}{2}}}\notag \\
%&=\frac{n^{E_{0}-N_{0}+1}}{\left ( 2^{n-1} \right )^{E_{0}}}\left ( \frac{n+1}{2} \right )^{N_{0}}\times \prod_{i=2}^{N_{0}}\lambda_{i}.
%\end{align}
%Thus, Eq. \eqref{d6}is the same as Eq.\eqref{d7}. Therefore, Eq.\eqref{d4} holds by Eqs.\eqref{d5} and \eqref{d6}.\\
%%\begin{displaymath}
%%N_{st}\left ( G\left ( n \right ) \right )=\left (n+1  \right )^{N_{0}-1} n^{E_{0}-N_{0}+1}\cdot  N_{st}\left ( G \right )
%%\end{displaymath}
Recalling the definition  of the iterated $n$-polygon graph  with generation $g$ (see Definition \ref{define2}), for any $g\geq 1$, % $\tau_n^g(G)$ is defined as the graph obtained through the iteration
$\tau_n^g(G)=\tau_n(\tau_n^{g-1}(G))$.
Therefore, we can obtain from from Eq.\eqref{d4} that, for any $g\geq 1$,
\begin{align}\label{NSTR}
N_{st}\left ( \tau_n^g(G) \right )&=\left (n+1  \right )^{N_{g-1}-1}n^{E_{g-1}-N_{g-1}+1}\cdot  N_{st}\left ( \tau_{n}^{g-1}(G) \right ).
\end{align}
By using  Eq.\eqref{NSTR} recursively, we obtain Eq.\eqref{FNSTNG}.

\section{Conclusion}
Given the normalized Laplacian spectrum of an arbitrary connected graph $G$, we have obtained the normalized Laplacian spectrum of the graph $\tau_n(G)$, which is obtained by replacing each edge of $G$ with  a $(n+1)$-polygon $(n\geq 2)$. Then, the  normalized Laplacian spectrum of the graphs  $\tau_n^g(G)$ ($g\geq 0$) can also be obtained. As applications, we also calculated the  multiplicative degree-Kirchhoff index, the Kemeny's const and the number of spanning trees for graphs $\tau_n^g(G)$ $(g\geq 0)$. Note that `edge replacing' is a common graph operator which were widely used to construct different meaningful networks. The results  obtained here would be helpful for network optimization.

\section*{Acknowledgement(s)}

The authors  would like to express their sincere gratitude to all former researchers who present original results they quote in this manuscript.

\section*{Disclosure statement}

No potential conflict of interest was reported by the authors.

\section*{Funding}

The work is supported by the National Natural Science Foundation of China (Grant No. 61873069 and 61772147).

\section*{ORCID}

Junhao Peng,  https://orcid.org/0000-0001-8839-7070

%\bibliographystyle{tfnlm}
%
%\bibliography{GZPref}

%\section*{Acknowledgment}
%
%The authors are grateful to Elena Agliari for valuable suggestions. JHP is supported by the National Natural Science Foundation of China (Grant No. 61873069 and 61772147) and the National Key R\&D Program of China (Grant No. 2018YFB0803604). TS was supported by the Alexander von Humboldt Foundation. \\
%
%
%
%
\bibliographystyle{elsarticle-num}
\bibliography{GZPref}

\appendix

\section{The proof of the Lemma~\ref{Pro_Rec_an}}
\label{Proof_Pro_Rec_an}

Here, we  use `$a_{k}$' to represent `$a_{k}(\mu)$' to lighten the notations.
Firstly, we solve the recursive relation which series $\{a_{n}\}$ satisfies and present the general formula of $\{a_{n}\}$.
Let $\beta=2(1-\mu)$. Then the recursive  relation of series $\{a_{n}(\mu)\}$ can be rewritten as %, \mu\in[0,2], then $\beta\in[-2,2]$
  \begin{equation}\label{Re05}
  a_{n}=\beta a_{n-1}-a_{n-2},
    \end{equation}
with initial conditions $a_{-1}=0 $, $a_{0}=1$.
The characteristic equation of equation (\ref{Re05}) is $x^{2}-\beta x+1=0$, with roots $x=\frac{\beta \pm i\sqrt{4-\beta ^{2}}}{2}$. %, where $i$ is imaginary number, we have \\
  Then $\{a_{n}(\mu)\}$ can be rewritten as
    \begin{equation}\label{a05}
       a_{n}=c_{1}\left (\frac{\beta +i\sqrt{4-\beta ^{2}}}{2}\right )^{n}+c_{2}\left (\frac{\beta -i\sqrt{4-\beta ^{2}}}{2}\right )^{n}.
       \end{equation}
  By using the initial conditions $a_{-1}=0 $, $a_{0}=1$, we have\\
       \begin{equation}\label{a06}
      \left\{ \begin{array}{l} a_{-1}=c_{1}\left (\frac{\beta +i\sqrt{4-\beta ^{2}}}{2}\right )^{-1}+c_{2}\left (\frac{\beta -i\sqrt{4-\beta ^{2}}}{2}\right )^{-1}=0\\
       a_{0}=c_{1}+c_{2}=1 \end{array} \right.,
       \end{equation}
       which yields
     %  \begin{equation}\label{a7}
%       a_{0}=c_{1}+c_{2}=1.
%       \end{equation}
%      Solving Eqs.\eqref{a06} and \eqref{a7} gives\\
       \begin{equation}\label{a07}
      \left\{ \begin{array}{l}  c_{1}=\frac{1}{2}-\frac{i\beta }{2\sqrt{4-\beta ^{2}}}\\
       c_{2}=\frac{1}{2}+\frac{i\beta }{2\sqrt{4-\beta ^{2}}} \end{array} \right..
       \end{equation}
%       \begin{displaymath}
%        c_{1}=\frac{1}{2}-\frac{i\beta }{2\sqrt{4-\beta ^{2}}},
%       \end{displaymath}
%        \begin{displaymath}
%       c_{2}=\frac{1}{2}+\frac{i\beta }{2\sqrt{4-\beta ^{2}}}.
%       \end{displaymath}
    Therefore,%   Substitute the result into Eq.\eqref{a05}, we can get
\begin{eqnarray}\label{a1}
&&a_{n}(\mu)\nonumber\\%=\alpha_n(\beta)
&=&(\frac{1}{2}-\frac{i\beta }{2\sqrt{4-\beta ^{2}}})\left (\frac{\beta +i\sqrt{4-\beta ^{2}}}{2}\right )^{n}+(\frac{1}{2}+\frac{i\beta }{2\sqrt{4-\beta ^{2}}})\left (\frac{\beta -i\sqrt{4-\beta ^{2}}}{2}\right )^{n}\nonumber\\
&\equiv&a_n(\beta),
\end{eqnarray}
where $\beta=2(1-\mu)$.\\
Then, we prove the results of Lemma~\ref{Pro_Rec_an} item-by-item.

(i) Note that $\beta=2(1-\mu)$. Then $2[1-(2-\mu)]=2(\mu-1)=-\beta$. Therefore, %By Eq.\eqref{a1}, we have
\begin{align*}
&a_{n}(2-\mu)=a_n(-\beta)\\
=&(\frac{1}{2}-\frac{-i\beta }{2\sqrt{4-\beta ^{2}}})\left (\frac{-\beta +i\sqrt{4-\beta ^{2}}}{2}\right )^{n}+(\frac{1}{2}+\frac{-i\beta }{2\sqrt{4-\beta ^{2}}})\left (\frac{-\beta -i\sqrt{4-\beta ^{2}}}{2}\right )^{n}\\
=&\left ( -1 \right )^{n}\left [(\frac{1}{2}+\frac{i\beta }{2\sqrt{4-\beta ^{2}}})\left (\frac{\beta -i\sqrt{4-\beta ^{2}}}{2}\right )^{n} +(\frac{1}{2}-\frac{i\beta }{2\sqrt{4-\beta ^{2}}})\left (\frac{\beta +i\sqrt{4-\beta ^{2}}}{2}\right )^{n}\right ]\\
=&\left ( -1 \right )^{n}a_n(\beta)=\left ( -1 \right )^{n}a_{n}(\mu).
\end{align*}

(ii) Note that $\beta=2(1-\mu)$. Then
\begin{equation}
  \beta=\left\{ \begin{array}{ll}2 & \text{if $\mu=0$}\\ -2 & \text{if $\mu=2$} \end{array} \right.\nonumber.
\end{equation}
 Replacing $\beta$ with $2$ and $-2$ respectively in Eq. \eqref{a1}, we get Eqs. \eqref{mu=0} and  \eqref{mu=2}.

(iii) Firstly, we prove $a_n(\mu)$  is a polynomial  of degree $n$ while $(n\geq 0)$ by mathematical induction.

% By Eq.\eqref{a1}, obviously, $a_n(\mu)$  is a polynomial in $\mu$ of degree $n$ while $(n\geq 0)$.

 {%\color{red}

%By Eq.\eqref{a1}, we get
%\begin{displaymath}
%a_{1}(\mu)=-2\mu+2,\quad a_{2}(\mu)=4\mu^2-8\mu+3.
%\end{displaymath}
%Thus, $a_{1}^{\left ( 0 \right )}=2, a_{2}^{\left ( 0 \right )}=3$. By the recursive relations of $a_{n}$, we have
%\begin{equation}\label{n=00}
%a_{n}^{\left ( 0 \right )}=2a_{n-1}^{\left ( 0 \right )}-a_{n-2}^{\left ( 0 \right )}.
%\end{equation}
%Solving the general term of Eq \eqref{n=00} yields Eq \eqref{n=0}.

Step $1$: for $n=0$ and $n=1$,
\begin{equation}\label{a0=1}
a_{0}(\mu)=1,
\end{equation}
\begin{equation}\label{a1=beta}
a_{1}(\mu)=-2\mu+2.
\end{equation}
It is  obvious that $a_0(\mu)$ and $a_1(\mu)$ are  polynomials of degree $0$ and $1$ respectively.

 Step $2$: for any $n\geq 2$, we assume $a_{n-1}(\mu)$ and $a_{n-2}(\mu)$ are  polynomials of degree $n-1$ and $n-2$ respectively.
% Replacing $n$ with $1$ and $2$ respectively in Eq. \eqref{a1}, we get

\begin{equation}\label{a1=beta1}
a_{2}(\mu)=4\mu^2-8\mu+3.
\end{equation}
%From the recursion relationship and Eq. \eqref{a1=beta} ,
 Note that $a_{n}(\mu)=2(1-\mu)a_{n-1}(\mu)-a_{n-2}(\mu)$. Then, the degree of polynomial $a_n(\mu)$ is just $1$ plus the degree of polynomial $a_{n-1}(\mu)$. Therefore
 $a_n(\mu)$  is a polynomial  of degree $n$, which can be expressed as
 \begin{equation}\label{P-an}
a_{n}(\mu)=\sum_{i=0}^na_{n}^{(i)}\mu^i,
 \end{equation}
 where $a_{n}^{(i)}$ is the coefficient of $\mu^i$.

Replacing  $a_n(\mu)$, $a_{n-1}(\mu)$ and  $a_{n-2}(\mu)$ from Eq. \eqref{P-an} in recursive equation $a_{n}(\mu)=2(1-\mu)a_{n-1}(\mu)-a_{n-2}(\mu)$, we find
$\{a_{n}^{\left ( 0 \right )}\}$ satisfies the following recurrence relation
%Thus, by Eqs. \eqref{a1=beta}-\eqref{a1=beta1}, we have $a_{1}^{\left ( 0 \right )}=2, a_{2}^{\left ( 0 \right )}=3$. And by the recursive relations of $a_{n}$, we have
\begin{equation}\label{n=00}
a_{n}^{\left ( 0 \right )}=2a_{n-1}^{\left ( 0 \right )}-a_{n-2}^{\left ( 0 \right )},
\end{equation}
with initial conditions $a_{0}^{\left ( 0 \right )}=1$ and $a_{1}^{\left ( 0 \right )}=2$.

By solving the recursive equation as shown in Eq. \eqref{n=00}, we obtain  Eq. \eqref{n=0}.

Similarly, we find $\{a_{n}^{\left ( 1 \right )}\}$ satisfies the following recurrence relation
   \begin{equation}\label{n=01}
         a_{n}^{\left ( 1 \right )}=-2a_{n-1}^{\left ( 0 \right )}+2a_{n-1}^{\left ( 1 \right )}-a_{n-2}^{\left ( 1 \right )}=-2n+2a_{n-1}^{\left ( 1 \right )}-a_{n-2}^{\left ( 1 \right )},
   \end{equation}
   with initial conditions $a_{1}^{\left ( 1 \right )}=-2, a_{2}^{\left ( 1 \right )}=-8$;
and $\{a_{n}^{\left ( n \right )}\}$ satisfies the following recurrence relation
         \begin{equation}\label{n=0n}
        a_{n}^{\left ( n \right )}=-2a_{n-1}^{\left ( n-1 \right )},
        \end{equation}
     with initial condition $a_{0}^{\left ( 0 \right )}=1$.
%Solving the general term of Eq \eqref{n=00} yields Eq. \eqref{n=0}.
}
By solving the recurrence relation as shown in Eqs. \eqref{n=01} and \eqref{n=0n}, we obtain  Eqs. \eqref{n=1} and  \eqref{n=n}.
%Similarly, $a_{1}^{\left ( 1 \right )}=-2, a_{2}^{\left ( 1 \right )}=-8$. By the recursive relations of $a_{n}$, we have
%        \begin{equation}\label{n=01}
%         a_{n}^{\left ( 1 \right )}=-2a_{n-1}^{\left ( 0 \right )}+2a_{n-1}^{\left ( 1 \right )}-a_{n-2}^{\left ( 1 \right )}=-2n+2a_{n-1}^{\left ( 1 \right )}-a_{n-2}^{\left ( 1 \right )}.
%        \end{equation}
%Solving the general term of Eq \eqref{n=01} yields Eq \eqref{n=1}.
%
%       $a_{1}^{\left ( 1 \right )}=-2, a_{2}^{\left ( 2 \right )}=4$. By the recursive relations of $a_{n}$, we have
%        \begin{equation}\label{n=0n}
%        a_{n}^{\left ( n \right )}=-2a_{n-1}^{\left ( n-1 \right )}.
%        \end{equation}
%Solving the general term of Eq \eqref{n=0n} yields Eq \eqref{n=n}.

(iv) Here, we prove Eq.~\eqref{a3} by mathematical induction and   Eq.~\eqref{a4} can also be obtained similarly.

  Step $1$:  Eq.\eqref{a3} is true for $n=2$ and $n=3$.
$$ a_{2}=\beta a_{1}-a_{0}=a_{1}^{2}-a_{0}^{2}=(a_{1}+a_{0})(a_{1}-a_{0}),$$
and
$$ a_{3}=\beta a_{2}-a_{1}=\left (  a_{2}-a_{0}\right )a_{1}.  $$

 Step $2$: assume Eq.\eqref{a3} is true for $n=2t$ and $n=2t+1$ ($t$ is an arbitrary nonnegative integer), i.e.,
\begin{displaymath}
a_{2t}= \left ( a_{t}-a_{t-1} \right )\left ( a_{t}+a_{t-1} \right ),
\end{displaymath}
\begin{displaymath}
a_{2t+1}= \left ( a_{t+1}-a_{t-1} \right )a_{t}.
\end{displaymath}

%Now, prove they are true for $'2t+2'$ and $'2t+3'$
Then,
\begin{align*}
a_{2t+2}&=\beta a_{2t+1}-a_{2t}\\
&=\beta \left ( a_{t+1}-a_{t-1} \right )a_{t}-\left ( a_{t}^{2}-a_{t-1}^{2} \right )\\
&=\beta a_{t+1}a_{t}-\beta a_{t-1}a_{t}- a_{t}^{2}+a_{t-1}^{2} \\
&=\beta a_{t+1}a_{t}- a_{t}^{2}-a_{t-1}\left ( \beta a_{t}-a_{t-1} \right )\\
&=\beta a_{t+1}a_{t}- a_{t}^{2}-a_{t-1}a_{t+1}\\
&=a_{t+1}\left ( \beta a_{t}-a_{t-1} \right )- a_{t}^{2}\\
&=a_{t+1}^{2}- a_{t}^{2}\\
&=(a_{t+1}+a_{t})(a_{t+1}- a_{t}),
\end{align*}
and\\
\begin{align*}
a_{2t+3}&=\beta a_{2t+2}-a_{2t+1}\\
&=\beta \left ( a_{t+1}^{2}-a_{t}^{2} \right )-\left ( a_{t+1}-a_{t-1} \right )a_{t}\\
&=\beta a_{t+1}^{2}-\beta a_{t}^{2}- a_{t+1}a_{t}+a_{t-1}a_{t} \\
&=a_{t+1}\left ( \beta a_{t+1}-a_{t} \right )-a_{t}\left ( \beta a_{t}-a_{t-1} \right )\\
&=a_{t+1}a_{t+2}- a_{t}a_{t+1}\\
&=a_{t+1}\left (a_{t+2}- a_{t}  \right ).
\end{align*}
Therefore Eq.\eqref{a3} is true for $n=2t+2$ and $n=2t+3$.
Thus Eq.\eqref{a3} is true for any  $n\geq 2$.\\

Similarly, we can also prove  Eq.\eqref{a4} by  mathematical induction.

\section{Proof of Lemma \ref{Pro_an1}}
\label{Prov_Pro_an1}

(i) Let $n$ $(n\geq 3)$ be an arbitrary odd number and $\mu$ be an arbitrary  root of equation $a_{\frac{n-1}{2}}(\mu)-a_{\frac{n-1}{2}-1}(\mu)=0$.
Therefore,
% According to $a_{n}$'s recursive relationship, we have\\=a_{\frac{n-1}{2}-1}\left ( \mu \right )
        \begin{equation}\label{b2}
        a_{\frac{n-1}{2}-1}\left ( \mu \right )=a_{\frac{n-1}{2}}\left ( \mu \right )=2\left ( 1- \mu\right )a_{\frac{n-1}{2}-1}\left ( \mu \right )-a_{\frac{n-1}{2}-2}\left ( \mu \right ).
        \end{equation}
   Then%     Simplify Eq.\eqref{b2}, we have\\
        \begin{equation}\label{b3}
        \left [ 2\left ( 1- \mu\right )-1 \right ]a_{\frac{n-1}{2}-1}\left ( \mu \right )=a_{\frac{n-1}{2}-2}\left ( \mu \right ).
        \end{equation}
Replacing $a_{\frac{n-1}{2}+1}\left ( \mu \right )$ and $a_{\frac{n-1}{2}}\left ( \mu \right )$ with $2\left ( 1- \mu\right )a_{\frac{n-1}{2}}\left ( \mu \right )-a_{\frac{n-1}{2}-1}\left ( \mu \right )$ and $2\left ( 1- \mu\right )a_{\frac{n-1}{2}-1}\left ( \mu \right )-a_{\frac{n-1}{2}-2}\left ( \mu \right )$  in $1-\frac{ a_{\frac{n-1}{2}+1}-a_{\frac{n-1}{2}-1} }{ a_{\frac{n-1}{2}}-a_{\frac{n-1}{2}-2} }$ respectively, we obtain %and combining Eq.\eqref{b3},
\begin{align}\label{b4}
1-\frac{ a_{\frac{n-1}{2}+1}\left ( \mu\right )-a_{\frac{n-1}{2}-1}\left ( \mu \right ) }{ a_{\frac{n-1}{2}}\left ( \mu \right )-a_{\frac{n-1}{2}-2}\left ( \mu \right ) }
&=1-\frac{ 2\left ( 1- \mu\right )a_{\frac{n-1}{2}}\left ( \mu \right )-a_{\frac{n-1}{2}-1}\left ( \mu \right ) -a_{\frac{n-1}{2}-1}\left ( \mu \right ) }{ 2\left ( 1- \mu\right )a_{\frac{n-1}{2}-1}\left ( \mu \right )-a_{\frac{n-1}{2}-2}\left ( \mu \right )-a_{\frac{n-1}{2}-2}\left ( \mu \right ) }\notag\\
&=1-\frac{\left [   2\left ( 1- \mu\right )-2\right ]a_{\frac{n-1}{2}-1}\left ( \mu \right ) }{ \left \{ 2\left ( 1- \mu\right )-2\left [ 2\left ( 1- \mu\right )-1 \right ] \right \}a_{\frac{n-1}{2}-1}\left ( \mu \right ) }\notag\\
&=2.
\end{align}
One can easy obtain from Eq.\eqref{a3} that% and , we have\\ According to
        \begin{equation}\label{b0}
        a_{n-1}=\left ( a_{\frac{n-1}{2}}-a_{\frac{n-1}{2}-1}  \right )\left ( a_{\frac{n-1}{2}}+a_{\frac{n-1}{2}-1} \right )=0.
        \end{equation}
By using Eqs.\eqref{a3} and \eqref{a4}, we have\\
\begin{align}\label{b08}
a_{n-2}\left ( \mu \right )
&=a_{\frac{n-2-1}{2}}\left ( \mu\right )\left ( a_{\frac{n-2-1}{2}+1}\left ( \mu \right )-a_{\frac{n-2-1}{2}-1}\left ( \mu \right )  \right )\notag\\
&=a_{\frac{n-1}{2}-1}\left ( \mu \right )\left ( a_{\frac{n-1}{2}}\left ( \mu \right )-a_{\frac{n-1}{2}-2}\left ( \mu \right )  \right )\notag\\
&=a_{\frac{n-1}{2}}\left ( \mu \right )\left ( a_{\frac{n-1}{2}}\left ( \mu \right )-a_{\frac{n-1}{2}-2}\left ( \mu \right )  \right )\notag\\
&=1+a_{n-1}\left ( \mu \right )\notag\\
&=1.
\end{align}

(ii) Let $\mu$ be an arbitrary  root of equation $a_{\frac{n-1}{2}}(\mu)+a_{\frac{n-1}{2}-1}(\mu)=0$. Then
%We now proceed as in the proof of (i). According to $a_{n}$'s recursive relationship,We have
        \begin{equation}\label{b6}
        -a_{\frac{n-1}{2}-1}\left ( \mu \right )=a_{\frac{n-1}{2}}\left ( \mu \right )=2\left ( 1- \mu\right )a_{\frac{n-1}{2}-1}\left ( \mu \right )-a_{\frac{n-1}{2}-2}\left ( \mu \right )
        \end{equation}
 Thus%Simplify Eq.\eqref{b6}, we have\\
        \begin{equation}\label{b7}
        \left [ 2\left ( 1- \mu\right )+1 \right ]a_{\frac{n-1}{2}-1}\left ( \mu \right )=a_{\frac{n-1}{2}-2}\left ( \mu \right ).
        \end{equation}
Replacing $a_{\frac{n-1}{2}+1}\left ( \mu \right )$ and $a_{\frac{n-1}{2}}\left ( \mu \right )$ with $2\left ( 1- \mu\right )a_{\frac{n-1}{2}}\left ( \mu \right )-a_{\frac{n-1}{2}-1}\left ( \mu \right )$ and $2\left ( 1- \mu\right )a_{\frac{n-1}{2}-1}\left ( \mu \right )-a_{\frac{n-1}{2}-2}\left ( \mu \right )$   in $1-\frac{ a_{\frac{n-1}{2}+1}-a_{\frac{n-1}{2}-1} }{ a_{\frac{n-1}{2}}-a_{\frac{n-1}{2}-2} }$, we obtain%\\and combining Eq.\eqref{b7},
\begin{align}\label{b8}
1-\frac{ a_{\frac{n-1}{2}+1}\left ( \mu \right )-a_{\frac{n-1}{2}-1}\left ( \mu \right ) }{ a_{\frac{n-1}{2}}\left ( \mu \right )-a_{\frac{n-1}{2}-2}\left ( \mu \right ) }
&=1-\frac{ 2\left ( 1- \mu\right )a_{\frac{n-1}{2}}\left ( \mu \right )-a_{\frac{n-1}{2}-1}\left ( \mu \right ) -a_{\frac{n-1}{2}-1}\left ( \mu \right ) }{ 2\left ( 1- \mu\right )a_{\frac{n-1}{2}-1}\left ( \mu \right )-a_{\frac{n-1}{2}-2}\left ( \mu \right )-a_{\frac{n-1}{2}-2}\left ( \mu \right ) }\notag\\
&=1-\frac{-\left [   2\left ( 1- \mu\right )+2\right ]a_{\frac{n-1}{2}-1}\left ( \mu \right ) }{ \left \{ 2\left ( 1- \mu\right )-2\left [ 2\left ( 1- \mu\right )+1 \right ] \right \}a_{\frac{n-1}{2}-1}\left ( \mu \right ) }\notag\\
&=0
\end{align}
Note that $ a_{n-1}=\left ( a_{\frac{n-1}{2}}-a_{\frac{n-1}{2}-1}  \right )\left ( a_{\frac{n-1}{2}}+a_{\frac{n-1}{2}-1} \right )=0$. By using  Eqs.\eqref{a3} and \eqref{a4}, we have%\\and $n$ is odd number.
\begin{align}\label{b09}
a_{n-2}\left ( \mu \right )
&=a_{\frac{n-2-1}{2}}\left ( \mu \right )\left ( a_{\frac{n-2-1}{2}+1}\left ( \mu \right )-a_{\frac{n-2-1}{2}-1}\left ( \mu \right )  \right )\notag\\
&=a_{\frac{n-1}{2}-1}\left ( \mu \right )\left ( a_{\frac{n-1}{2}}\left ( \mu \right )-a_{\frac{n-1}{2}-2}\left ( \mu \right )  \right )\notag\\
&=-a_{\frac{n-1}{2}}\left ( \mu \right )\left ( a_{\frac{n-1}{2}}\left ( \mu \right )-a_{\frac{n-1}{2}-2}\left ( \mu \right )  \right )\notag\\
&=-\left (  1+a_{n-1}\left ( \mu \right ) \right )\notag\\
&=-1
\end{align}

(iii)  Let $\mu$ be an arbitrary  root of equation $a_{\frac{n-1}{2}}(\mu)=0$. Then
%Let $\acute\mu_{i},i=1,2,\cdots,\frac{n-1}{2}$ denote the solutions of equation $a_{\frac{n-1}{2}}-a_{\frac{n-1}{2}-2} =0$.\\
%
%Let $\hat\mu_{i},i=1,2,\cdots,\frac{n-1}{2}$ denote the solutions of equation $a_{\frac{n-1}{2}}=0$,namely\\
%\begin{equation}\label{b12}
%a_{\frac{n-1}{2}}\left ( \hat\mu_{i} \right )=0,i=1,2,\cdots,\frac{n-1}{2}.
%\end{equation}
%According to $a_{n}$'s recursive relationship, we have\\
        \begin{equation}\label{bb6}
        0=a_{\frac{n-1}{2}}\left ( \mu \right )=2\left ( 1- \mu\right )a_{\frac{n-1}{2}-1}\left ( \mu \right )-a_{\frac{n-1}{2}-2}\left ( \mu \right ),
        \end{equation}
 which yields %Thus, %Simplify Eq.\eqref{bb6}, we have\\
 \begin{displaymath}
 2\left ( 1- \mu\right )a_{\frac{n-1}{2}-1}\left ( \mu \right )=a_{\frac{n-1}{2}-2}\left ( \mu \right ).
 \end{displaymath}
Replacing $a_{\frac{n-1}{2}+1}\left ( \mu \right )$ and $a_{\frac{n-1}{2}-2}\left ( \mu \right )$ with $2\left ( 1- \mu\right )a_{\frac{n-1}{2}}\left ( \mu \right )-a_{\frac{n-1}{2}-1}\left ( \mu \right )$ and $2\left ( 1- \mu\right )a_{\frac{n-1}{2}-1}\left ( \mu \right )$  in $1-\frac{ a_{\frac{n-1}{2}+1}-a_{\frac{n-1}{2}-1} }{ a_{\frac{n-1}{2}}-a_{\frac{n-1}{2}-2} }$ respectively, we obtain\\
\begin{align}\label{bu1}
1-\frac{ a_{\frac{n-1}{2}+1}\left ( \mu \right )-a_{\frac{n-1}{2}-1}\left ( \mu \right ) }{ a_{\frac{n-1}{2}}\left ( \mu \right )-a_{\frac{n-1}{2}-2}\left ( \mu \right ) }
&=1-\frac{ 2\left ( 1- \mu \right )a_{\frac{n-1}{2}}\left ( \mu \right )-a_{\frac{n-1}{2}-1}\left ( \mu \right ) -a_{\frac{n-1}{2}-1}\left ( \mu \right ) }{ -2\left ( 1- \mu \right )a_{\frac{n-1}{2}-1}\left ( \mu \right ) }\notag\\
&=1-\frac{1}{1- \mu }.
\end{align}
Since $\mu \in \left ( 0,2 \right ) $, we have $1-\frac{1}{1- \mu}\in \left ( -\infty ,0 \right )\cup \left ( 2,+\infty  \right )$.

%According to $a_{n}$'s recursive relationship and $n$ is odd number, and
It is easy to verify from Eq.\eqref{a3} that%, we have
%\begin{align}\label{b15}
$$a_{n}\left ( \mu  \right )%&=2\left ( 1- \mu  \right )a_{n-1}\left ( \mu  \right )-a_{n-2}\left ( \mu  \right )\notag\\
=\left (  a_{\frac{n-1}{2}+1}\left ( \mu  \right )-a_{\frac{n-1}{2}-1}\left ( \mu  \right ) \right )a_{\frac{n-1}{2}}\left ( \mu  \right )=0.$$
%\end{align}
%Replacing $a_{n}\left ( \mu  \right )$ with $0$ in $1+a_{n}\left ( \mu  \right )$, and according to Eq.\eqref{a4} and $n$ is odd number, we have
Thus
\begin{align}\label{b16}
1&=1+a_{n}\left ( \mu  \right )\notag\\
&=\left ( a_{\frac{n+1}{2}-1}\left ( \mu  \right )-a_{\frac{n+1}{2}-2}\left ( \mu \right ) \right )\left ( a_{\frac{n+1}{2}}\left ( \mu  \right )+a_{\frac{n+1}{2}-1}\left ( \mu \right ) \right )\notag\\
&=\left ( a_{\frac{n-1}{2}}\left (\mu  \right )-a_{\frac{n-1}{2}-1}\left ( \mu  \right ) \right ) \left [  2\left ( 1- \mu  \right )a_{\frac{n-1}{2}}\left ( \mu  \right )-a_{\frac{n-1}{2}-1}\left ( \mu  \right )+a_{\frac{n-1}{2}}\left ( \mu  \right )  \right ]\notag\\
&=\left ( a_{\frac{n-1}{2}-1}\left ( \mu  \right ) \right )^{2} .
\end{align}

Therefore%According to Eq.\eqref{b0} and combining Eqs \eqref{b16}, we have
\begin{align}%\label{b17}
a_{n-1}\left ( \mu  \right )
&=\left ( a_{\frac{n-1}{2}}\left ( \mu  \right )-a_{\frac{n-1}{2}-1}\left ( \mu  \right ) \right )\left ( a_{\frac{n-1}{2}}\left ( \mu  \right )+a_{\frac{n-1}{2}-1}\left ( \mu  \right ) \right )\notag\\
&=-\left ( a_{\frac{n-1}{2}-1}\left ( \mu  \right ) \right )^{2}\notag\\
&=-1.
\end{align}

Replacing $a_{n}\left ( \mu  \right )$  and $a_{n-1}\left ( \mu  \right )$ with $0$  and $-1$ in $a_{n}\left ( \mu  \right )=2(1-\mu)a_{n-1}\left ( \mu  \right )-a_{n-2}\left ( \mu  \right )$,  we have
\begin{equation}%\label{b18}
a_{n-2}\left ( \mu  \right )=-2\left ( 1- \mu  \right ).
\end{equation}
%\end{proof}

\section{Proof of Lemma \ref{lambda_0_2_2}}
\label{Prov_Pro_an2}

%We shall adopt the same procedure as in the proof of Lemma \ref{Pro_an1}.
%\begin{lemma} \label{lambda_0_2_2}
%$n$ is even number.
%Let $\tilde\mu_{i},i=1,2,\cdots,\frac{n}{2}$ denote the solutions of equation $a_{\frac{n}{2}}+a_{\frac{n}{2}-1}=0$,by substituting $\tilde\mu_{i},i=1,2,\cdots,\frac{n}{2}$ into $1-\frac{ a_{\frac{n}{2}+1}-a_{\frac{n}{2}-1} }{ a_{\frac{n}{2}-1}-a_{\frac{n}{2}-2} }$, $a_{n-1}$ and $a_{n-2}$,respectively,we have
%\begin{equation}\label{even5}
%1-\frac{ a_{\frac{n}{2}}\left ( \tilde\mu_{i} \right )-a_{\frac{n}{2}-1}\left (\tilde\mu_{i} \right ) }{ a_{\frac{n}{2}-1}\left ( \tilde\mu_{i} \right )-a_{\frac{n}{2}-2}\left ( \tilde\mu_{i} \right ) }=1-\frac{1}{1- \tilde\mu_{i}}\in \left ( -\infty ,0 \right )\cup \left ( 2,+\infty  \right )
%\end{equation}
%\begin{equation}\label{even6}
%a_{n-1}\left ( \tilde\mu_{i}  \right )-1
%\end{equation}
%\begin{equation}\label{even7}
%a_{n-2}\left ( \tilde\mu_{i}  \right )=-2\left ( 1- \tilde\mu_{i}  \right )
%\end{equation}
%\end{lemma}
%\begin{proof}
%If $n$ $(n\geq2)$ be an  arbitrary even number, and $\{a_n(\mu)\}_{n\geq 0}$ is a series defined in Lemma \ref{Pro_Rec_an}, we have the following results.
(i)   Let $n$ $(n\geq2)$ be an  arbitrary even number and $\mu$ be an arbitrary  root of equation $a_{\frac{n}{2}}(\mu)-a_{\frac{n-2}{2}}(\mu)=0$. Then $a_{\frac{n}{2}}(\mu)=a_{\frac{n-2}{2}}(\mu)=0$ and
%According to Eq.\eqref{a3} and $n$ is even number,we have\\
        \begin{equation}\label{c0}
        a_{n-1}(\mu)=\left ( a_{\frac{n}{2}}(\mu)-a_{\frac{n}{2}-2}(\mu)  \right )a_{\frac{n}{2}-1}(\mu)=0.
        \end{equation}

 Replacing $a_{\frac{n}{2}-2}\left ( \mu  \right )$ with $a_{\frac{n}{2}}\left ( \mu  \right )$ in  $1-\frac{ a_{\frac{n}{2}+1}-a_{\frac{n}{2}-1} }{ a_{\frac{n}{2}-1}-a_{\frac{n}{2}-2} }$ , we obtain
\begin{equation}\label{c4}
1-\frac{ a_{\frac{n}{2}}\left (\mu \right )-a_{\frac{n}{2}-1}\left ( \mu \right ) }{ a_{\frac{n}{2}-1}\left ( \mu \right )-a_{\frac{n}{2}-2}\left ( \mu \right ) }
=1-\frac{ a_{\frac{n}{2}}\left (\mu \right )-a_{\frac{n}{2}-1}\left ( \mu \right ) }{ a_{\frac{n}{2}-1}\left ( \mu \right )-a_{\frac{n}{2}}\left ( \mu \right ) }=2.
\end{equation}

By using  Eqs.\eqref{a3} and \eqref{a4}, we have %and $n$ is even number. And note that $a_{n-1}\left ( \mu \right )=0$
\begin{align}\label{c08}
a_{n-2}\left (\mu \right )
&=\left ( a_{\frac{n}{2}-1}\left ( \mu \right )-a_{\frac{n}{2}-2}\left ( \mu \right )  \right )\left ( a_{\frac{n}{2}-1}\left ( \mu \right )+a_{\frac{n}{2}-2}\left ( \mu \right )  \right )\notag\\
&=\left ( a_{\frac{n}{2}-1}\left ( \mu \right )-a_{\frac{n}{2}-2}\left ( \mu \right )  \right )\left ( a_{\frac{n}{2}-1}\left ( \mu \right )+a_{\frac{n}{2}}\left ( \mu \right )  \right )\notag\\
&=1+a_{n-1}\left ( \mu \right )\notag\\
&=1.
\end{align}

(ii) Let $\mu$ be an arbitrary  root of equation $a_{\frac{n}{2}-1}(\mu)=0$.
Replacing $a_{\frac{n}{2}}\left ( \mu \right )$ with $2\left ( 1- \mu \right )a_{\frac{n}{2}-1}\left ( \mu \right )-a_{\frac{n}{2}-2}\left ( \mu \right )$ in $1-\frac{ a_{\frac{n}{2}}-a_{\frac{n}{2}-1} }{ a_{\frac{n}{2}-1}-a_{\frac{n}{2}-2} }$, we obtain %and note that $a_{\frac{n}{2}-1}\left (\mu  \right )=0$,
\begin{align}\label{c8}
&1-\frac{ a_{\frac{n}{2}}\left ( \mu \right )-a_{\frac{n}{2}-1}\left (\mu \right ) }{ a_{\frac{n}{2}-1}\left ( \mu \right )-a_{\frac{n}{2}-2}\left ( \mu \right ) }\notag\\
=&1-\frac{ 2\left ( 1-\mu \right )a_{\frac{n}{2}-1}\left ( \mu \right )-a_{\frac{n}{2}-2}\left ( \mu \right )-a_{\frac{n}{2}-1}\left (\mu \right ) }{ a_{\frac{n}{2}-1}\left ( \mu \right )-a_{\frac{n}{2}-2}\left ( \mu \right ) }\notag\\
=&1-\frac{ -a_{\frac{n}{2}-2}\left (\mu \right ) }{ -a_{\frac{n}{2}-2}\left ( \mu \right ) }\notag\\
=&0.
\end{align}
By using  Eqs.\eqref{a3} and \eqref{a4} and noticing that $a_{n-1}\left ( \mu \right )=0$, we have
%According to Eqs.\eqref{a3} and \eqref{a4} and $n$ is even number. And note that
\begin{align}\label{c09}
a_{n-2}\left (\mu \right )
&=\left ( a_{\frac{n}{2}-1}\left ( \mu \right )-a_{\frac{n}{2}-2}\left ( \mu \right )  \right )\left ( a_{\frac{n}{2}-1}\left ( \mu \right )+a_{\frac{n}{2}-2}\left ( \mu \right )  \right )\notag\\
&=-\left ( a_{\frac{n}{2}-1}\left ( \mu \right )-a_{\frac{n}{2}-2}\left ( \mu \right )  \right )\left ( a_{\frac{n}{2}-1}\left ( \mu \right )+a_{\frac{n}{2}}\left ( \mu \right )  \right )\notag\\
&=-\left (  1+a_{n-1}\left ( \mu \right ) \right )\notag\\
&=-1.
\end{align}

(iii)  Let $\mu$ be an arbitrary  root of equation of equation $a_{\frac{n}{2}}(\mu)+a_{\frac{n}{2}-1}(\mu)=0$.
We have %According to $a_{n}$'s recursive relationship,
\begin{equation}\label{c12}
-a_{\frac{n}{2}-1}\left (\mu \right )=a_{\frac{n}{2}}\left (\mu \right )=2(1-\mu)a_{\frac{n}{2}-1}\left (\mu \right )-a_{\frac{n}{2}-2}\left (\mu \right ),
\end{equation}
which yields %Simplify Eq. \eqref{c12} ,we have
\begin{displaymath}
a_{\frac{n}{2}-2}\left (\mu \right )=[2(1-\mu)+1]a_{\frac{n}{2}-1}\left (\mu \right ).
\end{displaymath}
Replacing $a_{\frac{n}{2}}\left ( \mu \right )$ and $a_{\frac{n}{2}-2}\left (\mu \right )$ with $a_{\frac{n}{2}-1}\left (\mu \right )$ and $[2(1-\mu)+1]a_{\frac{n}{2}-1}\left (\mu \right )$ in $1-\frac{ a_{\frac{n}{2}}-a_{\frac{n}{2}-1} }{ a_{\frac{n}{2}-1}-a_{\frac{n}{2}-2} }$ respectively, we have
\begin{align}\label{cu1}
1-\frac{ a_{\frac{n}{2}}\left ( \mu \right )-a_{\frac{n}{2}-1}\left (\mu \right ) }{ a_{\frac{n}{2}-1}\left ( \mu \right )-a_{\frac{n}{2}-2}\left ( \mu \right ) }
&=1-\frac{ -2a_{\frac{n}{2}-1}\left ( \mu \right ) }{ a_{\frac{n}{2}-1}\left ( \mu \right )-\left [ 2\left ( 1-\mu \right )+1 \right ]a_{\frac{n}{2}-1}\left ( \mu \right ) }\notag\\
&=1-\frac{1}{1- \mu}.
\end{align}
Since $\mu \in \left ( 0,2 \right )$, we have $1-\frac{1}{1- \mu}\in \left ( -\infty ,0 \right )\cup \left ( 2,+\infty  \right )$.

%By Eq.\eqref{a3} and $n$ is even number, and note that $a_{\frac{n}{2}}\left ( \mu  \right )+a_{\frac{n}{2}-1}\left ( \mu  \right )=0$, we have
It is easy to verify
%\begin{align}\label{c15}
$$a_{n}\left ( \mu  \right )
=\left [    a_{\frac{n}{2}}\left ( \mu  \right )-a_{\frac{n}{2}-1}\left ( \mu  \right ) \right ]\left [    a_{\frac{n}{2}}\left ( \mu  \right )+a_{\frac{n}{2}-1}\left ( \mu  \right ) \right ]=0.$$
%\end{align}
Then %By Eq \eqref{a4} and $n$ is even number, note that $a_{\frac{n}{2}}\left ( \mu  \right )=-a_{\frac{n}{2}-1}\left ( \mu  \right )$ and $a_{n}\left ( \mu  \right )=0$, we have
\begin{align}\label{c17}
a_{n-1}\left ( \mu   \right )
&=\left ( a_{\frac{n}{2}}\left ( \mu   \right )-a_{\frac{n}{2}-2}\left ( \mu   \right ) \right )a_{\frac{n}{2}-1}\left ( \mu   \right )\notag\\
&=-\left ( a_{\frac{n}{2}}\left ( \mu   \right )-a_{\frac{n}{2}-2}\left ( \mu   \right ) \right )a_{\frac{n}{2}}\left ( \mu   \right )\notag\\
&=-\left ( 1+ a_{n}\left ( \mu   \right )\right )\notag\\
&=-1.
\end{align}

Replacing $a_{n}\left ( \mu  \right )$  and $a_{n-1}\left ( \mu  \right )$ with $0$  and $-1$ in $a_{n}\left ( \mu  \right )=2(1-\mu)a_{n-1}\left ( \mu  \right )-a_{n-2}\left ( \mu  \right )$ respectively,  we have
\begin{equation}\label{c18}
a_{n-2}\left ( \mu  \right )=-2\left ( 1- \mu  \right ).
\end{equation}

The proof is completed.

\end{document}